# Several variable $p$-adic families of Siegel-Hilbert cusp eigensystems and their Galois representations

J.Tilouine[*]   E.Urban[†]

November 10, 1998

## Abstract

**Français:** Soit $F$ un corps totalement réel et $G = GSp(4)_{/F}$. Dans cet article, nous montrons sous une hypothèse faible qu'étant donné un système $\lambda$ de valeurs propres de Hecke $(p, P)$-ordinaire (pour un parabolique $P$ quelconque fixé de $G$), il existe une famille $\underline{\lambda}$ à plusieurs variables de systèmes de valeurs propres de Hecke quasi-ordinaires en $p$ qui le contient. L'hypothèse est que $\lambda$ intervienne dans la cohomologie d'un système de coefficients régulier. Si $F = \mathbf{Q}$ (le nombre de variables est alors 3), nous construisons la famille $p$-adique à trois variables de représentations galoisiennes $\rho_{\underline{\lambda}}$ associée à $\underline{\lambda}$. Enfin, sous des hypothèses géométriques (qui seront satisfaites si l'on montre que les représentation galoisiennes de la famille proviennent de motifs de Grothendieck) nous montrons que $\rho_{\underline{\lambda}}$ est quasi-ordinaire pour le parabolique dual de $P$.

**English:** Let $F$ be a totally real field and $G = GSp(4)_{/F}$. In this paper, we show under a weak assumption that, given a Hecke eigensystem $\lambda$ which is $(p, P)$-ordinary for a fixed parabolic $P$ in $G$, there exists a several variable $p$-adic family $\underline{\lambda}$ of Hecke eigensystems (all of them $(p, P)$-nearly ordinary) which contains $\lambda$. The assumption is that $\lambda$ is cohomological for a regular coefficient system. If $F = \mathbf{Q}$, the number of variables is three. Moreover, in this case, we construct the three variable $p$-adic family $\rho_{\underline{\lambda}}$ of Galois representations associated to $\underline{\lambda}$. Finally, under geometric assumptions (which would be satisfied if one proved that the Galois representations in the family come from Grothendieck motives), we show that $\rho_{\underline{\lambda}}$ is nearly ordinary for the dual parabolic of $P$.

[*]LAGA-UMR 7539, Institut Galilée, Université de Paris-Nord, Av. J.B.Clément, 93430 Villetaneuse, France, et Institut Universitaire de France

[†]CNRS-UMR 7539, Institut Galilée, Université de Paris-Nord, Av. J.B.Clément, 93430 Villetaneuse, France and UCLA, Department of Math. Los Angeles CA 90095 USA

# Contents

# 0 Introduction

For the group of symplectic similitudes $G = GSp_4$ over a totally real field $F$ of degree $d$, we develop Hida Theory in a manner similar to [16] and [17]. This means that we introduce a big $p$-adic cuspidal Hecke algebra $\mathbf{h}$ defined as the inverse limit of commutative algebras of Hecke correspondences acting on the cohomology of Shimura varieties whose $p$-level tends to infinity. In particular, the study focuses on the direct factor $\mathbf{h}^{n.o}$ of this algebra cut out by the so-called $p$-nearly ordinary idempotent. Hida theory seeks to establish

1. the independence of the weight for $\mathbf{h}^{n.o}$

2. the control of $\mathbf{h}^{n.o}$ when localized at "arithmetic" codimension one primes of the relevant Iwasawa algebra $\Lambda$

3. its finiteness over $\Lambda$ and what is its Krull dimension (it is equal to that of $\Lambda$ once one has established the torsion-freeness of $\mathbf{h}^{n.o}$ over $\Lambda$)

4. the existence of associated Galois representations over suitable local components thereof.

Actually a more precise result is desirable for those Galois representations. Namely, their restriction to the decomposition groups at primes above $p$ should take values in some parabolic subgroups determined by the type of near ordinarity imposed to the Hecke algebra (*i.e.* to the corresponding automorphic forms). A Galois representation satisfying this local condition is called nearly ordinary (see [36], Section 5). For $GL_2$ over a totally real field, it is a theorem of Wiles [47], see also Hida [15]) that if a Hilbert modular form is nearly ordinary, then its Galois representation is nearly ordinary. In more general cases, however, this fact is not known. A study of this question has been our main motivation (see Section 7 below). The method for dealing with points 1-3 is cohomological, using the natural faithful representation of $\mathbf{h}$ on the cohomology of a Siegel manifold of level divisible by $p^\infty$. This approach follows the ideas of Hida's recent works concerning the case of $GL(n)$ (see [17] and [18]). A good deal of our intermediate results, like independence and control theorems for the full cohomology and probably for the boundary cohomology, are very general; they (should) hold at least for split connected reductive groups of type A,B,C or D over any totally real field. Our work is motivated by three reasons. The first is that it allows one to speak of Hida families for nearly ordinary cohomological genus two Siegel-Hilbert cusp eigensystems $\lambda$. Thus, it should provide new families of $p$-adic $L$-functions (although several variables $p$-adic Siegel-Hilbert Eisenstein measures should be constructed first) and new Main Conjectures. Second, for $F = \mathbf{Q}$ it provides big Galois representations

$$\rho_\lambda : \text{Gal}(\overline{\mathbf{Q}}/\mathbf{Q}) \to GSp_4(\mathbf{Z}_p[[T_1, T_2, T_3]])$$

associated to a family $\lambda$ (Section 7 below). The specializations of $\rho_\lambda$ at almost all arithmetic primes $\mathbf{P}_\theta$ (see Definition 6.2.4 in the text) are the Galois representations $\rho_{\lambda_\theta}$ associated to the Hecke eigensystem $\lambda_\theta$ defined by $\lambda$ in weight $\theta$. These representations $\rho_{\lambda_\theta}$ have been constructed by Shimura [32],Chai-Faltings [7], R. Taylor [34] and more recently by Laumon [23] and Weissauer [46] independently. We have attempted in [36] Section 10 at formulating a generalization of the Langlands correspondence predicting for these Galois representations a specific behaviour at $p$, namely, $\rho(D_v) \subset P_v$ where $P_v$ is in the conjugacy class of the Langlands dual of the parabolic giving the condition of near ordinarity on the automorphic side. This we can prove, if we assume that the Newton polynomial at $p$ associated to crystalline $\rho_{\lambda_\theta}$'s is given by the $p$-Euler factor of the automorphic $L$ function of $\lambda_\theta$. Our proof then makes use of the observation that in an analytic family of nearly ordinary forms, many have level prime to $p$. Such forms give rise conjecturally to crystalline representations; for those, one can compare the Hodge and Newton polygons and their ordinarity follows. Then we use a density argument to conclude it holds for the family, hence for *all* its members. Actually, for $F \neq \mathbf{Q}$, assuming that 4-dimensional Galois representations associated to the $\lambda_\theta$ exist one can still show $\rho(D_v) \subset P_v$ despite the fact the Newton and the Hodge polygons may never meet (cf. Remark at the end of Appendix B1); indeed, assuming a rather natural assumption, called $t$-separability, we obtain the ordinarity even for $F \neq \mathbf{Q}$. The third motivation is to apply this ordinarity result to the Galois representation associated by R. Taylor to modular forms on $GL(2, K)$ where $K$ is imaginary quadratic, providing in this way a lacking ingredient in a paper by one of the authors (see Conjectures 3.3.1 and 3.3.2 of [39]).

In order to state a striking particular case of the present work, we introduce some notations. Let $V$ be a finite dimensional complex vector space and $\rho : U^*(2) \to GL(V)$ be an irreducible algebraic representation of the maximal compact real Lie subgroup of $Sp(4, \mathbf{R})$; let $\rho_{\mathbf{C}} : GL(2, \mathbf{C}) \to GL(V)$ be its complexification. Let $f$ be a holomorphic $V$-valued Siegel cusp form of weight $\rho$ and level group $\Gamma$: $f((AZ + B)(CZ + D)^{-1}) = \rho_{\mathbf{C}}(CZ + D)f(Z)$ for $\begin{pmatrix} A & B \\ C & D \end{pmatrix} \in \Gamma$; assume $f$ is eigen for all Hecke operators $T_\ell$, $R_\ell$, $S_\ell$ (see Def.6.1) for all rational primes $\ell$ prime to the level $N$ of $\Gamma$. The eigenvalue of $f$ for the Hecke operator $T$ is denoted by $\lambda(T; f)$. Let $\mathcal{O}_0$ be a Dedekind ring finite over $\mathbf{Z}$ containing the eigenvalues of $f$. The highest weight $\chi = (a_1, a_2) \in \mathbf{Z}^2$ of $\rho$ satisfies $a_1 \geq a_2 \geq 3$ (noe that classical Siegel modular forms correspond to $a_1 = a_2 = k \geq 3$). Let $\chi' = (a, b) = (a_1 - 3, a_2 - 3)$; let $B'$ be the standard Borel subgroup of $Sp_4$ (see Section 1.1); put

$$L(\chi'; \mathbf{Z}) = L(a, b; \mathbf{Z}) = \{f : Sp(4, \mathbf{Z}) \to \mathbf{Z}; f \text{ is regular, defined over } \mathbf{Z} \text{ and}$$

$$f(gtu) = {\chi'}^{-1}(t)f(g) \text{ for } t \text{ semisimple in } B' \text{and } u \text{ unipotent in } B'\}$$

viewed as left $\Gamma$-module by $\gamma.f(g) = f(\gamma^{-1}g)$. For any ring $A$, let $L(\chi'; A) = L(\chi'; \mathbf{Z}) \otimes A$. It is known (see [33] p.323) that the eigensystem associated to $f$ occurs in $H^3(\Gamma, L(a, b; \mathbf{C}))$. We assume that $\chi'$ is regular: $a - b > 0$ and $b > 0$. This excludes classical Siegel modular forms (for which $a = b = k - 3$).

Let $p$ be a rational prime, relatively prime to $N$ and to the order of the torsion subgroup of $H^q(\Gamma, L(a, b; \mathbf{Z}))$ ($q = 2, 3$). Let $v$ be a non trivial valuation of $\mathcal{O}_0$ such that $v(p) > 0$; assume that $v(\lambda(T_p; f)) = 0$ and $v(\lambda(R_p; f)) = b$. Let $\mathcal{O}$ be the completion of $\mathcal{O}_0$ at $v$.

Let $\Lambda' = \mathcal{O}[[T_1, T_2]]$. Let $u = 1 + p$. For any pair $(a, b)$ with $a > b > 0$, we define the arithmetic prime $\mathcal{P}_{a,b}$ of $\Lambda'$ as the kernel of the homomorphism $\Lambda \to \mathcal{O}$ given by $T_1 \mapsto u^a - 1$, $T_2 \mapsto u^b - 1$. Let $\Gamma_0(p)$ be the subgroup of $\Gamma$ consisting in matrices whose reduction mod. $p$ fall in $B'$.

**Theorem 0.1** *Under the previous assumptions, the ordinary part $\mathbf{h}^o$ of the cuspidal Hecke $\mathcal{O}$-algebra is finite torsion-free over $\Lambda'$; moreover, there exists a finite flat extension $\mathbf{J}$ of $\Lambda'$ and a $\Lambda'$-algebra homomorphism $\lambda : \mathbf{h}^o \to \mathbf{J}$ such that for any arithmetic prime $\mathcal{P}_{a',b'}$ with $a' \equiv a$, $b' \equiv b \ mod.(p-1)$, $a' \geq b' \geq 0$, and any prime $\mathbf{P}$ in $\mathbf{J}$ above $\mathcal{P}_{a',b'}$, the reduction of $\lambda \ mod.\mathbf{P}$ "corresponds" to a Hecke eigenclass $c_{a',b'}$ in $H^3(\Gamma_0(p), L(a', b'; \mathcal{O}))$. For $(a', b') = (a, b)$, this Hecke eigenclass is deduced from $c_{a,b}$ by the p-stabilization isomorphism:*

$$Res : H^3_{ord}(\Gamma, L(a, b; \mathcal{O})) \cong H^3_{ord}(\Gamma_0(p), L(a, b; \mathcal{O}))$$

*(see Proposition 3.2 of Section 3.5 below).*

*Secondly, assuming that multiplicity one holds for $GSp_4$; then, there exists a finite extension $\mathcal{F}$ of $\mathrm{Frac}(\mathbf{J})$ and a continuous Galois representation $\rho_\lambda : Gal(\overline{\mathbf{Q}}/\mathbf{Q}) \to GSp_4(\mathcal{F})$ associated to the $\Lambda'$-algebra homomorphism $\lambda$; namely, it is unramified outside $Np$, and for* any prime $\ell$ relatively prime to $Np$, the characteristic polynomial of $\rho_\lambda(Frob_\ell)$ is given by $\lambda(Q_\ell(X))$ where

$$Q_\ell(X) = X^4 - T_\ell X^3 + \ell(R_\ell + (\ell^2+1)S_\ell)X^2 - \ell^3 T_\ell S_\ell X + \ell^6 S_\ell^2$$

**Comments:** 1) In the case $\mathbf{J} = \Lambda'$, the $p$-adic analyticity of the family of Hecke eigensystems amounts to the Kummer congruences: if $a' \equiv a''$ and $b' \equiv b''$ mod.$(p-1)p^n$, then the eigensystems $\lambda_{a',b'}$ and $\lambda_{a'',b''}$ are congruent mod.$p^{n+1}$.

2) Roughly speaking, one can say that $\rho_\lambda$ lifts the Galois representation $\rho_f : Gal(\overline{\mathbf{Q}}/\mathbf{Q}) \to GL_4(\mathbf{Q}_p)$ constructed by Laumon [23] and Weissauer [46]. Even without assuming multiplicity one for $GSp_4$, R. Weissauer has associated a Galois representation $\rho_\pi : Gal(\overline{\mathbf{Q}}/\mathbf{Q}) \to GL_4(\overline{\mathbf{Q}}_p)$ to any cohomological cuspidal representation $\pi$ on $GSp_4$ (this result is only written up for $\pi_\infty$ in the holomorphic discrete series). Then, in our theorem, we can remove the multiplicity one assumption and the irreducibility assumption; the conclusion is then only that there exists a representation $\rho_\lambda : Gal(\overline{\mathbf{Q}}/\mathbf{Q}) \to GL_4(\mathcal{F})$ with the correct characteristic polynomials.

3) If one assumes multiplicity one, the main ingredient to show that $\rho_\lambda$ falls in $GSp_4(\mathcal{F})$ is that all smooth irreducible representations $\pi$ of $GSp_4(\mathbf{A})$ are autodual. This comes by proving that for a given local constituent $\pi_v$ of $\pi$, the traces of the Hecke operators acting on the space of $\pi_v$ and that of $\pi_v^\vee \otimes \omega_{\pi_v} \circ \nu$ are the same; this last point results then from the fact there exists an element $g_0 \in GSp_4$ such that ${}^t g^{-1} \nu(g) = g_0 g g_0^{-1}$, here $g_O = J = \begin{pmatrix} 0_2 & 1_2 \\ -1_2 & 0_2 \end{pmatrix}$ works. We are grateful to Clozel for showing us this argument.

4) If multiplicity one holds and if reduction of $\rho_f$ modulo the maximal ideal of $\overline{\mathbf{Z}}_p$ is still absolutely irreducible, then, one can see that $\rho_\lambda$ takes values in $GSp_4(\mathbf{J})$.

5) Moreover, we can prove that the image by $\rho_\lambda$ of a decomposition group at $p$ is contained in $B(\mathbf{J})$ (up to conjugation in $GSp_4(\mathbf{J})$ if the two statements **S1** and **S2** below hold. Consider a Zariski dense set of arithmetic primes $(\mathcal{P}_{a',b'})$ with $(a',b') \equiv (a,b)$ mod.$(p-1)$ and $a' > b' > 0$ and fix a prime $\mathbf{P}$ above $(\mathcal{P}_{a',b'})$ in $\mathbf{J}$. We know from Proposition 3.2 and Comment below Conjecture 7.2, that $\rho_\lambda$ mod.$\mathbf{P}$ is crystalline at $p$.

**S1** For a Zariski dense subset of arithmetic primes $(\mathcal{P}_{a',b'})$ as above, $\rho_\lambda$ mod.$\mathbf{P}$ has four distinct Hodge-Tate weights

**S2** The slopes of the Newton polygon of the $\phi$-filtered module attached to $\rho_\lambda$ mod.$\mathbf{P}$ are equal to the valuations of the eigenvalues of the roots of $\lambda_\mathbf{P}^{(p)}(Q_p^{(p)}(X))$. Here $\lambda_\mathbf{P}^{(p)} = \lambda$ mod.$\mathbf{P}$ is viewed as a character of the Hecke algebra of level group $\Gamma$ (prime to $p$) and $Q_p^{(p)}(X)$ is defined as "$Q_p$" since now $p$ is prime to the level of the group $\Gamma$.

**S1** is implied by the stability of the $L$-packets at infinity; this is known if $\lambda_\mathbf{P}^{(p)}$ is supercuspidal at some finite place and seems accessible to specialists even without this assumption.

**S2** is harder. An evidence towards it seems more accessible, namely that the slopes of the Newton polygon occur among the valuations of the eigenvalues of the roots of $\lambda_{\mathbf{P}}^{(p)}(Q_p(X))$.

6) The continuity of the Galois representation $\rho_\lambda$ means that it preserves a lattice $T$ in $\mathcal{F}^4$ and $\rho_\lambda : Gal(\overline{\mathbf{Q}}/\mathbf{Q}) \to GL(T)$ is continuous for the natural topologies; hence the theorem implies that for a Zariski dense family of arithmetic primes $\mathcal{P}_{a',b'}$ (with $(a',b') \equiv (a,b)$ mod.$p-1$), the reduction of $\rho_\lambda$ modulo $\mathcal{P}_{a',b'}$ exists and is the Galois representation associated to the eigenclass $c_{a',b'}$.

7) In the text below, one can even study a bigger Hecke algebra $\mathbf{h}^{n.o}$ over a bigger Iwasawa algebra, called the Hida-Iwasawa algebra $\Lambda$ (see Section 6, Def.6.2) obtained by adjoining to the "semisimple variables" of $\Lambda'$ the "central variable(s)". The Hida-Iwasawa algebra is isomorphic to the algebra of a finite group over a ring of formal power series in $3 = 2 + 1$ variables if $F = \mathbf{Q}$, and in $2d + 1 + \delta$ variables, if $F$ is totally real of degree $d$. Therefore, $p$-adic families of cuspidal Hecke eigensystems indexed by these variables exist, as well as corresponding families of Galois representations (assuming Conjecture 1 of Section 7).

An important application of our result is currently investigated by one of the authors (see [40]). Namely, one can use this theory to study congruences between families of Eisenstein-Klingen series and families of Siegel cusp eigenforms in order to produce non semisimple Galois representations with values in the Siegel parabolic. This basic idea can be used to study Greenberg's Main Conjecture for the symmetric square of a $p$-ordinary elliptic curve (see [40] and [20]). An important step of this investigation is to prove that the Galois representation $\rho_\lambda$ attached to an ordinary cusp eigensystem $\lambda$ is absolutely irreducible, provided that $\lambda$ is not globally endoscopic and that the residual representation, if reducible, admits an irreducible two-dimensional subquotient which is modular in Serre's sense and satisfies the assumptions of the theorem of Wiles and Taylor ([48] and [35]).

We give now a short survey of the contents of the present paper. In Section 1, we recall the standard Notation and definitions used throughout this article. We give the definitions of local Hecke operators and their action on flags in section 2. In section 3, we prove first a theorem of independence of the weight for the ordinary cohomology of $\Gamma_1$-type congruence subgroups with level $p^\infty$ subgroups (Cor.3.2); then, we establish control theorems for the ordinary cohomology group of "bottom degree" (that is, the first degree of non-vanishing for regular coefficients, namely, by Franke's theorem, the middle degree $3d$). There are two such control theorems: weak (Th. 3.2) and exact (Th.3.3). Up to this point, it appears clearly that the proofs work for any reductive group, after appropriate translation. Then, after a study of the strata cohomology in Section 4, we show in Section 5 control theorems for the boundary cohomology by studying the degeneracy of the nearly ordinary part of the spectral sequences attached to the parabolic subgroups defining the Borel-Serre compactification; there, restriction to genus two is required to insure the degeneracy of the spectral sequences at $E^2$ (th.5.2 for $F = \mathbf{Q}$ and Lemma 5.2 in general). From these control theorems, we deduce our main results (Corollary 5.1 for $F = \mathbf{Q}$ and Theorem 5.8 in general) for the ordinary interior cohomology of congruence subgroups of $\Gamma_1$-type, of level $p^\infty$ in middle degree $q = 3d$. We apply these results in Section 6, to obtain their counterpart for the nearly-ordinary cohomology of the $p^\infty$-ramified Siegel threefolds (Th.6.2 and Th.6.3). Then one deduces a control theorem for the nearly ordinary cuspidal Hecke algebra for $GSp_4$ (Cor.6.3). This result is enough to insure the existence of a $(2d+1+\delta)$-variables Hida family interpolating a Siegel-Hilbert cusp eigensystem of given level and weight, as well as the specialization at any arithmetic primes of such a family into such eigensystems of given level, weight and central character (Cor.6.7). Finally in Section 7, we construct, assuming standard conjectures proven in some cases, the Galois representation attached to a Hida family, and we discuss its near ordinarity (Th.7.1). It will be obvious to the reader that our approach follows closely Hida's recent paper [17] devoted to the $GL(n)$ case. We learned his method in a course at Paris-Nord in March 1994. We are glad to acknowledge our debt to him here. We also benefited from several conversations with R. Weissauer whom we thank for his invaluable explanations. Part of this paper has been written during visits of the first author at the Tata Institute of Fundamental Research in Bombay and at Mannheim University; let these institutions be thanked for their kind hospitality.

# 1 Notation and preliminaries

## 1.1 The symplectic group

Let $J = \begin{pmatrix} 0_2 & 1_2 \\ -1_2 & 0_2 \end{pmatrix} \in GL_4$. Let $G = GSp_4 = \{g \in GL_4; {}^t gJg = \nu(g)J\}$ be the group of symplectic similitudes of $J$; the character $\nu : G \to \mathbf{G}_m$ is called the multiplier of $GSp_4$: The group $G^1 = Sp_4$ of symplectic isometries is the kernel of $\nu$. These groups are smooth group schemes over $\mathbf{Z}$. The symplectic module acted on by $G$ is denoted by $(W, <\ ,\ >)$ and we write its canonical basis (over $\mathbf{Z}$) as $\{e_1, e_2, f_1, f_2\}$ with $<e_1, e_2>\ =\ <f_1, f_2>\ =\ 0$ and $<e_i, f_j>\ =\ \delta_{i,j}$. Consider the standard maximal torus $T$ of $G$ consisting of diagonal matrices. We identify it to $\mathbf{G}_m^3$ by the isomorphism: $\tau : \mathbf{G}_m^3 \cong T$,

$$(t_1, t_2; x) \mapsto diag(t_1, t_2, xt_1^{-1}, xt_2^{-1}).$$

Its character group $X^*(T)$ is then identified to the sublattice of $\mathbf{Z}^3$ consisting of $(a_1, a_2; b)$'s such that $a_1 + a_2 \equiv b$ mod. 2. We write $\chi = (a_1, a_2; b)$ for the character $\chi$ defined by

$$\chi(\tau(t_1, t_2; x)) = t_1^{a_1} t_2^{a_2} x^{(b-a_1-a_2)/2}.$$

We denote by $B$ the standard Borel of $G$ consisting in matrices of the form

$$\begin{pmatrix} * & * & * & * \\  & * & * & * \\  &  & * &  \\  &  & * & * \end{pmatrix}$$

Its unipotent radical is denoted by $B^+$. The roots associated to $(G, B, T)$ are $(t_1/t_2)^{\pm 1}$ and $(x^{-1}t_i t_j)^{\pm 1}$ for $1 \leq i \leq j \leq 2$; we denote their set by $R_\mathbf{Z}$, these are characters defined over $\mathbf{Z}$. The positive roots are the four ones given by the exponent $+1$; their set is $R_\mathbf{Z}^+$; the simple roots are $\alpha_1 = t_1/t_2$ and $\alpha_2 = \nu^{-1}t_2^2$ (the long root). Each conjugacy class in $G$ of parabolic subgroups (over $\mathbf{Z}$) has a unique representative containing $B$; there are exactly three such parabolic subgroups:

- the Siegel parabolic:
$$P = \{\begin{pmatrix} A & B \\ 0_2 & D \end{pmatrix} \in G\};$$

  the conditions on the $2 \times 2$-blocks are that ${}^t AD = \nu.1_2$ for some scalar $\nu$ and $A^{-1}B$ is symmetric. The Levi subgroup is:
$$M = \{g = \begin{pmatrix} A & 0_2 \\ 0_2 & x\ {}^t A^{-1} \end{pmatrix}; A \in GL_2,\ x \in \mathbf{G}_m\}$$

It is therefore isomorphic to $GL_2 \times \mathbf{G}_m$. We write this isomorphism $(A,x) \in GL_2 \times \mathbf{G}_m \mapsto \mu(A;x)$ this isomorphism; it maps the center $\mathbf{G}_m^2$ of $GL_2 \times \mathbf{G}_m$ onto the center $Z_M$ of $M$. Let $M^1$ be the derived group of $M$ and $C_M = M/M^1$ be its cocenter. The natural map $Z_M \to C_M$ is an isogeny of degree 2 which we identify to $(z,x) \mapsto (z^2, x)$. The unipotent radical of $P$ is

$$P^+ = \{ \begin{pmatrix} 1_2 & S \\ 0_2 & 1_2 \end{pmatrix}; {}^t S = S \}$$

- the Klingen parabolic:

$$P^* = \{ \begin{pmatrix} a & * & * & * \\ 0 & \alpha & * & \beta \\ 0 & 0 & b & 0 \\ 0 & \gamma & * & \delta \end{pmatrix}; \alpha\delta - \beta\gamma = ab \} \cap G;$$

The Levi subgroup is:

$$M^* = \{ \begin{pmatrix} a & & & \\ & \alpha & & \beta \\ & & b & \\ & \gamma & & \delta \end{pmatrix}; \alpha\delta - \beta\gamma = ab \neq 0 \};$$

It is therefore isomorphic to $\mathbf{G}_m \times GL_2$. We write

$$\mu^* : (a, \begin{pmatrix} \alpha & \beta \\ \gamma & \delta \end{pmatrix}) \mapsto \begin{pmatrix} a & & & \\ & \alpha & & \beta \\ & & b & \\ & \gamma & & \delta \end{pmatrix}$$

with $ab = \alpha\delta - \beta\gamma$. By this isomorphism, the center $\mathbf{G}_m^2$ of $\mathbf{G}_m \times GL_2$ is mapped to the center $Z_{M^*}$. The natural isogeny from the center to the cocenter of $M^*$ becomes $(x,z) \mapsto (x, z^2)$. The unipotent radical is

$$P^{*+} = \{ \begin{pmatrix} 1 & \alpha & \beta & \gamma \\ 0 & 1 & \gamma & 0 \\ 0 & 0 & 1 & 0 \\ 0 & 0 & -\alpha & 1 \end{pmatrix}; \alpha, \beta, \gamma \text{ arbitrary} \}$$

- the Borel subgroup $B = TB^+$, which has been described already.

Note that $P$ and $P^*$ are maximal and that $B = P \cap P^*$. For $Q = P$ or $P^*$ with Levi subgroup $M_Q$ (often abbreviated as $M$ when the context is clear) and Levi decomposition $Q = MQ^+$, we put $\tilde{Q} = M^1 Q^+$ and $C = M/M^1$. Consider the natural map $T \to C$ induced by the inclusions $B \subset Q$ and $B^+ \subset \tilde{Q}$. It identifies the group $X^*(C)$ of characters of $C$ to the sublattice of $X^\star(T) \subset \mathbf{Z}^3$ consisting of

- the $(a, a; b) \in \mathbf{Z}^3$ such that $2a \equiv b$ mod. 2, if $Q = P$,
- the $(a, 0; b) \in \mathbf{Z}^3$ such that $a \equiv b$ mod. 2, if $Q = P^*$.

Let $Q$ be any parabolic subgroup with Levi subgroup $M$ and unipotent radical $Q^+$. We denote by $\Delta_Q$ the set roots of $(M, T)$ and $\Delta_Q^\pm$ the set of positive resp. negative roots for $(M, B \cap M, T)$; similarly, we write $R_Q^\pm$ the set of roots for $T$ acting on $Q^\pm$. We denote by $W$ the Weyl group of $(G, B, T)$, $W_M$ that of $(M, B \cap M, T)$. The set of positive (resp. negative) roots of $(G, B, T)$ is then decomposed as $R_\mathbf{Z}^\pm = \Delta_Q^\pm \coprod R_Q^\pm$

## 1.2 Local systems over Siegel-Hilbert modular varieties

Let $F$ be a totally real field of degree, say, $d$; let $\mathbf{r}$ be its ring of integers, with discriminant $D$. We note $I_F$ the set of embeddings of $F$ in $\bar{\mathbf{Q}}$. For $f : \text{Spec } \mathbf{r} \to \text{Spec } \mathbf{Z}$, we consider the scheme $H = f_* f^* G$, restriction of scalars of the base change of $G$ to $\mathbf{r}$; for any ring $A$, $H(A) = G(A \otimes \mathbf{r})$; $H$ is a group scheme over $\mathbf{Z}$, smooth over $\mathbf{Z}[1/D]$. The multiplicator $\nu$ induces a natural morphism from $H$ to $f_* f^* \mathbf{G}_m$. Its kernel is denoted by $H^1$; for any ring $A$, we have $H^1(A) = G^1(A \otimes O_F)$. In this paper, we denote by $Q$ a parabolic subgroup of $G$ and $Q' = Q \cap G^1$. We know that there are three conjugation classes of parabolic subgroups in $G$ containing respectively $P$, $P^*$ and $B$. Let us fix $Q \in \{B, P, P^*\}$. We write the Levi decomposition of $Q$ as $MQ^+$. Let $M' = M \cap G^1$, $M^1 = M'^1$ the derived group of $M$ and $\tilde{Q} = M^1 Q^+$.

Let us fix an odd prime $p$; let $S_p = \{v \text{ place of } F; v|p\}$. Let $J$ be a subset of $S_p$.

**Definition 1.1** *A $J$-proper standard parabolic subgroup $Q$ of $H \otimes \mathbf{Z}_p$ (abbreviated as $J$-proper SPS) is a product $Q = \prod_{v|p} \text{Res}_{\mathbf{Z}_p}^{\mathbf{r}_v} Q_v$ of parabolic subgroups $Q_v$ of $G \otimes \mathbf{r}_v$ such that for each $v \in J$, $Q_v \in \{P, P^*, B\}$, and for $v \notin J$, $Q_v = G$. Let $Q' = \prod_{v|p} Q'_v$, $M = \prod_{v|p} M_v$ its Levi subgroup and $Q^+ = \prod_{v|p} Q_v^+$ its unipotent radical system; $Q$ is a semi-direct product $MQ^+$. Let $M^1 = \prod_{v|p} M_v^1$ be the derived group of $M$ and $\tilde{Q} = M^1 Q^+ = \prod_{v|p} \tilde{Q}_v$. For any $\mathbf{Z}_p$-algebra $A$, let $A_v = A \otimes \mathbf{r}_v$ for each $v|p$; then $Q(A) = \prod_{v|p} Q_v(A_v)$ via the canonical isomorphism $H(A) \cong \prod_{v|p} G(A_v)$. Similarly for $Q'(A)$, $Q^+(A)$ and $\tilde{Q}(A)$.*

We fix an arbitrary $J$-proper SPS $Q$ with Levi components $(M, Q^+)$. Most of our results will be valid without further assumption on $Q$, although results of Section 6 will require that $J = S_p$ and $Q$ is of the form $Q = \Pi_H = f_* f^* \Pi$ for $\Pi \in \{B, P, P^*\}$.

Let **A** be the ring of rational adeles and $\mathbf{Q}_f$ resp. $\mathbf{Q}_\infty$ be its finite, resp. infinite part. Let $H_\mathbf{A}$, resp. $H_f$, $H_\infty$ be the group of **A**-points of $H$, resp. of $\mathbf{Q}_f$-, $\mathbf{Q}_\infty$-points. Let $U_\infty$ be the stabilizer in $H_\infty = \prod_{v|\infty} G(F_v)$ of the map $h : \mathbf{C}^\times \longrightarrow H_\infty$ whose $v$-component is

$$h_v(x+iy) = \begin{pmatrix} x1_2 & y1_2 \\ -y1_2 & x1_2 \end{pmatrix} \in G(F_v).$$

Let us fix once for all a compact open subgroup $U$ of $H(\hat{\mathbf{Z}})$ of levelprime to $p$. For any $r \geq 1$ We denote by $U_{\tilde{Q}}(p^r)$ the compact open subgroup of $H_f$ defined as

$$U_{\tilde{Q}}(p^r) = \{g \in U; g \bmod p^r \in \tilde{Q}(\mathbf{Z}/p^r\mathbf{Z})\}$$

**Remark:** Despite the Notation, the level of these subgroups is concentrated on the $J$-part of $S_p$ since the condition is void at places $v \notin J$. If $J = \emptyset$, these groups are of level trivial at $p$.

We form the Shimura varieties:

$$S_r(U) = H_\mathbf{Q} \backslash H_\mathbf{A} / U_{\tilde{Q}}(p^r) U_\infty$$

Their connected components are Siegel-Hilbert modular varieties of dimension $3d$. For $s \geq r$, there is a natural finite morphism $S_s(U) \to S_r(U)$; the varieties $S_r(U)$ together with these transition maps form an inverse system.

Let $\mathcal{O}$ be the valuation ring of a finite extension $K$ of $\mathbf{Q}_p$. Let $L_K$ be a finite dimensional $K$-vector space with a rational action of $H_p = H \otimes \mathbf{Q}_p$; we assume it is "adapted to $Q$", i.e it is (algebraically) induced from $Q$ to $G$ by a rational representation $\rho$ of $M$. Note that this condition is automatic if $L_K$ is induced from $B_H$ to $H$ by a character of $T_H$, by transitivity of the induction ($B_H$ to $Q$ and $Q$ to $H$). Let $I_Q = \{g \in H(\mathbf{Z}_p); g \bmod. p \in Q(\mathbf{Z}/p\mathbf{Z})\}$. In the text, we shall consider a suitable lattice $L$ of $L_K$, stable by $I_Q$, hence by the $p$-component of ideles in $U_{\tilde{Q}}(p^r)$ for any $r \geq 1$; let us put $R = L_K/L$; it is a discrete $\mathcal{O}$-module with a continuous action of $U_{\tilde{Q}}(p^r)$ (through its $p$-adic component); we consider the sheaf $\tilde{R}$ of locally constant sections of the covering

$$H_\mathbf{Q} \backslash H_\mathbf{A} \times R / U_{\tilde{Q}}(p^r) U_\infty$$

$$\downarrow$$

$$H_\mathbf{Q} \backslash H_\mathbf{A} / U_{\tilde{Q}}(p^r) U_\infty = S_r(U)$$

where the action defining the covering space is given by

$$\gamma(h,v)uu_\infty = (\gamma h u u_\infty, u_p^{-1}.v) \text{ for } \gamma \in H_\mathbf{Q},\ uu_\infty \in U_{\tilde{Q}}(p^r)U_\infty \text{ and } v \in R.$$

# 2 Flags and Hecke operators

In this section, we view $H^1$ and its subgroups as group schemes over $\mathbf{Z}_p$; $G^1$ is viewed over $\mathbf{r}_v$ for each $v|p$.

## 2.1 Parahorics and flags

Let $Q$ be a $J$-proper SPS; let $M'$ be the Levi subgroup of $Q' = Q \cap H^1$. In the case of the standard Borel subgroup $Q = f_* f^* B$ of $H$, we write $Q' = B'_H$. We have a Levi decomposition $Q' = M'Q^+$; let $Z_M$ and $C_M = M/M^1$ be the center resp. the "cocenter" of $M$, and $Z_{M'}$ resp. $C_{M'} = M'/M^1$ be the center, resp. the "cocenter" of $M'$; we put $\tilde{Q} = \tilde{Q}' = M^1 Q^+$. Let us put, for each $v \in S_p$, $X_v = G^1/Q'_v$ and $Y_v = G^1/Q_v^+$. Let $\pi_v : Y_v \to X_v$ be the structural map; it is a morphism of $\mathbf{r}_v$-scheme. Note that $X_v$ is a projective scheme and $Q_v$ is an $M'_v$-bundle over $X_v$. We write $\underline{0}_v$ for the marked point on $X_v$ given by the trivial class. Put $Y_{Q'} = \prod_{v|p} Y_v$, $X_{Q'} = \prod_{v|p} X_v$, $\pi_{Q'} = \prod_{v \in J} \pi_v$ and $\underline{0}_{Q'} = \prod_{v|p} \underline{0}_v$.

Let $Q^- = \prod_{v|p} \operatorname{Res}_{\mathbf{Z}_p}^{\mathbf{r}_v} Q_v^-$ be the unipotent subgroup opposite of $Q^+$ (hence $Q_v^- = 1$ for $v \notin J$). For each $r \geq 1$, we put

$$I_r = \{h \in H(\mathbf{Z}_p); h \bmod p^r \in Q(\mathbf{Z}/p^r\mathbf{Z})\}.$$

and

$$I'_r = I_r \cap H^1(\mathbf{Z}_p) = I'_r = Q'^-(p^r \mathbf{Z}_p) Q'(\mathbf{Z}_p)$$

We introduce also

$$\tilde{I}_r = \{h \in H(\mathbf{Z}_p); h \bmod p^r \in \tilde{Q}(\mathbf{Z}/p^r\mathbf{Z})\}.$$

and

$$\tilde{I}'_r = \tilde{I}_r \cap H^1(\mathbf{Z}_p)$$

For $r = 1$, we drop the index 1: $I = I_1$, $I' = I'_1$ and $\tilde{I} = \tilde{I}_1$.

**Terminology:** We refer to the subgroups $I_r$ and $I'_r$, resp. $\tilde{I}_r$ and $\tilde{I}'_r$ as the level $p^r$ parahoric subgroups of $H(\mathbf{Z}_p)$, resp. $H^1(\mathbf{Z}_p)$ associated to $Q$. We call them Iwahori subgroups when $Q$ is the standard Borel subgroup $B_H$ of $H$.

For $v|p$, let

$$I_{v,r} = \{g \in G(\mathbf{r}_v); g \bmod p^r \in Q_v(\mathbf{r}_v/(p^r))\}.$$

and $I'_{v,r} = I_{v,r} \cap G^1(F_v)$. Then, $I'_r = \prod_{v|p} I'_{v,r}$. Similarly, we have $\tilde{I}'_r = \prod_{v|p} \tilde{I}'_{v,r}$.

We put

$$X^a_{v,r} = I'_{v,r}/Q'_v(\mathbf{r}_v),$$

and

$$X^a_{Q',r} = I'_r/Q'(\mathbf{Z}_p) \text{ and } X^a_{Q'} = X^a_{Q',1},$$

so, $X^a_{Q',r} = \prod_{v|p} X^a_{v,r}$; similarly:

$$Y^a_{v,r} = I'_{v,r}/Q^+_v(\mathbf{r}_v) = Q_v^-(p^r\mathbf{r}_v) \times M'_v(\mathbf{r}_v)$$

Let

$$Y^a_{Q',r} = I'_r/Q^+(\mathbf{Z}_p) \text{ and } Y^a_{Q'} = Y^a_{Q',1}$$

hence, $Y^a_{Q',r} = \prod_{v|p} Y^a_{v,r}$.

For any $r \geq 1$, $Y^a_{Q',r}$ is the inverse image by $\pi_{Q'}$ of $X^a_{Q',r}$; these sets are p-adic open in $Y_{Q'}(\mathbf{Q}_p)$, resp. in $X_{Q'}(\mathbf{Q}_p)$.

They can be viewed as $\mathbf{Z}_p$-points of schemes whose functorial description is as follows: for any $\mathbf{Z}_p$-algebra $A$, let $A_v = A \otimes \mathbf{r}_v$ for each $v|p$.

- for $Q'_v = P'$, $X^a_{v,r}(A_v)$ is the set of maximal isotropic direct factors $E$ in $W \otimes A_v$ whose reduction modulo $p^r$ is the standard lagrangian $<e_1, e_2>$; $Y^a_{v,r}(A_v)$ is the set of pairs $(E, \phi)$ of a maximal isotropic $A_v$-submodule $E$, direct factor in $W \otimes A_v$ and an isomorphism $\phi : E \cong A_v^2$ whose reduction modulo $p^r$ is the standard isomorphism: $e_1 \mapsto (1,0)$, $e_2 \mapsto (0,1)$. The morphism $\pi_v$ consists in forgetting $\phi$.

- for $Q'_v = P'^*$, $X^a_{v,r}(A_v)$ is the set of rank 1 free direct factors $E_1$ in $W \otimes A_v$ whose reduction modulo $p^r$ is the line $<e_1>$; $Y^a_{v,r}(A_v)$ is the set of pairs $(E_1, \phi_1)$ of an $E_1$ as above and an isomorphism $\phi_1 : E_1 \cong A_v$ such that $\phi_1$ mod.$p^r$ coincides with the standard isomorphism $e_1 \mapsto 1$. The morphism $\pi_v$ consists in forgetting $\phi_1$.

- for $Q'_v = B'$, $X^a_{v,r}(A_v)$ is the set of flags $(E_1, E)$ of isotropic submodules which are free direct factors in $W \otimes A_v$ and whose reduction modulo $p^r$ is the standard flag $(<e_1>, <e_1, e_2>)$; $Y^a_{v,r}(A_v)$ is the set of pairs $((E_1, E), (\phi_1, \phi))$ consisting in a flag $(E_1, E)$ as above, and isomorphisms $\phi_1 : E_1 \cong A_v$ and $\phi : E \cong A_v^2$ such that $\phi$ restricted to $E_1$ coincides with $\phi_1$ and $\phi$ modulo $p^r$ coincide with the standard isomorphisms: $e_1 \mapsto (1,0)$, $e_2 \mapsto (0,1)$. The morphism $\pi_v$ consists in forgetting $(\phi_1, \phi)$.

- for $Q'_v = G$, $X^a_{v,r}(A_v) = \{.\}$ and $Y^a_{v,r}(A_v)$ is the set of symplectic bases of $W \otimes A_v$. It is in canonical bijection with $G^1(\mathbf{r}_v)$.

Then, we have $Y^a_{Q',r}(A) = \prod_{v|p} Y^a_{v,r}(A_v)$ and similarly for $X^a_{Q',r}(A)$.

Using these $p$-adic open subsets, we can define lattices stable by our parahoric subgroups, as follows.

Let $\mathcal{O}$ be the valuation ring of a finite extension $K$ of $\mathbf{Q}_p$ containing the Galois closure of a compositum $\prod_{v|p} F_v$ (and sufficiently big when necessary); we fix a family of finite free $\mathcal{O}$-modules $(V_v)_{v|p}$ and a family $(\rho_v)_{v|p}$ of group scheme morphisms $\rho_v : M_v \to GL_{\mathcal{O}}(V_v)$ defined over $\mathbf{Z}_p$ as an algebraic representation such that $\rho_v \otimes_{\mathcal{O}} K$ is absolutely irreducible;

let us put $V = \bigotimes_{v|p} V_v$; $\rho = \bigotimes_{v|p} \rho_v : M \to GL(V)$ is an algebraic representation of $M$ over $\mathbf{Z}_p$; it is absolutely irreducible. Note here that we have fixed representations of the Levi subgroups $M_v$ of $G$; that is, we have prescribed the central action (instead of imposing only their behavior on $M'_v = M_v \cap G^1$). This will be necessary for the control theorem for the Hecke algebra, although throughout the study of group cohomology this does not play any role. The maximal torus of $M_v$ is unchanged, equal to the maximal torus $T$ of $B$, hence $\rho_v \otimes K$ corresponds to a character $\delta_v \in X^*(T) = Hom_{\mathcal{O}}(Res_{\mathbf{Z}_p}^{\mathbf{r}_v} T, \mathbf{G}_m)$ which is dominant for the ordering defined by $B \cap M_v$ (actually the ordering is on $X^*(T')$ hence we mean the restriction of $\delta_v$ to $T'$).

**Definition 2.1** *We say that $\rho$ is $H$-compatible if the characters $\delta_v$ are actually dominant in $X^*(T')$ for the ordering defined by $B$. This means that the induced representation of $\delta_v$ over $K$ from $B$ to $G$ is absolutely irreducible.*

¿From now on, we fix such a representation $\rho = \bigotimes_{v|p} \rho_v$. Note that a character $\chi_v$ of $C_{M_v} = M_v/M_v^1 = Q_v/\tilde{Q}_v$ can be viewed in $X^*(T)$ via $T = B/B^+ \to C_v$. Recall that we have fixed a $J$-proper SPS $Q$.

**Definition 2.2** *An algebraic character $\chi = \bigotimes_{v|p} \chi_v$ of $C_M = \prod_{v|p} C_v$ is called dominant with respect to $\rho$ if $\delta_v \otimes \chi_v$ is dominant in $X^*(T)$ for each $v|p$. For each place $v|p$, the representation of $M_v$ associated to $\delta_v \otimes \chi_v$ restricted to $T$ is $\rho_v \otimes \chi_v$. For each such $\chi$, we form the representation $\rho \otimes \chi : M \to GL(V)$.*

Then we form a representation of $H_p^1 = H^1(\mathbf{Q}_p)$ by parabolic induction:

$$L(\rho; K) = \{f : Y_{Q'} \to V \otimes \mathbf{Q}_p;\ f \text{ polynomial},\ f(ym) = \rho(m^{-1})f(y) \text{ for } m \in M\}$$

For each system $\chi \in X^*(C_M)$, we can form the representation $L(\rho \otimes \chi; K)$. It is easy to see that $\chi$ is dominant with respect to $\rho$ if and only if $L(\rho \otimes \chi; K)$ is an irreducible $H^1$-module.

In it, we fix the following sublattice:

$$L^a(\rho \otimes \chi; \mathcal{O}) = \{f \in L(\rho \otimes \chi; K);\ f(Y_{Q'}^a) \subset V\}$$

One lets $I'$ act on this module by left translation: $(g.f)(y) = f(g^{-1}y)$. This representation of $I'$ will be called in the sequel the algebraic induction from $M'(\mathbf{Z}_p)$ to $I'$ of $\rho \otimes \chi$. These lattices tensor with $K/\mathcal{O}$ will be used as coefficients of our cohomology groups. Before defining global Hecke operators on these cohomology groups, we study the local Hecke algebra.

## 2.2 Local Hecke algebra

For each $v|p$, we choose a uniformizing parameter $\varpi_v$ of $F_v$. Consider the subset $\Delta_{v,r} = I'_{v,r} D_v I'_{v,r}$ of $G(F_v)$ where $D_v$ is the set of matrices

- $\mu(diag_2(\varpi_v^a); \varpi_v^b)$ for $0 \leq 2a \leq b$ if $Q' = P'$,

- $\mu^*(\varpi_v^a, diag_2(\varpi_v^b))$ with $0 \leq a \leq b$, if $Q' = P'^*$,

- $\tau(\varpi_v^{a_1}, \varpi_v^{a_2}; \varpi_v^b)$ for $0 \leq a_1 \leq a_2 \leq \frac{1}{2}b$, if $Q' = B'$.

When $r = 1$, we simply write $\Delta_v$ for $\Delta_{v,1}$. Let $C$ be any subgroup between $\tilde{I}'_{v,r}$ and $I'_{v,r}$.

**Remark:** The local Hecke algebra

$$\mathcal{H}_v(C) = \mathcal{O}[C \backslash \Delta_{v,r} / C]$$

depends on the choice of $\varpi_v$; however the main object of interest to us, the nearly ordinary idempotent, won't depend on that choice.

We use Shimura's notation $[C\xi C]$ for an element of the canonical basis of $\mathcal{H}_v(C)$. The multiplication is given by $[C\xi C].[C\xi' C] = \sum_\eta c_\eta [C\eta C]$ where $C\xi C \xi' C = \bigcup_\eta C\eta C$ is a disjoint union, and $c_\eta = \sharp\{(i,j); C\xi_i \xi'_j = C\eta\}$ (it depends only on the double class of $\eta$).

**Proposition 2.1** $\Delta_{v,r}$ *is a semigroup; the ring $\mathcal{H}_v(C)$ is isomorphic to a polynomial ring over the $\mathcal{O}$-algebra of a finite abelian group:*

$$\mathcal{O}[I'_{v,r}/C][T_{v,i}; i \in A]$$

*where*

- *the group $I'_{v,r}/C$ acts through the double classes $[CgC]$, $g \in I'_{v,r}$,*

- *$A$ is a subset of $\{0, 1, 2\}$: $J = \{0, 1\}$ for $Q = P$, $A = \{0, 2\}$ for $Q = P^*$, and $J = \{0, 1, 2\}$ if $Q = B$,*

- *$T_{v,0} = [C\varpi_v C]$, $T_{v,1} = [C\mu(1_2; \varpi_v)C]$, and $T_{v,2} = [C\mu^*(1; \varpi_v.1_2)C]$*

*In particular, this algebra is commutative.*

**Proof:** The key is to establish the formula $[C\xi C].[C\xi' C] = [C\xi\xi' C]$ which relies on the multiplicativity of the degree map: $\xi \mapsto d(\xi) = h$ where $C\xi C = \bigcup_{i=1}^h C\xi_i$. The computation is similar to that of [17], Proposition .

Let $Q'$ be a parabolic subgroups of $H^1$. Consider its associated parahoric subgroup $I'_r$ and $\tilde{I}'_r$ of level $p^r$. For $\Delta_r = \prod_{v|p} \Delta_{v,r}$ (or simply $\Delta$ if $r = 1$), we consider the $p$-semilocal Hecke algebra $\mathcal{H}_{p,r} = \mathcal{O}[\tilde{I}'_r \backslash \Delta_r / \tilde{I}'_r]$ (or $\mathcal{H}_p$ when $r = 1$); there is a canonical isomorphism $\mathcal{H}_{p,r} =$ $\otimes_{v|p} \mathcal{H}_v(\tilde{I}'_{v,r})$. Hence, $\mathcal{H}_{p,r}$ is a polynomial algebra over the group algebra of $C_{M'}(\mathbf{Z}/p^r\mathbf{Z})$ where $M'$ is the Levi system of $Q'$. Note also that one can write $\mathcal{H}_{p,r} = \mathcal{H}_{J,r} \times \mathcal{H}_r^{(J)}$ where one defines the $J$-semilocal Hecke algebra as $\mathcal{H}_{J,r} = \otimes_{v \in J} \mathcal{H}_v(\tilde{I}'_{v,r})$ and its complement (the non-$J$-part) by $\mathcal{H}_r^{(J)} = \otimes_{v \notin J} \mathcal{H}_v(Sp_4(\mathbf{r}_v))$.

## 2.3 Action on flags

Let $Q \subset H$ as above. We keep the Notation of Section 2.1. For any $r \geq s \geq 1$, we have an obvious action of $\tilde{I}'_r = \prod_{v|p} \tilde{I}'_{v,r}$ on $Y^a_{Q',s}$ by left translation and of $M'(\mathbf{Z}_p)$ by right translation. We wish to extend the action on the left to the semi-group $\Delta_r = \prod_{v|p} \Delta_{v,r}$. It is enough to define the action of $D_v$ on $I'_{v,s}/Q_v^+(\mathbf{r}_v)$ for each $v|p$.

- For $Q = P$, let $d_v = \mu(t.1_2; x) \in D_v$ (so $0 \leq 2.\text{ord}_v(t) \leq \text{ord}_v(x)$); let $(E_v; \phi_v) \in Y_{v,s}$; note that the image $d_v E_v$ of $E_v$ can be written $tE'_v$ for some direct factor $E'_v$ in $W \otimes \mathbf{r}_v$ and $\phi'_v = \phi_v \circ (td_v^{-1})$ is an isomorphism $E'_v \cong \mathbf{r}_v^2$, whose reduction mod. $\overline{\varpi}_v^s$ is the canonical isomorphism $e_1 \mapsto (1,0)$, $e_2 \mapsto (0,1)$; so we put $d_v.(E_v, \phi_v) = (E'_v, \phi'_v)$.

- For $Q = P^*$, let $d_v = \mu^*(x, y.1_2) \in D_v$ (so $0 \leq \text{ord}_v(x) \leq \text{ord}_v(y)$); let $(E_{v,1}; \phi_{v,1}) \in Y_{v,s}$; similarly, $d_v E_{v,1} = xE'_{v,1}$ for a direct factor $E'_{v,1}$ in $W \otimes \mathbf{r}_v$ and $\phi'_{v,1} = \phi_{v,1} \circ (xd_v^{-1})$ is an isomorphism $E'_{v,1} \cong \mathbf{r}_v$, whose reduction mod. $\overline{\varpi}_v^s$ is the canonical isomorphism $e_1 \mapsto 1$; so we can put $d_v.(E_{v,1}, \phi_{v,1}) = (E'_{v,1}, \phi'_{v,1})$.

- For $Q = B$, let $d_v = \tau(t, t'; x) \in D_v$ with $0 \leq \text{ord}_v(t) \leq \text{ord}_v(t') \leq \frac{1}{2}\text{ord}_v(x)$); let $y = ((E_{v,1}, E_v); \phi_{v,1}, \phi_v) \in Y_{v,s}$; again, $d_v E_{v,1} = tE'_{v,1}$ and let $E'_v = W \cap (d_v E_v \otimes_{\mathbf{r}_v} F_v)$; we have $E'_{v,1} \subset E'_v \subset W \otimes \mathbf{r}_v$, each term being direct factor in the next one; to define $\phi'_v : E'_v \cong \mathbf{r}_v^2$; we take a basis $(g_1, g_2)$ of $E_v$ such that $\phi_v(g_1) = (1,0)1$ and $\phi_v(g_2) = (0,1)$, then we can pick a basis $(g'_1, g'_2)$ of $E'_v$ such that $d_v g_1 = tg'_1$ (so, we put $\phi'_{v,1}g'_1 = 1$) and $d_v g_2 = t'g'_2 + \lambda g'_1$ for some $\lambda \in F_v$, then we put $\phi'_v(g'_2) = (0,1)$. The reduction of this isomorphism mod. $\overline{\varpi}_v^s$ coincides with the canonical isomorphism; hence we can put $d_v.y = y' = (E'_{v,1}, E'_v); \phi'_{v,1}, \phi'_v)$.

Then, one can let any element $\delta_v = ud_v u' \in \Delta_{v,r}$ act on $Y_{v,s}$ by composing the actions of $u$, $d_v$ and $u'$. This yields the desired left action of $\Delta_r$ on $Y_{Q',s}$.

Actually, for further calculations, it will be also useful to give the group-theoretic description of this action. It simply amounts to the following: for $d_v \in D_v$, for any $g \in I_{v,s}$, resp. $I'_{v,s}$ there exists $g_1 \in I_{v,s}$ (resp. $I'_{v,s}$) such that

$$d_v g d_v^{-1} \equiv g_1 \mod Q^+(F_v)$$

moreover the coset $g_1 Q^+(\mathbf{r}_v) \in Y_{v,s}$ is uniquely determined by this congruence. Then, we have:

$$d_v.(gQ^+(\mathbf{r}_v)) = g_1 Q^+(\mathbf{r}_v)$$

Consider now $d_{v,1} = \mu(1_2; \varpi_v)$ and $d_{v,2} = \mu^*(1, \varpi_v.1_2)$, $d_{v,3} = d_{v,1}d_{v,2}$ and $d_Q = \prod_{v|p} d_{v,i}$. Write $P_1 = P$, $P_2 = P^*$ and $P_3 = B$. Let us establish, following Hida [17], the crucial Contraction Property of this action. We give two proofs, one in the language of buildings, the other being group-theoretic. In the following proposition, the map $\pi_{Q'}$ is the one defined in paragraph 3.1.

**Proposition 2.2** For $r \geq s \geq 1$, for $Q = P_i$, the operator $d_Q^{r-s}$ contracts $Y_{Q',s}^a(\mathbf{Z}/p^r\mathbf{Z})$ to $\pi_{Q'}^{-1}(\underline{0}_{\mathbf{Z}/p^r\mathbf{Z}})$

**Proof:** It is enough to consider each $v|p$ separately. Let $e = ord_v(p)$.

- for $Q = P$, take a lagrangian submodule $E$ with basis $(g_1, g_2)$ congruent to the standard basis $(e_1, e_2)$ modulo $p^s$. Then $d_{v,1}^{e(r-s)}(g_1, g_2)$ is congruent modulo $p^r$ to $\alpha(e_1, e_2)$ where $\alpha \in GL(E)$ and $det\ \alpha \equiv 1\ mod.p^s$, as desired.

- for $Q = P^*$, take an isotropic line $E_1$ with basis $g_1$, congruent to $e_1$ modulo $p^s$. Then $d_{v,2}^{e(r-s)}g_1$ is congruent to $(1 + p^s \star)e_1$ modulo $p^r$.

- for $Q = B$, take a lagrangian flag $(E_1, E)$ with basis $(g_1, g_2)$ congruent to $(e_1, e_2)$ modulo $p^s$. Then by case one, by applying $d_{v,2}^{e(r-s)}$ one can assume $g_1 = (1 + p^s \star)e_1$; then, applying $d_{v,1}^{r-s}$ does not change $g_1$ and sends $g_2$ to $ae_1 + (1 + p^s \star)e_2$ with $a \in \mathbf{r}_v$. So, $d_{v,3}^{e(r-s)}(E_1, E)$ is the standard flag $(<e_1>, <e_1, e_2>)$.

A group-theoretic proof of this Proposition is to notice that for $Q_v = P_i$, for any root $\alpha$ in $R_Q^+$, one has $v(\alpha(d_{v,i})) < 0$, hence if $g \in I'_{v,s}$, then the element $g_1$ defined by 2.3 belongs to $I'_{v,s+1}$. Compare Lemma 4.3 below. ∎

For $r \geq 1$, one deduces from the left action of $\Delta_r$ on $Y_{Q'}^a$ a left action of $\Delta_r^{-1}$ on $L_{Q'}^a(\rho, \chi; \mathcal{O})$ given by $(\delta^{-1}.f)(y) = f(\delta.y)$; by tensorization by $K/\mathcal{O}$, it extends to $R = L_{Q'}^a(\rho, \chi; \mathcal{O}) \otimes K/\mathcal{O}$. This will be a typical coefficient module for the group cohomology we have in mind.

## 2.4 Nearly ordinary part of the cohomology groups.

Given $Q \subset H_{/\mathbf{Z}_p}$ as above and $U$ a level group unramified at $p$ we have defined in 1.2 the level groups $U_{\tilde{Q}}(p^r)$. Similarly, we consider

$$U_Q(p^r) = \{h \in H(\hat{\mathbf{Z}}); h\ mod.\ p^r \in Q(\mathbf{Z}/p^r\mathbf{Z})\}$$

For $N \geq 1$, we put $\Gamma = (U \times H_\infty) \cap H_\mathbf{Q}^1$; for each $r \geq 1$ we define

$$\Gamma_1(p^r) = \Gamma \cap (U_{\tilde{Q}}(p^r) \times H_\infty^1) \quad \text{and} \quad \Gamma_0(p^r) = \Gamma \cap (U_Q(p^r) \times H_\infty^1)$$

In the sequel, these congruence subgroups are viewed as embedded either in the archimedean component $H^1_\infty$ or in the $p$-adic component $H^1_p$ of $H^1_{\mathbf{A}}$. The context should make it clear which embedding is used. We assume

$$\Gamma \text{ is torsion-free}$$

Let $\mathcal{Z} = H^1_\infty/(U_\infty \cap H^1_\infty)$ be the Siegel-Hilbert space of genus 2 over $F$; it is a global hermitian domain of dimension $3d$; then, $H_\infty/U_\infty$ is the union of two copies of $\mathcal{Z}$. For any $\Gamma$-module $R$, for any subgroup $\Gamma' \subset \Gamma$ let $\tilde{R}$ be the local system associated to $R$ on the $3d$-dimensional complex manifold $\Gamma'\backslash\mathcal{Z}$. There are canonical isomorphisms: $H^q(\Gamma', R) \cong H^q(\Gamma'\backslash\mathcal{Z}, \tilde{R})$

We fix a $p$-adic field $K$ with ring of integers $\mathcal{O}$. Let $R$ be a discrete $\mathcal{O}$-module with left action of $\Delta_r^{-1}$. In particular, the group $\Gamma_*(p^r)$ ($* = 0, 1$) acts on $R$ by its $p$-adic embedding $\Gamma_*(p^r) \subset I'_r$. For each integer $0 \leq q \leq 6d$, consider the cohomology group $H^q(\Gamma_*(p^r), R)$ for $r \geq 1$.

As in the beginning of Section 4 of [17], we introduce global versions $T_Q$ of the local Hecke operators introduced in Sect.2.2 and let them act on $H^q(\Gamma_*(p^r), R)$; they are defined using global double classes $\Gamma_*(p^r)\xi_Q\Gamma_*(p^r)$ ($* = 0$ or $1$) for suitable elements $\xi_Q \in H_{\mathbf{Q}}$. Recall that (cf. [16] Sect.1.10 or [17] Sect.4): for $\xi \in H_{\mathbf{Q}}$, and $R$ a left module over the semigroup generated by $\Gamma'$ and $\xi^{-1}$, the double class $[\Gamma'\xi\Gamma']$ acts on $H^q(\Gamma', R)$ by sending a $q$-homogeneous cocycle $u$ to the cocycle $v$ defined as follows. If $\Gamma'\xi\Gamma' = \bigcup \Gamma'\xi_j$, for $\gamma \in \Gamma'$ let $\gamma_j \in \Gamma'$ such that $\xi_j\gamma = \gamma_j\xi_{j'}$ for some $j'$; then $v(\gamma^{(0)}, ..., \gamma^{(q)}) = \sum_j \xi_j^{-1}u(\gamma_j^{(0)}, ..., \gamma_j^{(q)})$. One sees easily that this does not depend on the choice of the representatives $\xi_j$ and that it is well defined on cohomology. Consider the endomorphism $T_Q = [\Gamma_*(p^r)\xi_Q\Gamma_*(p^r)]$ of $H^q(\Gamma_*(p^r), R)$; although the operator $T_Q$ will depend on several choices, it won't be the case for the nearly ordinary idempotent $e_Q$ defined by

$$e_Q = \lim_{n\to\infty} T_Q^{n!}$$

The image of the cohomology by this idempotent is denoted by

$$H^q_{Q-no}(\Gamma_*(p^r), R).$$

Note that these groups form an inductive system with respect to the restriction maps when $r$ grows.

It remains to construct $\xi_Q$.

Let $x_Q$ be the element of $H^1(\mathbf{Q} \otimes \hat{\mathbf{Z}})$ whose components are 1 outside $J$ and is given at each $v \in J$ by

- $\mu(p^{-2}.1_2; 1)$ if $Q_v = P_1$,

- $\mu^*(p^{-2}, 1_2)$ if $Q_v = P_2$,

- $\tau(p^{-2}, p^{-1}; 1) = \mu(p^{-1}.1_2; 1)\mu^*(p^{-1}, 1_2)$, if $Q_v = B = P_3$,

By Strong Approximation Theorem, there exists $\zeta_Q \in H^1_{\mathbf{Q}}$ such that $\zeta_Q \equiv x_Q \mod U_{\tilde{Q}}(p^r)$. By finiteness of the class group of $F$, there exists $h \geq 1$ and $\pi \in F^\times$ such that $\pi p^{-h} \equiv 1 \mod. \, p^r O_{F,v}$ for each $v \in J$ and $\operatorname{ord}_v(\pi) = 0$ for any other finite place $v$. Define

$$\xi_Q = \pi^2 \zeta_Q^h$$

Then, via the embedding $F \subset F_v$, one has

$$\xi_Q \equiv d_{v,i}^{2h} \mod. \, I_{v,r}$$

Hence, the double class $\Gamma_*(p^r)\xi_Q \Gamma_*(p^r)$ has for $v$-adic completion

- $I'_{v,r} d_{v,i} I'_{v,r}$ for $v \in J$,

- $G'(\mathbf{r}_v)$ for $v \in S_p \backslash J$,

- $U'_\ell \pi^2 U'_\ell = U'_\ell \pi^2$ for $U' = U \cap H^1_f$ for $\ell \notin S_p$

Note that for any unit $u \in \mathbf{r}_v$, for $n$ large enough $u^{n!} \equiv 1 \mod. \, p^r$. Hence for $n$ large enough, the local double class of $\xi_Q^{n!}$ is that of $1$ at $v$ outside $J$, and that of $d_Q^{2hn!}$ at $v \in J$. In particular, the idempotent $e_Q$ associated to $T_Q$ does not depend on the choice of $\pi$.

# 3 Control theorem for the full cohomology

## 3.1 Induced modules

As before, we fix a $J$-proper SPS $Q$ with Levi subgroup $M_{/\mathbf{Z}_p}$ and a representation $V = \bigotimes_{v|p} V_v$ of $M_{/\mathbf{Z}_p}$, given by $\rho: M \to GL_{\mathcal{O}}(V)$. When we want to put emphasis on the set $J$ of places where $Q_v$ is a proper parabolic of $G_v$, we write $V = V_J \otimes V^J$ where $V_J = \bigotimes_{v \in J} V_v$ and $V^J = \bigotimes_{v \notin J} V_v$. Let $A_r = p^{-r}\mathcal{O}/\mathcal{O}$ ($r = 1, ...$), and $A = A_\infty = K/\mathcal{O}$. For each character $\chi$ of the cocenter $C_{M'}$ of $M'$, dominant with respect to $\rho$, recall that in 2.1 we have defined an $\mathcal{O}$-module of finite type $L^a(\rho \otimes \chi; \mathcal{O})$ which we decompose when needed as $L^a(\rho \otimes \chi; \mathcal{O}) = L^a_J(\rho \otimes \chi; \mathcal{O}) \otimes V^J$ where $L^a_J(\rho \otimes \chi; \mathcal{O}) = \bigotimes_{v \in J} \operatorname{ind}_{Q'_v}^{I'_v}(\rho_v \otimes \chi_v)$. We call $L^a(\rho \otimes \chi; \mathcal{O})$ the algebraic induction of $\rho \otimes \chi$ from $Q'(\mathbf{Z}_p)$ to $I'$. we form

$$L^a(\rho \otimes \chi; A_r) = L^a(\rho \otimes \chi; \mathcal{O}) \otimes A_r$$

**Remark:** Let us write $L(Y^a_{Q'}(\mathbf{Z}/p^s\mathbf{Z}), \rho \otimes \chi; A_r)$ for the submodule of $L^a(\rho \otimes \chi; A_r)$ consisting of functions which factor through $Y^a_{Q'}(\mathbf{Z}/p^s\mathbf{Z})$; then, for any finite $r \geq 1$, there exists an integer $s \geq 1$ such that

$$L^a(\rho \otimes \chi; A_r) = L(Y^a_{Q'}(\mathbf{Z}/p^s\mathbf{Z}), \rho \otimes \chi; A_r).$$

This follows easily from the finiteness of $L^a(\rho \otimes \chi; \mathcal{O})$ as $\mathcal{O}$-module.

We consider also another induced module associated to $\rho$:

$$\mathcal{C}(\rho^1; A) = \{f : Y_{Q'}^a \to V \otimes A; f \text{ locally constant}; f(xm) = \rho(m^{-1})f(x) \text{ for any } m \in M^1(\mathbf{Z}_p)\}$$

We shall refer to it as the smooth induction from $M^1$ to $I'$ of $\rho|_{M^1}$. This module is not of cofinite type over $\mathcal{O}$. Note also a difference between the two definitions, besides the nature of the functions, namely that for the latter we impose an equivariance condition only for $m \in M^1$. These modules carry a left action of $\Delta^{-1}$ given by $(\delta^{-1}f)(y) = f(\delta.y)$ where $\Delta$ acts on the left on $Y_{Q'}^a$ as defined in Section 2.3. The latter is a smooth admissible representation of $I'$ in the following sense: let

$$H^1(p^r) = Ker(H^1(\mathbf{Z}_p) \to H^1(\mathbf{Z}/p^r\mathbf{Z}))$$

then, for any $r \geq 1$: $\mathcal{C}(\rho; A)^{H^1(p^r)}$ consists in functions on the finite set $Y_{Q'}^a(\mathbf{Z}/p^r\mathbf{Z})$, hence is of (co)finite type, and

$$\mathcal{C}(\rho; A) = \bigcup_r \mathcal{C}(\rho; A)^{H^1(p^r)}$$

Let us write $\mathcal{C}(Y_{Q'}^a(\mathbf{Z}/p^r\mathbf{Z}), \rho^1; A)$, resp. $\mathcal{C}(M'(\mathbf{Z}/p^r\mathbf{Z}), \rho^1; A)$ for the module of $V \otimes A$-valued functions satisfying the equivariance condition for $M^1(\mathbf{Z}/p^r\mathbf{Z})$ and defined on $Y_{Q'}^a(\mathbf{Z}/p^r\mathbf{Z})$, resp. on $M'(\mathbf{Z}/p^r\mathbf{Z})$.

We consider now two short exact sequences which will permit the comparison of these two notions of induction. Namely, for all $r \geq 1$,

$$0 \to K_r \to L^a(\rho \otimes \chi; A_r) \xrightarrow{\phi_r} V \otimes \chi \otimes A_r \to 0 \tag{3.1}$$

where $\phi_r$ is the evaluation map at the marked point $O_Y$ of $Y_{Q'}^a$. Observe that for any $s \geq r$ defined as in Remark of section 4.1, it is an exact sequence of $\Delta_s^{-1}$-module, hence of $\Gamma_0(p^s)$-modules; however we view it only as a sequence of $\Gamma_1(p^s)$-modules. The reason is that we want to let vary the characters $\chi$; actually as a $\Gamma_1(p^s)$-module, we have $V \otimes \chi \otimes A_r = V \otimes A_r$.

Similarly, for $s \geq 1$,

$$0 \to \mathcal{K}_s \to \mathcal{C}(Y_{Q'}^a(\mathbf{Z}/p^s\mathbf{Z}), \rho^1; A_s) \xrightarrow{\psi_s} \mathcal{C}(M'(\mathbf{Z}/p^s\mathbf{Z}), \rho^1; A_s) \to 0 \tag{3.2}$$

where $\psi_s$, is the restriction map from the $\mathbf{Z}/p^s\mathbf{Z}$-points of $Y_{Q'}^a$ to $\pi_{Q'}^{-1}(O_{X,\mathbf{Z}/p^s\mathbf{Z}})$, that is, to the $\mathbf{Z}/p^s\mathbf{Z}$-points of $Q'/Q'^+ = M'$.

We view 3.2 as a sequence of $\Delta_{Q,s}^{-1}$-modules. One should remark that although the middle terms of 3.1 and 3.2 are $\Delta_Q^{-1}$-modules, the maps $\phi_r$ and $\psi_s$ are only $\Delta_{Q,s}^{-1}$-linear; hence the restriction to the above-mentioned semigroups.

## 3.2 Taking $p$-adic limit in Shapiro's lemma

¿From (3.1) and (3.2), we deduce long exact sequences of cohomology: for $s \geq r$,

$$... \to H^q(\Gamma_1(p^s), K_r) \to H^q(\Gamma_1(p^s), L^a_{Q'}(\rho \otimes \chi; A_r)) \xrightarrow{\phi_r}$$

$$\to H^q(\Gamma_1(p^s), V \otimes A_r) \to H^{q+1}(\Gamma_1(p^s), K_r) \to ...$$

and

$$... \to H^q(\Gamma_0(p^r), \mathcal{K}_r) \to H^q(\Gamma_0(p^r), \mathcal{C}(Y^a_{Q'}(\mathbf{Z}/p^r\mathbf{Z}); \rho^1; A_r)) \xrightarrow{\psi_r}$$

$$H^q(\Gamma_0(p^r), \mathcal{C}(M'(\mathbf{Z}/p^r\mathbf{Z}), \rho^1; A_r)) \to H^{q+1}(\Gamma_0(p^r), \mathcal{K}_r) \to ...$$

If we apply the idempotent $e_Q$ to this sequence, it remains exact. We shall prove

**Proposition 3.1** *For any $q \geq 0$, $r \geq 1$, and for $s \geq r$ associated to $r$ as in the Remark of Section 4.1, the nearly ordinary idempotent $e_Q$ annihilates $H^q(\Gamma_1(p^s), K_r)$ and $H^q(\Gamma_0(p^r), \mathcal{K}_r)$.*

**Comment:** It is important to note that only the continuity of the functions $f: Y^a_{Q'} \to A_r$ in $K_r$ (and $\mathcal{K}_r$) is used in the proof of this lemma. The very fact that a power of $T_Q$ contracts $Y^a_{Q'}(\mathbf{Z}/p^s\mathbf{Z})$ to $\pi^{-1}_{Q'}(O_{X,\mathbf{Z}/p^s\mathbf{Z}})$ will do.

**Proof:** Let $u$ be an homogeneous $q$-cocycle with values in $K_r$ or $\mathcal{K}_r$. Since the proof is very similar in the two cases, let us deal with $K_r$ only. By Remark of the Section 4.1, we know that $K_r$ is contained in $L(Y^a_{Q'}(\mathbf{Z}/p^s\mathbf{Z}), V \otimes_\mathcal{O} A_r)$. Let $(\gamma^{(0)}, ..., \gamma^{(q)}) \in \Gamma_1(p^s)^{q+1}$; by Proposition 2, $\xi_Q^{s-1}$ sends $Y^a_{Q'}(\mathbf{Z}/p^s\mathbf{Z})$ into $\pi^{-1}_{Q'}(O_{X,\mathbf{Z}/p^s\mathbf{Z}})$. Let us put $v = u|T_Q^{s-1}$. If $\Gamma_1(p^s)\xi_Q^{s-1}\Gamma_1(p^s) = \bigcup_j \Gamma_1(p^s)\eta_j$, we have for any $y \in Y^a_{Q'}$:

$$v(\gamma^{(0)}, ..., \gamma^{(q)})(y) = \sum_j \eta_j^{-1} u(\gamma_j^{(0)}, ..., \gamma_j^{(q)})(y)$$

where $\eta_j \gamma^{(i)} = \gamma_j^{(i)} \eta_{j'}$. That is,

$$v(\gamma^{(0)}, ... \gamma^{(q)}; y) = \sum_j u(\gamma_j^{(0)}, ..., \gamma_j^{(q)}; \eta_j.y)$$

and $\eta_j$ is of the form $\xi_Q^{s-1} u_j$ where $u_j \in \tilde{Q}'(\mathbf{Z}_p)$. Hence $\eta_j.y \in \pi^{-1}_{Q'}(O_{X,\mathbf{Z}/p^s\mathbf{Z}})$. This proves $u|T_Q^{s-1} \equiv 0$ mod. $p^s$. By letting $s$ grow to infinity, one can then conclude.

Let us put

$$H^q_{Q-no}(\Gamma_1(p^\infty), L^a(\rho \otimes \chi; A)) = \varinjlim_s H^q_{Q-no}(\Gamma_1(p^s), L^a(\rho \otimes \chi; A_r))$$

$$H^q_{Q-no}(\Gamma_1(p^\infty), V \otimes A) = \varinjlim_{s \geq r} H^q_{Q-no}(\Gamma_1(p^s), V \otimes A_r))$$

and
$$H^q_{Q-no}(\Gamma_0(p^\infty), \mathcal{C}(\rho^1; A)) = \varinjlim_r H^q_{Q-no}(\Gamma_0(p^r), \mathcal{C}(Y^a_{Q'}(\mathbf{Z}/p^r\mathbf{Z}); \rho^1; A_r))$$

**Corollary 3.1** *The maps $\phi_r$, resp. $\psi_r$ form a system of compatible isomorphisms (when $r \geq 1$ and $s \geq r$ grow)*

$$H^q_{Q-no}(\Gamma_1(p^s), L^a(\rho \otimes \chi; A_r)) \cong H^q_{Q-no}(\Gamma_1(p^s), V \otimes A_r)$$

*resp.*

$$H^q_{Q-no}(\Gamma_0(p^r), \mathcal{C}(Y^a_{Q'}(\mathbf{Z}/p^r\mathbf{Z}), \rho^1; A_r)) \cong H^q_{Q-no}(\Gamma_0(p^r); \mathcal{C}(M'(\mathbf{Z}/p^r\mathbf{Z}), \rho^1; A_r))$$

*By taking their inductive limit we obtain isomorphisms*

$$\iota_1 : H^q_{Q-no}(\Gamma_1(p^\infty), L^a(\rho \otimes \chi; A)) \cong H^q_{Q-no}(\Gamma_1(p^\infty), V \otimes A) \tag{3.3}$$

*and*

$$H^q_{Q-no}(\Gamma_0(p^\infty), \mathcal{C}(\rho^1; A)) \cong \varinjlim_r H^q_{Q-no}(\Gamma_0(p^r), \mathcal{C}(M'(\mathbf{Z}/p^r\mathbf{Z}), \rho^1; A))$$

We can now state the main result of this subsection:

**Theorem 3.1** *For any character $\chi$ of $C_{M'}$ dominant with respect to $\rho$, there is a canonical isomorphism*

$$\iota_2 : H^q_{Q-no}(\Gamma_0(p), \mathcal{C}(\rho^1; A)) \cong H^q_{Q-no}(\Gamma_1(p^\infty), V \otimes \chi \otimes A) \tag{3.4}$$

**Comment:**

The two ingredients of the proof are:

1) a version of Shapiro's lemma involving smooth induction instead of finite induction, since $\mathcal{C}(\rho; A)$ is a direct limit of the inductions from $\Gamma_1(p^r)$ to $\Gamma_0(p^r)$,

2) the lowering of the level on the $\Gamma_0$-side from $p^\infty$ to $p$, by the so-called Hida's lemma (Lemma 3.1 below).

**Proof:**

We apply Shapiro's lemma, noticing that

$$\mathcal{C}(M'(\mathbf{Z}/p^r\mathbf{Z}), \rho^1; A_r) = Ind_{\Gamma_1(p^r)}^{\Gamma_0(p^r)} V \otimes A_r$$

It yields the canonical identification:

$$H^q_{Q-no}(\Gamma_0(p^r), \mathcal{C}(M'(\mathbf{Z}/p^r\mathbf{Z}), \rho^1; A_r)) = H^q_{Q-no}(\Gamma_1(p^r), V \otimes A_r)$$

Using Corollary 3.1, we obtain isomorphisms:
$$H^q_{Q-no}(\Gamma_0(p^r), \mathcal{C}(Y^a_{Q'}(\mathbf{Z}/p^r\mathbf{Z}); \rho^1; A_r)) \cong H^q_{Q-no}(\Gamma_1(p^r), V \otimes A_r) =$$
$$H^q_{Q-no}(\Gamma_1(p^r), V \otimes \chi \otimes A_r)$$

We observe then that
$$\varinjlim_r H^q_{Q-no}(\Gamma_0(p^r), \mathcal{C}(Y^a_{Q'}(\mathbf{Z}/p^r\mathbf{Z}), \rho^1; A_r)) \cong \varinjlim_r H^q_{Q-no}(\Gamma_0(p^r), \mathcal{C}(\rho^1; A_r)) =$$
$$H^q_{Q-no}(\Gamma_0(p^\infty), \mathcal{C}(\rho^1; A_r))$$

Next, we show that the transition maps in the inductive system defining $H^q_{Q-no}(\Gamma_0(p^\infty), \mathcal{C}(\rho^1; A))$ are isomorphisms, therefore showing that the inductive limit coincides with the first group of the inductive system. For that purpose, we use a variant of a lemma used extensively by Hida [13], Lemma.4.3 (and due to Shimura): Let $\xi_Q \in H(\mathbf{Q})$ be the element constructed in Section 2.4. For $s \geq r \geq 1$, we put:
$$T = [\Gamma_0(p^r)\xi_Q^{s-r}\Gamma_0(p^r)] \,, T' = [\Gamma_0(p^s)\xi_Q^{s-r}\Gamma_0(p^r)] \text{ and } T'' = [\Gamma_0(p^s)\xi_Q^{s-r}\Gamma_0(p^s)]$$

**Lemma 3.1** *With Notation as above, for any $\Delta^{-1}$-module L, one has a commutative diagram:*
$$\begin{array}{ccc} H^q(\Gamma_0(p^r), L) & \to & H^q(\Gamma_0(p^s), L) \\ \downarrow T & \swarrow T' & \downarrow T'' \\ H^q(\Gamma_0(p^r), L) & \to & H^q(\Gamma_0(p^s), L) \end{array}$$
*where the horizontal maps are restriction maps.*

The important consequence of this is that on the $Q$-nearly ordinary part of these cohomology groups, the restriction map admits an inverse because $T$ and $T''$ become automorphisms; therefore the restriction maps induce isomorphisms for any $r \geq 1$:
$$H^q_{Q-no}(\Gamma_0(p^r), L) \cong H^q_{Q-no}(\Gamma_0(p^\infty), L)$$

**Proof:** Let us show the commutativity of the upper triangle. It is a consequence of the equality $\Gamma_0(p^s)\xi_Q^{s-r}\Gamma_0(p^r) = \Gamma_0(p^r)\xi_Q^{s-r}\Gamma_0(p^r)$ which follows from a calculation establishing:
$$\xi_Q^{r-s}\Gamma_0(p^s)\xi_Q^{s-r} \cap \Gamma_0(p^r) = \xi_Q^{r-s}\Gamma_0(p^r)\xi_Q^{s-r} \cap \Gamma_0(p^r)$$

Let us check it. It is a local computation at each $v$ dividing $p$; we may therefore replace $\xi_Q$ by $d_{v,i}$ defined before Proposition 2.2. Since the three cases are similar, let us treat only the case $Q_v = B$; we write a matrix of $I_{v,s}$ as
$$\begin{pmatrix} * & * & * & * \\ (p^s) & * & * & * \\ (p^s) & (p^s) & * & (p^s) \\ (p^s) & (p^s) & * & * \end{pmatrix}$$

it becomes after conjugation by $d_{v,3}^{s-r} = diag(1, p^{s-r}, p^{3(s-r)}, p^{2(s-r)})$:

$$\begin{pmatrix} * & (p^{s-r}) & (p^{3(s-r)}) & (p^{2(s-r)}) \\ (p^r) & * & (p^{2(s-r)}) & (p^{s-r}) \\ (p^{3r-2s}) & (p^{2r-s}) & * & (p^r) \\ (p^{2r-s}) & (p^r) & (p^{s-r}) & * \end{pmatrix}$$

The point of this computation is that the blocks corresponding to the lower unipotent $Q^-$ have $p$-order less than $r$, therefore the intersection of the group of such matrices with $I_{v,r}$ is made of the matrices whose blocks constituting the lower unipotent of $Q$ have $p$-order $\geq r$ and those constituting the upper unipotent satisfy divisibility conditions identical to those defining $d_{v,3}^{r-s} I_{v,r} d_{v,3}^{s-r} \cap I_{v,r}$. This shows the equality. For the lower triangle, one observes that

1. the restriction map from $\Gamma_0(p^r)$ to $\Gamma_0(p^s)$ is given by the double class action $[\Gamma_0(p^r) 1_2 \Gamma_0(p^s)]$ and for any $\xi \in H_{\mathbf{Q}}$:

$$[\Gamma_0(p^s) \xi \Gamma_0(p^r)][\Gamma_0(p^r) 1_2 \Gamma_0(p^s)] = [\Gamma_0(p^s) \xi \Gamma_0(p^r)]$$

2. the computation above shows the equality of double cosets:

$$\Gamma_0(p^s) \xi_Q^{s-r} \Gamma_0(p^r) = \Gamma_0(p^s) \xi_Q^{s-r} \Gamma_0(p^s)$$

These two observations yield the desired commutativity. This concludes the proof of the lemma.

We apply it to the left hand side of 3.3 for the pair $(r, s) = (1, s)$. By taking the inductive limit over $r$ of these isomorphisms, we obtain

$$H^q_{Q-no}(\Gamma_0(p), \mathcal{C}(\rho^1; A)) \cong H^q_{Q-no}(\Gamma_1(p^\infty), V \otimes A)$$

## 3.3 Independence of the weight

For each element $\zeta \in C_{M'}(\mathbf{Z}_p)$, we define an automorphism $<\zeta>_{\rho, \chi}$ of

$$H^q_{Q-no}(\Gamma_1(p^\infty), L^a(\rho \otimes \chi; A))$$

as follows: first, for each $r \geq 1$, one views $\zeta$ as an element of $C_{M'}(\mathbf{Z}/p^r\mathbf{Z}) = \Gamma_0(p^r)/\Gamma_1(p^r)$ (compatibly when $r$ grows) and we choose for each $r \geq 1$ a lifting $\gamma_\zeta$ of $\zeta$ in $\Gamma_0(p^r)$; then, one puts for $u \in Z^q(\Gamma_1(p^r), L^a(\rho \otimes \chi; A))$ and $y \in Y^a_{Q'}$:

$$(<\zeta>_{\rho,\chi} u)(\gamma_0, ..., \gamma_q)(y) = u(\gamma_\zeta^{-1} \gamma_0 \gamma_\zeta, ..., \gamma_\zeta^{-1} \gamma_q \gamma_\zeta)(\gamma_\zeta^{-1} y)$$

These automorphisms are compatible when $r \geq 1$ grows: their inverse limit defines $<\zeta>_{\rho,\chi}$.

On the other hand, $Z_{M'}(\mathbf{Z}_p)$ acts on $\mathcal{C}(\rho^1; A)$ by: $(z.f)(y) = f(yz^{-1})$. This action commutes with the action of $\Gamma_0(p)$ hence defines automorphisms in cohomology.

Let $i$ be the natural isogeny $i : Z_{M'} \to C_{M'}$. Let $\omega$ denote the central character of $\rho$ and $\omega_\chi = \omega \times (\chi \circ i)$; it is the central character of $\rho \otimes \chi$.

With these Notation, we can draw from Theorem 3.1 the following corollary

**Corollary 3.2** *For any character $\chi$ of $C_{M'}$ dominant with respect to $\rho$, there exists a canonical isomorphism*

$$\iota_\chi : H^q_{Q-no}(\Gamma_1(p^\infty), L^a(\rho \otimes \chi; A)) \cong H^q_{Q-no}(\Gamma_0(p), \mathcal{C}(\rho^1; A))$$

*such that for any $z \in Z_{M'}(\mathbf{Z}_p)$ and for $\zeta = i(z)$ in $C_{M'}$,*
*for any $c \in H^q_{Q-no}(\Gamma_1(p^\infty), L^a(\rho \otimes \chi; A))$, we have:*

$$z.\iota_\chi(c) = \omega_\chi(z).\iota_\chi(<\zeta^{-1}>_{\rho,\chi} c)$$

**Comment:** This statement can be viewed as a $p$-adic version of the Matsushima-Murakami theorem (see [3], chapter 7). It interprets the infinite level cohomology of finite (co)-rank modules in terms of a finite level cohomology of an (infinite corank) smooth admissible module. Moreover, this module does not depend on $\chi$.

**Proof:** We define $\iota_\chi$ by $\iota_2^{-1} \circ \iota_1$. Let us check the (twisted) equivariance for the central action. Recall that

$$eH^{\cdot}(\Gamma_0(p), \mathcal{C}(\rho^1, A)) = \varinjlim eH^{\cdot}(\Gamma_0(p^r), Ind_{\Gamma_1(p^r)}^{\Gamma_0(p^r)} V \otimes A_r))$$

and

$$e\mathcal{W}_\chi = \varinjlim eH^{\cdot}(\Gamma_1(p^r), V)$$

and that $\iota_\chi^{-1}$ is induced by the canonical isomorphism

$$eH^{\cdot}(\Gamma_0(p^r), Ind_{\Gamma_1(p^r)}^{\Gamma_0(p^r)} V \otimes A_r)) \to eH^{\cdot}(\Gamma_1(p^r), V)$$

Now, let us consider an element $z \in Z_{M'}(\mathbf{Z}_p)$; let $\zeta = i(z)$ and $\gamma_\zeta$ be the corresponding element of $\Gamma_0(p^r)/\Gamma_1(p^r)$; note that one can assume that $\gamma_\zeta \equiv z$ mod. $p^r$. For $[u] \in eH^q_{Q-no}(\Gamma_0(p^r), \mathcal{C}(\rho^1; A))$, one knows that the cocycle

$v : (\gamma_0, \ldots, \gamma_q) \mapsto (y \in M'(\mathbf{Z}/p^r\mathbf{Z}) \mapsto u(\gamma_\zeta^{-1}\gamma_0\gamma_\zeta, \ldots, \gamma_\zeta^{-1}\gamma_q\gamma_\zeta)(\gamma_\zeta^{-1}y))$ is cohomologous to $u$. On the other hand, since $z \in Z_{M'}$, one has $\gamma_\zeta^{-1} y = yz^{-1}$; hence

$$v(\gamma_0, \ldots, \gamma_q) = (z.u)(\gamma_\zeta^{-1}\gamma_0\gamma_\zeta, ..., \gamma_\zeta^{-1}\gamma_q\gamma_\zeta);$$

thus, by applying $\iota_\chi^{-1}$ to the classes of these cocycles, one sees that

$$\omega_\chi(z)^{-1} < \zeta >_{\rho,\chi} \iota_\chi^{-1}(z.[u]) = \iota_\chi^{-1}([u])$$

which yields 3.2.

### 3.4 Weak Control theorem

Since the $J$-proper SPS $Q$ contains the Borel subgroup $B_H$, its Levi subgroup $M$ contains the standard maximal torus $T_H = Res_{\mathbf{Z}}^{O_F}T$ of $H$. Recall that we have fixed an absolutely irreducible representation $\rho$ of $M$ defined over $\mathcal{O}$; we denote by $\delta \in X^*(T_H) = X^*(T)^{S_p}$ the dominant character associated with the fixed $H$-compatible representation $\rho$ (see Definition 2.1).

**Definition 3.1** *A character $\chi \in X^*(T_H) = X^*(T)^{S_p}$ is called regular dominant, if for all $v \in S_p$, one has $\chi_v = (a_1, a_2; b)$ with $b \geq a_1$ and $a_1 > a_2 > 0$.*

For any $\alpha \geq 1$, let

$$L_\alpha^a(\rho; \mathcal{O}) = \{f : Y_{Q',\alpha}^a(\mathbf{Z}_p) \to V; \ f \text{ polynomial,}$$

$$f(ym) = \rho(m)^{-1}f(y) \text{ for } m \in M'(\mathbf{Z}_p)\}$$

resp.

$$\mathcal{C}_\alpha(\rho^1; \mathcal{O}) = \{f : Y_{Q',\alpha}^a(\mathbf{Z}_p) \to V; \ f \text{ continuous,}$$

$$f(ym) = \rho(m)^{-1}f(y) \text{ for } m \in M^1(\mathbf{Z}_p)\}$$

These are left $I'_\alpha$-modules for $g.f(y) = f(g^{-1}y)$. Note that $\mathcal{C}_\alpha(\rho^1; A) = \mathcal{C}_\alpha(\rho^1; \mathcal{O}) \otimes A$ can also be defined as the space of locally constant functions with values in $A$, by density of the space of $V$-valued locally constant functions in the space of $V$-valued continuous functions.

Let us define also twisted versions of $L_\alpha^a(\rho; \mathcal{O})$. For that purpose we introduce arithmetic characters.

**Definition 3.2** *An arithmetic character of $C_{M'}$ is a product $\varepsilon\chi$ where $\varepsilon : C_{M'}(\mathbf{Z}_p) \to \mathcal{O}^\times$ is a finite order character and $\chi$ is an algebraic character of $C_{M'}$. a level for $\varepsilon\chi$ is a p-power $p^r$ such that $\varepsilon$ factors through $C_{M'}$. $\varepsilon\chi$ is said dominant with respect to $\rho$ (in brief, $\rho$-dominant), resp. regular with respect to $\rho$, if $\delta \otimes \chi$ is dominant, resp. regular.*

Note that if the character $\delta$ associated with $\rho$ is regular, then any arithmetic character $\varepsilon\chi$ of $C_{M'}$ dominant with respect to $\rho$, is regular with respect to $\rho$.

For any $\rho$-dominant arithmetic character $\varepsilon\chi$, let $\alpha \geq 1$ such that $\varepsilon$ factors through $C_{M'}(\mathbf{Z}/p^\alpha\mathbf{Z})$; we write $L_\alpha^a(\rho \otimes \varepsilon\chi; \mathcal{O})$ for $L_\alpha^a(\rho \otimes \chi; \mathcal{O})$ viewed as an $I'_\alpha$-module for the following twisted action: $(g.f)(y) = \varepsilon(g)f(g^{-1}y)$, where $\varepsilon(g)$ means that one applies $\varepsilon$ to the element of
$$C_{M'}(\mathbf{Z}/p^\alpha\mathbf{Z}) = Q'(\mathbf{Z}/p^\alpha\mathbf{Z})/\tilde{Q}(\mathbf{Z}/p^\alpha\mathbf{Z})$$
congruent to $g \bmod p^\alpha \in Q'(\mathbf{Z}/p^\alpha\mathbf{Z})$ modulo $\tilde{Q}(\mathbf{Z}/p^\alpha\mathbf{Z})$.

Let $\omega_{\varepsilon\chi}$ be the character of $Z_{M'}(\mathbf{Z}_p)$ given by $\omega \times (\varepsilon\chi \circ i)$. As already noticed, there is an action of $Z_{M'}(\mathbf{Z}_p)$ on $\mathcal{C}_\alpha(\rho^1; A)$, given by $(z.f)(g) = f(gz^{-1})$; hence one can speak of the largest submodule $\mathcal{C}_\alpha(\rho; A)[\omega_{\varepsilon\chi}]$ on which $Z_{M'}(\mathbf{Z}_p)$ acts via $\omega_{\varepsilon\chi}$. One can now state the main theorem of this Section. A linear map between $\mathcal{O}$-modules is called an isogeny if its kernel and cokernel are finite.

**Theorem 3.2** *For any $\rho$-dominant arithmetic character $\varepsilon\chi$ regular with respect to $\rho$, for any $\alpha$ such that $p^\alpha$ is a level of $\varepsilon\chi$, one has*

*(i) For any $0 \leq q < 3d$, $H^q(\Gamma_0(p^\alpha), L_\alpha^a(\rho \otimes \varepsilon\chi; K)) = 0$,*

*(ii) for any $q$, there are natural isomorphisms $\iota_{\varepsilon\chi}^q$:*

$$\iota_{\varepsilon\chi}^q : H^q_{Q-no}(\Gamma_1(p^\infty), L_\alpha^a(\rho \otimes \varepsilon\chi; A)) \cong H^q_{Q-no}(\Gamma_0(p^\alpha), \mathcal{C}_\alpha(\rho^1; A))$$

*defined in a way similar to 3.4; with the Notation of Corollary 2, they satisfy*

$$z.\iota_{\varepsilon\chi}^q(c) = \omega_{\varepsilon\chi}(z)\iota_{\varepsilon\chi}^q(<i(z)>_{\rho,\varepsilon\chi}^{-1} c)$$

*(iii) for $q = 3d$, $\iota_{\varepsilon\chi}^q$ induces, for any $r \geq \alpha$, an isogeny:*

$$H^{3d}_{Q-no}(\Gamma_0(p^r), L_\alpha^a(\rho \otimes \varepsilon\chi; A)) \to H^{3d}_{Q-no}(\Gamma_0(p^\alpha), \mathcal{C}_\alpha(\rho^1; A))[\omega_{\varepsilon\chi}].$$

**Proof:** (i) This statement is a corollary of a result by J. Franke [8]. The calculation is detailed in Appendix A below. The idea is the following: for any congruence subgroup $\Phi$ of $G'(F)$, Franke defines a spectral sequence abutting to the cohomology groups $H^q(\Phi, L(\rho \otimes \psi; \mathbf{C}))$ (Corollary 4.8 of [45] and Proof of Th.18 of Franke's paper). Its $E_1$-term is expressed in terms of the $(\text{Lie}(M), K_M)$-cohomology of unitary automorphic representations of $M$. Then, if for each archimedean place $v$ of $F$ the highest weight $\delta_v\psi_v$ of $L(\rho_v \otimes \psi_v; \mathbf{C})$ is regular, it follows from the classification of Vogan-Zuckerman ([44] that $E_1^{s,t} = 0$ if $s + t < 3d$. In fact, one only needs the explicitation of Vogan-Zuckerman calculations for $Sp_4$ and $SL_2$. It is done in the $Sp_4$-case by R. Taylor (p.293 of [34]), for $SL_2$ it is a classical result that the cohomology of an infinite dimensional unitary representation can be non-zero only in degree 1. We first found these vanishing results in [27], for $F = \mathbf{Q}$ where the author assumed that the weight of the local system is sufficiently regular; they were confirmed to us for an arbitrary totally real field $F$ by J. Schwermer in a letter [28]. We thank D. Blasius for pointing out to us the reference of Franke [8].

We shall deduce *(iii)* from *(i)* and *(ii)* by using a control criterion due to Hida (Lemma 5.1 of [17]). First, we perform several reductions.

1) For any $\mathbf{Z}_p$-algebra $A$, let
$$M_A^z = Z_{M'}(A)M^1(A);$$

it is a subgroup of $M'(A)$ of exponent two (it is trivial if $Q = B_H$). Let $\Gamma'_0(p^r)$ be the subgroup of $\Gamma_0(p^r)$ of elements congruent to an element of $M^z_{\mathbf{Z}/p^r\mathbf{Z}}Q^+(\mathbf{Z}/p^r\mathbf{Z})$ modulo $p^r$. Consider
$$Ind_{M^z_{\mathbf{Z}/p^r\mathbf{Z}}}^{M'(\mathbf{Z}/p^r\mathbf{Z})} V \otimes \varepsilon\chi \otimes A_r = \{f : M'(\mathbf{Z}/p^r\mathbf{Z}) \to V \otimes A_r;$$
$$f(gh^{-1}) = \rho \otimes \varepsilon\chi(h)f(g) \text{ for all } h \in M^z_{\mathbf{Z}/p^r\mathbf{Z}}\}$$

Shapiro's lemma yields:
$$H^q_{Q-no}(\Gamma_0(p^r), Ind_{M^z_{\mathbf{Z}/p^r\mathbf{Z}}}^{M'(\mathbf{Z}/p^r\mathbf{Z})} V \otimes \varepsilon\chi \otimes A_r) \cong H^q_{Q-no}(\Gamma'_0(p^r), V \otimes \varepsilon\chi \otimes A_r)$$

For any $r \geq \alpha$, let $f : Y^a_{Q',\alpha}(\mathbf{Z}/p^r\mathbf{Z}) \to V \otimes A_r$ be a function such that $f(xm^{-1}) = \rho(m)f(x)$ for $m \in M^1(\mathbf{Z}/p^r\mathbf{Z})$ and $f(xz^{-1}) = \omega_{\varepsilon\chi}(z)f(x)$ for $z \in Z_{M'}(\mathbf{Z}/p^r\mathbf{Z})$; its restriction to $M'(\mathbf{Z}/p^r\mathbf{Z})$ has a prescribed behavior on $M^z_{\mathbf{Z}/p^r\mathbf{Z}}$, hence belongs to $Ind_{M^z_{\mathbf{Z}/p^r\mathbf{Z}}}^{M'(\mathbf{Z}/p^r\mathbf{Z})} V \otimes \varepsilon\chi \otimes A_r$. We consider a variant of 3.2; namely, for each $r \geq \alpha$ the exact sequence of $\Delta^{-1}_{Q,r}$-modules:

$$0 \to \mathcal{K}'_r \to \mathcal{C}(Y^a_{Q',\alpha}(\mathbf{Z}/p^r\mathbf{Z}), \rho^1; A_r)[\omega_{\varepsilon\chi}] \xrightarrow{\psi_r}$$
$$\to Ind_{M^z_{\mathbf{Z}/p^r\mathbf{Z}}}^{M'(\mathbf{Z}/p^r\mathbf{Z})} V \otimes \varepsilon\chi \otimes A_r \to 0$$

where the map $\psi_r$ is the restriction from $Y^a_{Q',\alpha}(\mathbf{Z}/p^r\mathbf{Z})$ to $M'(\mathbf{Z}/p^r\mathbf{Z})$.

2) We take the long exact sequence of cohomology for $\Gamma_0(p^r)$; Again, one checks that for any $r \geq \alpha$,

1. $H^q_{Q-no}(\Gamma'_0(p^r), \mathcal{K}'_r) = 0$

2. for any $\Delta^{-1}_{Q,\alpha}$-module $L$, for any $r \geq \alpha$,
$$Res^{\Gamma_0(p^\alpha)}_{\Gamma_0(p^r)} : H^q_{Q-no}(\Gamma_0(p^r), L) \cong H^q_{Q-no}(\Gamma_0(p^\alpha), L)$$

Therefore, we obtain an isomorphism similar to $\iota_2$ of 3.4 (see Theorem 3.1 above):

$$H^q_{Q-no}(\Gamma_0(p^\alpha), \mathcal{C}(Y^a_{Q',\alpha}(\mathbf{Z}/p^r\mathbf{Z}); \rho^1; A_r)[\omega_{\varepsilon\chi}]) \stackrel{\iota'_2}{\cong} H^q_{Q-no}(\Gamma'_0(p^r), V \otimes \varepsilon\chi \otimes A_r) \quad (3.5)$$

for each $r \geq \alpha$.

3) Let $Q^z(A) = M_{\bar{A}}^z Q^+(A)$; for any $s \geq \alpha$ such that

$$L_\alpha^a(\rho \otimes \varepsilon\chi; A_r) = L(Y_{Q',\alpha}^a(\mathbf{Z}/p^s\mathbf{Z}), \rho \otimes \varepsilon\chi; A_r).$$

Consider the exact sequence of $\Delta_{Q^z,s}^{-1}$-modules (for $s \geq r$):

$$0 \to K_r' \to L_\alpha^a(\rho \otimes \varepsilon\chi; A_r) \xrightarrow{\psi_r} V \otimes \varepsilon\chi \otimes A_r \to 0 \tag{3.6}$$

We can take the cohomology of 3.6 for $\Gamma_1(p^s)$ or for $\Gamma_0'(p^s)$; with the first choice, one obtains for any $s \geq r$:

$$\iota_1' : H_{Q-no}^q(\Gamma_1(p^s), L_\alpha^a(\rho \otimes \varepsilon\chi; A_r)) \cong H_{Q-no}^q(\Gamma_1(p^s), V \otimes \varepsilon\chi \otimes A_r)$$

and by taking the direct limit over $s \geq r$:

$$\iota_1' : H_{Q-no}^q(\Gamma_1(p^\infty), L_\alpha^a(\rho \otimes \varepsilon\chi; A)) \cong H_{Q-no}^q(\Gamma_1(p^\infty), V \otimes \varepsilon\chi \otimes A) \tag{3.7}$$

With the second choice, one obtains for any $s \geq r$:

$$\iota_1'' : H_{Q-no}^q(\Gamma_0'(p^s), L_\alpha^a(\rho \otimes \varepsilon\chi; A_r)) \cong H_{Q-no}^q(\Gamma_0'(p^s), V \otimes \varepsilon\chi \otimes A_r) \tag{3.8}$$

Comparing 3.4 (in Th. 3.1) and 3.7, we obtain the isomorphism $\iota_{\varepsilon\chi}^q = \iota_2^{-1} \circ \iota_1'$:

$$\iota_{\varepsilon\chi}^q : H_{Q-no}^q(\Gamma_1(p^\infty), L_\alpha^a(\rho \otimes \varepsilon\chi; A)) \cong H_{Q-no}^q(\Gamma_0(p^\alpha), \mathcal{C}_\alpha(\rho^1; A))$$

This proves (ii). Moreover, by comparing 3.8 and 3.5, we obtain an isomorphism $\tilde{\iota}_{\varepsilon\chi}^q = \iota_2'^{-1} \circ \iota_1''$:

$$\tilde{\iota}_{\varepsilon\chi}^q : H_{Q-no}^q(\Gamma_0'(p^\alpha), L_\alpha^a(\rho \otimes \varepsilon\chi; A_r)) \cong H_{Q-no}^q(\Gamma_0(p^\alpha), \mathcal{C}_\alpha(\rho^1; A)[\omega_{\varepsilon\chi}]) \tag{3.9}$$

Exactly as in Corollary 3.2, one can show

$$z.\iota_{\varepsilon\chi}^q(c) = \omega_{\varepsilon\chi}(z)\iota_{\varepsilon\chi}(< \zeta^{-1} >_{\rho,\varepsilon\chi} c)$$

4) Next, we shall relate the two isomorphisms $\iota_{\varepsilon\chi}^q$ and $\tilde{\iota}_{\varepsilon\chi}^q$ for $q = 3d$ by using the natural map:

$$H_{Q-no}^q(\Gamma_0(p^\alpha), \mathcal{C}_\alpha(\rho^1; A)[\omega_{\varepsilon\chi}]) \to H_{Q-no}^q(\Gamma_0(p^\alpha), \mathcal{C}_\alpha(\rho^1; A))[\omega_{\varepsilon\chi}] \tag{3.10}$$

In order to show that for $q = 3d$ the map 3.10 is an isogeny, we apply Hida's control criterion (Lemma 5.1 of [17]). Namely, we consider the completed $\mathcal{O}$-algebra $\mathbf{\Lambda} = \mathcal{O}[[Z_{M'}(\mathbf{Z}_p)]]$ of the group $Z_{M'}(\mathbf{Z}_p)$. This last group is the direct product of a finite group $Z_0$ by a free $\mathbf{Z}_p$-module $\Phi$ of rank $m = \sum_{v \in J} d_v \dot{r}k(Z_{M'_v})$. Let $\Lambda = \mathcal{O}[[\Phi]]$ be the completed group $\mathcal{O}$-algebra of $\Phi$. Let $\mathbf{P}_{\varepsilon\chi}$ be the ideal of $\mathbf{\Lambda}$, kernel of $[z] \mapsto \omega_{\varepsilon\chi}(z)$; since $\Lambda$ is a regular local ring, $\mathbf{P}_{\varepsilon\chi} \cap \Lambda$ is generated by a regular sequence $(T_1, \ldots, T_m)$. Let $\mathcal{B}$, resp. $\mathcal{A}$, be the category of discrete $\Lambda$-modules, resp. of discrete $\mathcal{O}$-modules with left action of $\Delta_{Q,\alpha}^{-1}$, such that $I'_\alpha$ acts smoothly admissibly (that is, such that any point has an open stabilizer and invariant points by any open subgroup form a cofinite type submodule). Note that $\mathcal{A}$ is a subcategory of $\mathcal{B}$. Consider the cohomological functor

$$\mathcal{H}^{\cdot} : \mathcal{A} \mapsto \mathcal{B}, \quad N \mapsto H^{\cdot}_{\dot{Q}-no}(\Gamma_0(p^\alpha), N).$$

Let $\omega_0$, resp. $\omega_\Phi$ be the restriction of $\omega_{\varepsilon\chi}$ to $Z_0$ resp. to $\Phi$. Let $\mathcal{C}_0 = \mathcal{C}_\alpha(\rho^1; A)[\omega_0]$. Hida's control criterion says that the map 3.10 is an isogeny for $q = 3d$, provided the following four hypotheses are verified.

1. ($H_1$) The modules

    $$\mathcal{C}_\alpha(\rho^1; A), \quad ,\mathcal{C}_0[T_1, \ldots, T_j] \quad (j = 1, \ldots, m), \ \mathcal{C}_\alpha(\rho^1; A)[\omega_{\varepsilon\chi}] = \mathcal{C}_0[T_1, \ldots, T_m]$$

    are objects of $\mathcal{A}$,

2. ($H_2$) for each $j = 1, \ldots, m$, the endomorphism of multiplication by $T_j$ on $\mathcal{C}_0[T_1, \ldots, T_{j-1}]$ is surjective,

3. ($H_3$) $\mathcal{H}^q(\mathcal{C}_0)$ is finite for all $q < 3d$,

4. ($H_4$) $\mathcal{H}^{3d}(\mathcal{C}_\alpha(\rho^1; A)[\omega_{\varepsilon\chi}])$ is of cofinite type over $\mathcal{O}$ (that is, $\mathrm{Hom}_\mathcal{O}(-, A)$ sends it to a finitely generated $\mathcal{O}$-module).

Condition ($H_1$) is obvious; for ($H_2$), one writes $Y^a_{Q',\alpha}(\mathbf{Z}_p) = Q^-(p^\alpha \mathbf{Z}_p) \times M'(\mathbf{Z}_p)$, hence the Pontryagin dual $\mathcal{C}_\alpha(\rho^1; A)^* = \mathrm{Hom}_\mathcal{O}(\mathcal{C}_\alpha(\rho^1; A), A)$ can be identified to

$$\mathcal{O}[[Q^-(p^\alpha \mathbf{Z}_p)]] \hat{\otimes}_\mathcal{O} \mathcal{O}[[C_{M'}(\mathbf{Z}_p)]] \hat{\otimes}_\mathcal{O} \mathrm{Hom}_\mathcal{O}(V, \mathcal{O})$$

Since $\mathcal{O}[[C_{M'}(\mathbf{Z}_p)]] \cong \Lambda^k$, the multiplication by $T_j$ on $\mathcal{C}_0^*/(T_1, \ldots, T_{j-1})\mathcal{C}_0^*$ is injective; hence the conclusion by duality. For ($H_3$), we simply apply (i) (for $\chi$ regular with respect to $\rho$). The last condition ($H_4$) is verified by 3.9. ∎

## 3.5 Descent to prime-to-$p$ level

We consider in this section a proper standard parabolic subgroup $Q$ of $H$ defined over $\mathbf{Z}$; that is, $Q$ is the restriction of scalars from $\mathbf{r}$ to $\mathbf{Z}$ of $B$, $P$ or $P^*$. In particular, its base change to $\mathbf{Z}_p$ is an $S_p$-proper SPS for any prime number $p$. As usual, we write the Levi decomposition (over $\mathbf{Z}$) of $Q$ as $MQ^+$. Let $K_0 \subset \bar{\mathbf{Q}}$ be a number field containing all $\sigma(F)$ when $\sigma$ runs in $I_F$ and let $\mathcal{O}_0$ be its ring of integers; Fix a finite free $\mathcal{O}_0$-module $V_0$ and a representation $\rho_0 : M_{/\mathcal{O}_F} \to GL_{\mathcal{O}_0}(V_0)$ defined over $\mathcal{O}_F$. For any prime $p$, we fix a $p$-adic embedding $\iota_p$ of $\bar{\mathbf{Q}}$ in $\bar{\mathbf{Q}}_p$ and we denote by $\mathcal{O}$, resp. $K$, the $p$-adic completion of $\iota_p(\mathcal{O}_0)$, resp. of $\iota_p(K_0)$, and we put $V = V_0 \otimes_{\mathcal{O}_0} \mathcal{O}$, $\rho = \rho_0 \otimes \mathcal{O}$; we denote by $\pi$ a uniformizing parameter of $\mathcal{O}$.

Let

$$L_0(\rho_0; \mathcal{O}_0) = \{f : H^1_{\mathbf{Z}}/Q^+_{\mathbf{Z}} \to V_0;\ f \text{ regular, and } f(ym^{-1}) = \rho_0(m)f(y)\}$$

It is a lattice in $L_0(\rho; K_0)$, defined independently of $p$. Once $p$ is fixed, one can write $V = \bigotimes_{v \in J} V_v \otimes \bigotimes_{v \notin J} V_v$ and accordingly:

$$L_0(\rho; \mathcal{O}) = \bigotimes_{v \in J} L_0(\rho_v; \mathcal{O}) \otimes \bigotimes_{v \notin J} V_v$$

For each place $v$ of $F$ dividing $p$, let us consider the semi-group $\Delta_{v,0} = K'_v D_v K'_v$; where $D_v$ is defined in terms of the center $Z_v$ of $M_v = M_{Q_v}$ as in section 2.2, but with the hyperspecial subgroup $K'_v = G'(\mathbf{r}_v)$ instead of the parahoric subgroup $I'_v$. Recall that we denote by $\omega_v$ the central character of $\rho_v$, that is, the restriction of $\rho$ to the center $Z_{M'_v}$ of $M'_v$. We let $D_v^{-1}$ act on $L(\rho_v; K)$ by $d^{-1} \mapsto \omega_v(d^{-1}).\rho_v(d)$. Then we extend this as an action of all $\Delta_{v,0}^{-1}$ as in 2.3. Let $\Delta_0 = \prod_{v|p} \Delta_{v,0}$.

Let $e_0 = \lim_{n \to \infty} (T_Q^0)^{n!}$ be the idempotent associated to the Hecke operator $T_Q^0 = [\Gamma \xi_Q \Gamma]$ where $\xi_Q$ is defined as in Section 2.4.

A level subgroup $U \subset G(\hat{\mathbf{Z}})$ of level prime to $p$ is fixed. Let $\Gamma = UG_\infty \cap G^1(\mathbf{Q})$.

We want to see that if $\rho$ has regular weight, one can get rid of the $p$-part of the level in $H^q_{Q-no}(\Gamma_0(p), L^a(\rho; A))$ (this is the analogue of Lemma 7.2 of [17]). This fact will be of great importance in the exact control theorem in next section, as well as in the algebro-geometric considerations in Section 7 on the Galois representations associated to $p$-ramified cuspidal representations in the regular discrete series.

**Proposition 3.2** *Assume that the highest weight of $\rho$ is regular dominant for the ordering defined by $(M', B \cap M', T \cap M')$ and that $p$ is such that no completion $F_v$ for $v|p$ contains a $p$-th root of unity:*

(TF)  $\mu_p(F_v) = \{1\}$ *for $v|p$,*

*then*

(i) the lattice $L_0(\rho; \mathcal{O})$ is stable by the action of $\Delta_0$ described above,

(ii) for any $q \geq 0$, there is a natural isomorphism

$$H^q_{Q-no}(\Gamma_0(p), L(\rho; A)) \cong e_0 H^q(\Gamma, L_0(\rho; A))$$

**Comments:** 1) The first assumption is of course implied by the regularity of the highest weight of $\rho$ ( by the $H$-compatibility of $\rho$, see Def.2.1); it is strictly weaker if $Q_v$ is maximal at some place $v|p$.
2) the assumption $(TF)$ is fulfilled if $p$ does not ramify in $F$. $(TF)$ will be used to insure that $Z_{M'_v}(F_v)$ has no $p$-torsion.
3) By statement (i), the action of $T^0_Q$ on $H^q(\Gamma, L_0(\rho; A))$ is well defined, hence one can speak of $e_0 H^q(\Gamma, L_0(\rho; A))$.

**Proof:** For assertion (i), one needs to see that for any $d \in D_v$, $\omega_v(d^{-1}).\rho_v(d)$ preserves $L_0(\rho_0; \mathcal{O})$. For that, we observe first that by $(TF)$ and by Sect.2.5 and 2.11 of [22], we have a direct sum decomposition

$$L_0(\rho; \mathcal{O}) = \bigoplus_\eta V_\eta$$

where $\eta$ runs over weights of $Z_{M'_v}(\mathbf{r}_v)$ acting on $L_0(\rho; K)$ and $V_\eta$ is the eigenmodule corresponding to the eigenvalue $\eta$. For each $v \in J$, the highest weight of $\rho_v$ is $\delta_v$, therefore, if $\alpha$ and $\beta$ are respectively the short and long simple roots of $Sp_4$, one has for any weight $\theta$ of $T'(\mathbf{r}_v)$ acting on $L_0(\rho; \mathcal{O})$, $\delta_v \theta^{-1} = \alpha^{m_\eta} \beta^{n_\eta}$ with $m_\theta, n_\theta \geq 0$. Hence, for any $d \in D_v$,

$$ord_v(\theta \delta_v^{-1}(d)) \geq 0$$

Now, we remark that $\delta_v$ restricted to $Z_{M'_v}$ is equal to $\omega_v$, and that the map $\theta \mapsto \eta$ given by restriction from $T'$ to $Z_{M'_v}$ is surjective; we conclude that for any $d \in D_v$,

$$ord_v(\eta \omega_v^{-1}(d)) \geq 0$$

This proves (i).

For proving (ii), we need first a fact.

For $v \in J$ and $Q_v$ the parabolic subgroup determined by $Q$, $I'_v$ the corresponding parahoric subgroup; let $d_v \in D_v$ be the element associated to $Q_v$ as in Proposition 2.2; let $W_v$ be the Weyl group of $(G_v, B_v, T_v)$ resp. $W_{M_v}$ the Weyl group of $Z_{M_v}$ and $\tau$ an element of order two in $G'$ inducing the element of greatest length in $W_v$; actually one can take $\tau = \begin{pmatrix} 0 & 1_2 \\ -1_2 & 0 \end{pmatrix}$.
Let $d_v = d_{v,i}$ for $Q_v = P_i$ as in Sect.2.3; let us consider the group $J^o_v = K'_v \cap d_v^{-1} K'_v d_v$ and $J_v = \tau^{-1} J^o_v \tau$. We have

**Lemma 3.2**
$$K'_v d_v K'_v = \coprod_{w \in W_v/W_{M_v}} \coprod_{u \in ((\tau w)^{-1} J'_v \tau w \cap I'_v) \backslash I'_v} K'_v d_v wu$$

The computation is done in [33] p.316-318 (see also [17] p.467). Conventions there are slightly different, so we redo the calculation. We have

$$d_v^{-1} K'_v d_v \backslash d_v^{-1} K'_v d_v K'_v = \coprod_\varepsilon d_v^{-1} K'_v d_v \varepsilon, \quad \text{where } J_v^o \backslash K'_v = \coprod_\varepsilon J_v^o \varepsilon$$

Note that $J'_v$ contains a parahoric subgroup $I'_r$ (actually $r = 1, 2$ or $3$ according whether $Q'_v$ is the Siegel, Klingen or Borel parabolic). Hence, one has by the Bruhat decomposition:

$$K'_v = \coprod_{w \in W_v/W_{M_v}} Q'_v w I'_1 = \coprod_{w \in W_v/W_{M_v}} J'_v w I'_1 = \coprod_{w \in W_v/W_{M_v}} \coprod_{u \in w^{-1} J'_v w \cap I'_v) \backslash I'_v} J'_v wu$$

One can and will assume that the representatives $u$ of $w^{-1} J'_v w \cap I'_v) \backslash I'_v$ are in the unipotent radical $Q_v^+$ of $Q'_v$. This decomposition yields:

$$K'_v = \coprod_{w \in W_v/W_{M_v}} \coprod_{u \in (w^{-1} J'_v w \cap I'_v) \backslash I'_v} J_v^o \tau wu$$

hence,

$$d_v^{-1} K'_v d_v K'_v = \coprod_{w \in W_v/W_{M_v}} \coprod_{u \in w^{-1} J'_v w \cap I'_v) \backslash I'_v} d_v^{-1} K'_v d_v \tau wu$$

or, in conclusion:

$$K'_v d_v K'_v = \coprod_{w \in W_v/W_{M_v}} \coprod_{u \in ((\tau w)^{-1} J'_v \tau w \cap I'_v) \backslash I'_v} K'_v d_v wu$$

Now we can start the proof of the proposition. Recall that we denote by $\pi$ a uniformizing parameter of $\mathcal{O}$. For any $r \geq 1$, put $A_r = \pi^{-r}\mathcal{O}/\mathcal{O}$ and $A = K/\mathcal{O}$. We consider the homomorphism

$$\iota : L_0(\rho; A) \to L^a(\rho; A)$$

obtained from the inclusion of lattices $L_0(\rho; \mathcal{O}) \subset L^a(\rho; \mathcal{O})$ by tensoring with $A$. We shall study the homomorphism

$$j : e_0 H^q(\Gamma, L_0(\rho; A)) \to e H^q(\Gamma_0(p), L^a(\rho; A))$$

defined by $j = \iota_* \circ \ell$ where

$$\ell = e \circ \text{res} : e_0 H^q(\Gamma, L_0(\rho; A)) \to e H^q(\Gamma_0(p), L_0(\rho; A))$$

Let us show that the restrictions $\ell(\pi^r)$ to the $\pi^r$-torsion of these modules are isomorphisms. We start with $r = 1$. Consider the commutative square :

$$\begin{array}{ccc} e_0 H^q(\Gamma, L_0(\rho; A_1)) & \xrightarrow{\ell_1} & eH^q(\Gamma_0(p), L_0(\rho; A_1)) \\ \downarrow & & \downarrow \\ e_0 H^q(\Gamma, L_0(\rho; A))[\pi] & \xrightarrow{\ell(\pi)} & eH^q(\Gamma_0(p), L_0(\rho; A))[\pi] \end{array} \quad (3.11)$$

where $\ell_1$ is defined just like $\ell$, replacing $A$ by $A_1$ and where the vertical arrows are the natural surjective homomorphisms. We first show that $\ell_1$ is an isomorphism. For that purpose, let us consider the evaluation map

$$\mathbf{p}_0 : L_0(\rho, \psi; A_1) \to (V \otimes A_1), \quad f \mapsto f(Q^+)$$

for any polynomial function $f : H'/Q^+(\mathbf{Z}_p) \to V \otimes A_1$.

Note that $\mathbf{p}_0$ is equivariant for the action of $Z_{M'}(\mathbf{Z}/(p))$: $(z.f)(1) = f(z^{-1}) = \rho(z)f(1) = \omega(z)f(1)$. Similarly, the evaluation map at $wQ^+$ for $w \in W/W_M = \prod_{v \in J} W_v/W_{M_v}$ satisfies $(z.f)(w) = f(z^{-1}w) = \rho(wzw^{-1})f(w) = \omega^w(z)f(w)$. Thus, for a cocycle $c$ with values in $L_0(\rho; \mathcal{O})$, one has

$$\begin{cases} (d_v wu)^{-1} c(*, \ldots, *; Q^+) \equiv 0 \text{ mod. } \omega(d_v) & \text{if } w \in W_M \\ (d_v wu)^{-1} c(*, \ldots, *; Q^+) \equiv 0 \text{ mod. } \omega(d_v)\pi & \text{if } w \notin W_M \end{cases}$$

Hence, by 3.2 in the previous Lemma, we see that $T_Q^0$ coincides with the level $p$ Hecke operator $T_Q$ modulo $\pi$; thus, $e \circ \text{res} = \text{res} \circ e_0$. Let us consider the commutative diagram

$$\begin{array}{ccc} e_0 H^q(\Gamma, L_0(\rho; A_1)) & \xrightarrow{\ell_1} & eH^q(\Gamma_0(p), L_0(\rho; A_1)) \\ k_1 \downarrow & & \downarrow k_1' \\ eH^q(\Gamma_0(p), V \otimes A_1) & = & eH^q(\Gamma_0(p), V \otimes A_1) \end{array}$$

The proof of [12], Th.3.2 (second statement) proves similarly that $k_1 \circ e_0 = e \circ k_1$ and that $k_1$ and $k_1'$ are isomorphisms. Hence, $\ell_1$ is an isomorphism. Note that in the case of $k_1'$ the contraction Proposition 2.2 applies; however, for $k_1$ it does not apply because the strata of $K_v'/Q_v^+(\mathbf{r}_v)$ apart from the open one don't contract. Therefore, to prove the surjectivity of $k_1$, one requires instead to produce a section thereof (up to automorphism); it is provided by $e_0 \circ \text{Tr} \circ inj$ where $inj$ is the canonical injection of $V \otimes A_1 \to L_0(\rho; A_1)$ given by the highest weight, twisted by $\tau$, element of greatest length in the Weyl group, viewed in $\Gamma$ by weak approximation (see p.233-236 of [12]).

Note that the injectivity of $\ell_1$ is easy: $\text{Tr} \circ \text{res}$ is the endomorphism of multiplication by $(\Gamma : \Gamma_0(p))$ on $e_0 H^q(\Gamma, L_0(\rho; A_1))$, and this index is prime to $p$ (provided $p$ does not ramify in $F$; if it does, one should replace in the statement of the proposition $\Gamma_0(p)$ by $\Gamma_0(\mathfrak{p})$ where $\mathfrak{p}$ is the product of the prime ideals in $J$).

Now we can see that $\ell(\pi)$ itself is an isomorphism. The quasi-inverse of $\ell_1$ is

$$e_0 \circ \text{Tr} : eH^q(\Gamma_0(p), L_0(\rho; A_1)) \to e_0 H^q(\Gamma, L_0(\rho; A_1))$$

The formula $e_0 \circ \text{Tr} \circ e \circ \text{res} = (\Gamma : \Gamma_0(p))$ is true on $e_0 H^q(\Gamma, L_0(\rho; A_1))$ hence, by functoriality, it is also true on its surjective image $e_0 H^q(\Gamma, L_0(\rho; A))[\pi]$. It shows that $\ell(\pi)$ is injective, but it is surjective by 3.11 so it is bijective with quasi-inverse $\ell'(\pi)$ for $\ell' = e_0 \circ \text{Tr}$. Then, for $r \geq 1$ arbitrary, Nakayama's lemma shows that $\ell(\pi^r)$ is injective as well as $\ell'(\pi^r)$. Hence $\ell'(\pi^r) \circ \ell(\pi^r)$, resp. $\ell \circ \ell'(\pi^r)$, is injective on the finite set $e_0 H^q(\Gamma, L_0(\rho; A))[\pi^r]$ resp. on $eH^q(\Gamma_0(p), L_0(\rho; A))[\pi^r]$, hence these maps are bijective and $\ell(\pi^r)$ is bijective. We conclude that $\ell$ is bijective.

To conclude, we need to see that $\iota_*$ is an isomorphism. For that purpose one considers the commutative triangle:

$$\begin{array}{ccc} L_0(\rho; A_r) & \xrightarrow{\iota_r} & L^a(\rho; A_r) \\ & \searrow \quad \swarrow & \\ & (V \otimes A_r) & \end{array}$$

where the left map is $\mathbf{p}_0$ and the right map is $\mathbf{p}$. Using lemma 3.1 of Section 3.2, we see that the maps induced on cohomology by $\mathbf{p}_0$ and $\mathbf{p}$ are isomorphisms. Then, by the commutativity of the diagram:

$$\begin{array}{ccc} eH^q(\Gamma_0(p), L_0(\rho; A_r)) & \xrightarrow{\iota_{r,*}} & eH^q(\Gamma_0(p), L^a(\rho; A_r)) \\ & \searrow \quad \swarrow & \\ & eH^q(\Gamma_0(p), (V \otimes A_r)) & \end{array}$$

we obtain that $\iota_{r,*}$ is an isomorphism. Finally, taking the inductive limit over $r$, we see that $\iota_*$ is an isomorphism. In conclusion, $j = \iota_* \circ \ell$ is an isomorphism.

## 3.6 Exact Control

We fix the data $(U, \rho)$ as in Section 3.5. We want to prove the following theorem:

**Theorem 3.3** *Let $\rho_0$ as above with regular highest weight, there exists a finite set of primes $S_{U,\rho}$ depending only on the level group $U$ and on $\rho$ such that for any $p \notin S_{U,\rho}$, for any arithmetic character $\chi = \varepsilon\psi : C_{M'}(\mathbf{Z}_p) \to \mathcal{O}^\times$ congruent to 1 modulo $\pi\mathcal{O}$ of level $p^\alpha$, regular with respect to $\rho$, for any integer $q$ such that $0 \leq q < 3d$, one has:*

$$H^q_{Q-no}(\Gamma_0(p^\alpha), L^a(\rho \otimes \varepsilon\psi; A)) \cong H^q_{Q-no}(\Gamma_0(p^\alpha), \mathcal{C}_\alpha(\rho^1; A))[\omega_{\varepsilon\psi}] = 0$$

*and*

$$H^{3d}_{Q-no}(\Gamma_0(p^\alpha), L^a(\rho \otimes \varepsilon\psi; A)) \cong H^{3d}_{Q-no}(\Gamma_0(p^\alpha), \mathcal{C}_\alpha(\rho^1; A))[\omega_{\varepsilon\psi}]$$

**Comment:** The exceptional set $S_{U,\rho}$ can be described as follows: it consists of the set of $p$'s such that either

- $p$ divides the order of the torsion subgroup of $H^q(\Gamma, ind_{M'}^{H'}\rho_0)$ for some $q$ between 1 and $3d$, or

- $p$ divides the order of $\mu(F_p)$,

- $p$ is exceptional with respect to the boundary cohomology (see Section 5, Th.5.7)

For $F = \mathbf{Q}$ the two last conditions are void.

**Proof:** The proof is similar to that of Th.7.1 of [17]. It makes use of *(i)* and *(ii)* of Theorem 3.2 above and of a variant of Hida's control criterion, namely Lemma 7.1 of [17] (proof similar to that of Lemma 5.1 [17]). This variant is the following: one considers the same category of $\Delta_{Q,\alpha}^{-1}$-modules and the same cohomological functor $\mathcal{H}^\cdot = H_{Q-no}^\cdot$ on that category. One considers the same object $\mathcal{C}_\alpha(\rho^1; A)$ and we must verify four hypotheses about this object and its cohomology. These are identical to $(H_1)$ to $(H_4)$ except $(H_3)$ which is replaced by

$(H_3')$ for $p \notin S_{U,\rho}$, for $0 \leq q < 3d$, for all $\varepsilon\psi$ congruent to one modulo $\pi$ $H_{Q-no}^q(\Gamma_0(p^\alpha), L^a(\rho \otimes \varepsilon\psi; A)) = 0$.

In order to define $S_{U,\rho}$ and prove $(H_3')$, we use Proposition 3.2 ("of $p$-destabilisation"). We know that $(H_3')$ is equivalent to

$$H_{Q-no}^q(\Gamma_0(p^\alpha), \mathcal{C}_\alpha(\rho^1; A)[\omega_{\varepsilon\psi}]) = 0$$

But actually, an easy induction argument on $q < 3d$ (yielding exact control for $H_{Q-no}^q(\Gamma_0(p^\alpha), \mathcal{C}(\rho^1; A))$ for any arithmetic character $\omega_{\varepsilon\psi}$ as in the theorem, and for each $q < 3d$) shows that it is enough to check for all $q < 3d$ that

$$H_{Q-no}^q(\Gamma_0(p), \mathcal{C}(\rho^1; A)[\omega]) = 0$$

or in other words, for each $q < 3d$:

$$H_{Q-no}^q(\Gamma_0(p), L^a(\rho; A)) = 0$$

that is, by Proposition 3.2, for $q < 3d$:

$$e_0 H^q(\Gamma, L_0(\rho; K/\mathcal{O})) = 0$$

On one hand, this vanishing is valid for $p \notin S_{U,\rho}$ by Proposition 3.2.

**Comment:** It is however important to notice that one can give another criterion of vanishing for the localization of these groups at the maximal ideal $\mathbf{m}$ of the nearly ordinary

Hecke algebra (see Def.7.1 below) corresponding to the mod. $p$ eigensystem given by $\pi$. Let us say that the maximal ideal $\mathbf{m}$ is "non-Eisenstein" when the residual Galois representation associated to $\pi$ is absolutely irreducible. Then, we have the conjecture

**Conjecture 1** *For $F = \mathbf{Q}$; Assume that the residual Galois representation associated to $\pi$ is absolutely irreducible, then $H^2(\Gamma, L(\rho; K/\mathcal{O}))_{\mathbf{m}} = 0$.*

If one assumes this conjecture, exact control holds for $\mathcal{V}^3_{\mathbf{m}}$.

**Remarks:** (i) Let us assume the conjecture above. Using the Borel-Serre compactification and the corresponding boundary long exact sequence we find that $H^2_c(S(U), L_{\phi^*}(K/\mathcal{O})_{\mathbf{m}}$ is isomorphic to $H^1(\partial(S(U)), L_{\phi^*}(K/\mathcal{O})_{\mathbf{m}}$ which is cofree by Section 5 below (essentially by [16]); using now Poincaré duality, we see that $H^4(S(U), L_{\phi}(\mathcal{O})_{\mathbf{m}}$ is torsion-free, hence $H^3(S(U), L^*_{\phi}(K/\mathcal{O})_{\mathbf{m}}$ is cofree; we deduce from this, not only exact control, but also cofreeness of $e\mathcal{V}^3_{\mathbf{m}}$. In other words, for $F = \mathbf{Q}$, one obtains a better exact control theorem: if $p$ is ordinary and non Eisenstein for $\pi$, the $\Lambda$-module $e\mathcal{V}^3_{\mathbf{m}}$ is cofree.

(ii) Assuming various auxiliary assumptions, this conjecture have been proven recently by Mokrane-Tilouine (cf. [24]) and Urban (cf. [41]) independently and by completely different approaches. The method of Mokrane-Tilouine is quite general and should work for other Shimura varieties provided that an arithmetic compactification theory is available. It relies on the occurence of specific weights in the mod. $p$ crystalline cohomology of the Siegel variety. They need to require that $p$ is prime to the level in $\Gamma$, is bigger than a specific bound depending only on the highest weight of $\rho$ and that the image of $Gal(\overline{\mathbf{Q}}/\mathbf{Q})$ contains the $\mathbb{F}_p$-points of a reductive Chevalley group acting irreducibly on the space of the representation which seems conjecturally true for all but finitely many $p$. The method of Urban is more elementary and necessitates only some regularity condition modulo $p$ of the cohomological weight. However it works only for the case $GSp_{4\mathbf{Q}}$. It rests on the existence of a nice cycle in the Siegel threefold corresponding to the abelian surfaces which are products of two elliptic curves. This cycle had been already used by Weissauer to investigate the classes occuring in the degree two cohomology and coming from the so-called CAP representations.

# 4 p-Ordinary cohomology of boundary strata of the Borel-Serre compactification.

In this section, we generalize the theory of ramification of "cusps" of [16] to the group $GSp_4$ over a totally real number field. The generalization to other groups is straightforward. In fact, we use the fact that $H = GSp_4$ only in the paragraph 4.4.4. In the next section, we will give the complete calculation of the ordinary cohomology of the boundary of the Borel-Serre compactification. In a subsequent paper, we will give a partial generalization of these results for a general reductive $\mathbf{Q}$-group whose derived group is quasi-split at $p$.

## 4.1 Recall on Bruhat order

In this section, $A$ is a valuation ring, $\mathcal{M}_A$ is its maximal ideal and K its field of fractions. We denote by $val_A$ the valuation of A. Let $\mathcal{H}$ be a split and connected semi-simple algebraic $A$-group. We fix $T$ a maximal split torus in $\mathcal{H}$ and denote $W = N_\mathcal{H}(T(A))/T(A)$ the Weyl group associated with $T$. Let $(X, \Phi, X^\vee, \Phi^\vee)$ be a root datum of $(\mathcal{H}, T)$ fixed once for all. If $\alpha \in X$, we denote by $\alpha^\vee \in X^\vee$ the coroot associated to it, and (-,-) the canonical perfect pairing on $X \otimes X^\vee$. Let $R^+ \subset X$ be the corresponding subset of positive roots, $B$ the associated Borel subgroup and $\mathcal{U}_B$ its unipotent radical. For any root $\alpha$ wenote $s_\alpha$ an element of order two in $N_\mathcal{H}(T)$ inducing in the Weyl group the elementary reflexion corresponding to $\alpha$ and $u_\alpha : \mathbb{G}_a \to \mathcal{H}$ the corresponding one parameter subgroup. For any algebra $A$ we have

$$tu_\alpha(a)t^{-1} = u_\alpha(\alpha(t)a) \; \forall t \in T(A) \text{ and } \forall a \in A. \tag{4.1}$$

$$wu_\alpha(\epsilon x)w^{-1} = u_{w(\alpha)}(x) \; \forall x \in A, \forall w \in W \text{ and } \epsilon = \pm 1 \tag{4.2}$$

$$u_\alpha(a^{-1})s_\alpha = u_{-\alpha}(a)u_\alpha(-a^{-1})\alpha^\vee(a^{-1}). \tag{4.3}$$

We consider below the so-called Iwahori subgroup I defined by:

$$I = \{g \in \mathcal{H}^1(A) \text{ such that } g \text{ mod. } \mathcal{M}_A \in B(A/\mathcal{M}_A)\} = \mathcal{U}^-(\mathcal{M}_A)B(A)$$

where $\mathcal{U}^-$ is the unipotent group opposite to $\mathcal{U}_B$.

We think that the two following lemmas are well-known to specialists. However, we give them since we do not know any reference. We are grateful to Choucroun for explaining to us the following proof which is more general and simpler than our first one.

**Lemma 4.1** *Let $\alpha \in \Phi$, $a \in K^*$ and $w \in W$. Then we have:*

$$ws_\alpha u_{-\alpha}(a) \in IwB(K) \cup Iws_\alpha B(K).$$

**Proof** By 4.2, $ws_\alpha u_{-\alpha}(a) = u_{w(\alpha)}(\epsilon a)ws_\alpha$ with $\epsilon = \pm 1$.
If $w(\alpha) > 0$ and $a \in A$ then $ws_\alpha u_{-\alpha}(a) \in \mathcal{U}_B(A)ws_\alpha \subset Iws_\alpha B(K)$.
If $w(\alpha) < 0$ and $a \in \mathcal{M}_A$ then $ws_\alpha u_{-\alpha}(a) \in \mathcal{U}^-(\mathcal{M}_A)ws_\alpha \subset Iws_\alpha B(K)$.
If none of the terms of the alternative above hold, we make use of 4.3:

$$s_\alpha u_{-\alpha}(a) = s_\alpha u_\alpha(a^{-1})s_\alpha \alpha^\vee(a)u_\alpha(a^{-1}) \in s_\alpha u_\alpha(a^{-1})s_\alpha B(K).$$

Therefore by 4.2, $ws_\alpha u_{-\alpha}(a) \in u_{-w(\alpha)}(\epsilon a^{-1})wB(K)$ with $\epsilon = \pm 1$. Hence in the two remaining cases, we have:
If $w(\alpha) > 0$ and $a \notin A$ then $a^{-1} \in \mathcal{M}_A$ and thus $ws_\alpha u_{-\alpha}(a) \in U^-(\mathcal{M}_A)wB(K) \subset IwB(K)$.
If $w(\alpha) < 0$ and $a \notin \mathcal{M}_A$ then $a^{-1} \in A$ and thus $ws_\alpha u_{-\alpha}(a) \in \mathcal{U}_B(A)wB(K) \subset IwB(K)$. ∎

In order to state the next lemma, we recall some definitions about the Bruhat's order. For $w \in W$, we can write $w = s_{\alpha_1} \ldots s_{\alpha_n}$ with the $\alpha_i \in \Phi$ and $n$ minimal; n is then called the length of $w$. For $w' \in W$, we write $w' \prec w$ if $w' = s_{\alpha_{i_1}} \ldots s_{\alpha_{i_k}}$ with $1 \leq i_1 < \ldots < i_k \leq n$. Of course, if $w' \prec w$ we have $length(w') \leq length(w)$.

**Lemma 4.2** *Let $w \in W$, then we have the following inclusion:*

$$\mathcal{U}_B(K)w \subset \coprod_{w' \prec w} Iw'B(K)$$

*Moreover if $u \in \mathcal{U}_B(K) - \mathcal{U}_B(A)$ and $w^{-1}uw \notin \mathcal{U}_B(K)$, then $uw \in Iw'B(K)$ with $length(w') < length(w)$.*

**Proof** We proceed by induction on the length of $w$. Let us write $w = w_1 s_\alpha$ with $\alpha \in \Phi$ and $length(w_1) = length(w) - 1$ and fix $u \in \mathcal{U}_B(K)$. By induction, $uw_1 \in Iw_1'B(K)$ with $w_1' \prec w_1$. Thus there exist $g_0 \in I$, $t \in T(K)$ and $u' \in \mathcal{U}_B(K)$ such that $uw = g_0 w_1' u' s_\alpha t$. Now we can write $u's_\alpha = s_\alpha u_{-\alpha}(a)u''$ with $a \in K$ and $u'' \in \mathcal{U}_B(K)$. Thus by the previous lemma, $w_1' u' s_\alpha \in Iw_1 B(K) \cup Iw_1' s_\alpha B(K)$. This proves the lemma because $w_1' \prec w$ and $w_1' s_\alpha \prec w$. ∎

## 4.2 Ramification theory of the cusps

We denote by $H_p^1 = H^1 \otimes_{\mathbf{Z}} \mathbf{Z}_p = \prod_{v|p} Sp_{4/\mathbf{r}_v}$ and $W_p$ a subgroup of $H_p^1$ isomorphic to the Weyl group $(N_{H_p^1}(T)/T)_{/\mathbf{Z}_p}$ of $H_p^1$. We assume that $W_p = \prod_{v|p} W_v$ where each $W_v \subset G_{/\mathbf{r}_v}$ maps bijectively to the Weyl group of $Sp_4$ isomorphic to the dihedral group $D_8$. For any element $x$ in one of the previous groups, we will denote by $x(v)$ its projection on the $v$-component.

If $\mathbf{P}$ is a parabolic subgroup of $H_{\mathbf{Z}}^1$, we say that it is of Klingen type (resp. Siegel, Borel) if it is conjugate to the standard Klingen parabolic $P^*$ (resp. the standard Siegel parabolic $P$, or the standard Borel subgroup $B$). An arbitrary type is denoted by the letter $\Sigma$; we write $P_\Sigma$ for the corresponding standard parabolic subgroup and we put $\mathcal{P}_\Sigma = P_\Sigma \otimes \mathbf{Z}_p \subset H_p^1$.

Consider now the so-called Bruhat decomposition:

$$H^1(\mathbb{F}_p) = \coprod_{w \in W_p} \mathcal{U}_B(\mathbb{F}_p).w.B(\mathbb{F}_p)$$

The group $W_\Sigma = W_p \cap \mathcal{P}_\Sigma$ is isomorphic to the Weyl group of the Levi of $\mathcal{P}_\Sigma$ and we have:

$$\mathcal{P}_\Sigma(\mathbb{F}_p) = \coprod_{w \in W_\Sigma} B(\mathbb{F}_p).w.B(\mathbb{F}_p)$$

We deduce the Bruhat decomposition for $\mathcal{P}_\Sigma$:

$$H^1(\mathbb{F}_p) = \coprod_{w \in W/W_\Sigma} \mathcal{U}_B(\mathbb{F}_p).w.\mathcal{P}_\Sigma(\mathbb{F}_p) \qquad (4.4)$$

In this section, for any subgroup $\Gamma'$ with $\Gamma(p^r) \subset \Gamma' \subset \Gamma$ we consider the set $Cusp_\Sigma(\Gamma')$ of $\Gamma'$-conjugacy classes of type $\Sigma$ parabolic subgroups. For any parabolic subgroup $\mathbf{P}$ of type $\Sigma$, we note $s_\mathbf{P}$ the $\Gamma'$-conjugacy class of $\mathbf{P}$. Since $\Gamma$ is an arithmetic subgroup of level prime to $p$, we can write $Cusp_\Sigma(\Gamma) = \{s_{\mathbf{P}_1}, \ldots, s_{\mathbf{P}_t}\}$ where the $\mathbf{P}_i$'s are parabolic subgroups of $H^1$ such that $\mathbf{P}_i \otimes \mathbf{Q}_p = g_i(\mathcal{P}_\Sigma \otimes \mathbf{Q}_p) g_i^{-1}$ with $g_i \in H^1(\mathbf{Z}_p)$ and $g_i \equiv 1$ mod. $p^r$. Since $\mathcal{P}_\Sigma$ is a parabolic subgroup over $\mathbf{Z}_p$, we have by Iwasawa decomposition:

$$H^1(\mathbf{Q}_p)/\mathcal{P}_\Sigma(\mathbf{Q}_p) = H^1(\mathbf{Z}_p)/\mathcal{P}_\Sigma(\mathbf{Z}_p).$$

Then we can write:

$$Cusp_\Sigma(\Gamma') = \coprod_{i=1}^t \Gamma'\backslash\Gamma/\Gamma \cap \mathbf{P}_i = \coprod_{i=1}^t \bar{\Gamma}'\backslash H^1(\mathbf{Z}/p^r\mathbf{Z})/\mathcal{P}_\Sigma^1(\mathbf{Z}/p^r\mathbf{Z})\bar{E}_{\mathbf{P}_i}(\Gamma) \quad (4.5)$$

where the overline means the reduction modulo $p^r$ and $E_{\mathbf{P}_i}(\Gamma) = \Gamma \cap \mathbf{P}_i/\Gamma \cap \mathbf{P}_i^1$ and noting $\mathbf{P}^1 = M_\mathbf{P}^1 \mathcal{U}_\mathbf{P}$ for a Levi decomposition $\mathbf{P} = M_\mathbf{P}\mathcal{U}_\mathbf{P}$. Let $\mathcal{C}_\Sigma(\Gamma') = \bar{\Gamma}'\backslash H^1(\mathbf{Z}/p^r\mathbf{Z})/\mathcal{P}_\Sigma(\mathbf{Z}/p^r\mathbf{Z})$, then there is a canonical map $\pi$:

$$Cusp_\Sigma(\Gamma') \to \Gamma'\backslash H^1(\mathbf{Z}_p)/\mathcal{P}_\Sigma(\mathbf{Z}_p) \to \mathcal{C}_\Sigma(\Gamma') \quad (4.6)$$

¿From the decomposition 4.5, we can see that for $i \in \{1, \ldots, t\}$ $\pi$ restricted to $\Gamma'\backslash\Gamma/\Gamma \cap \mathbf{P}_i$ is one-to-one.

Let $g \in H^1(\mathbf{Z}_p)$, then by (4.4), we can write

$$g = g_0 w_g p_\Sigma \quad (4.7)$$

with $p_\Sigma \in \mathcal{P}_\Sigma(\mathbf{Z}_p)$, $w_g \in W$ and $g_0 \in I$. Let $r_g(v) = Sup \{i; g_0(v) \in I_{v,i}\}$; here the $I_{v,i}$'s are those defines in section 2.1. For $g \in H^1(\mathbf{Q}_p)$, we denote by $[g]$ its class in $\mathcal{C}_\Sigma(\Gamma')$. Let $s \in Cusp_\Sigma(\Gamma')$ and $g \in H^1(\mathbf{Z}_p)$ such that $\pi(s) = [g]$. Let us recall that we deal with two parabolic subgroups $\mathcal{P}_\Sigma$ and $\mathcal{Q}$. The first one defines the type of a stratum of the Borel-Serre compactification and the second one defines the ordinary idempotent $e_\mathcal{Q}$. Then if we let $W_Q$ the Weyl group of $Q$ ($W_Q = W \cap Q$), the class of $w_g$ in $\underline{W}_{Q\Sigma} = W_Q\backslash W_p/W_\Sigma$ is independent of the representative $g$ of $s$ and is denoted by $w_s$.

**Definition 4.1** *The class $w_s$ in $\underline{W}_{Q\Sigma} = W_Q\backslash W_p/W_\Sigma$ is called the Weyl type of $s$. We set $r_s(v) = Sup\{r_g(v)$ where $[g] = \pi(s)\}$; we call $r_s(v)$ the $v$-depth of $s$.*

**Remark** If $r_s(v) \geq r$ then $r_s(v) = \infty$ because by successive modifications by elements of $\Gamma' \supset \Gamma(p^r)$ one can let $r_g(v)$ grow to infinity. Therefore $r_s(v) \in \{1, 2, \ldots, r-1, \infty\}$.

For any $h, g$ elements in a group, we write $g^h = hgh^{-1}$. We denote by $R_Q^-$ the roots of $\mathcal{U}_Q^-$, the unipotent subgroup opposite to $\mathcal{U}_Q$.

**Lemma 4.3** *Let $v|p$ and $\Gamma'$ as above. Let $\beta \in H^1(\mathbf{Q})$ such that $\beta \equiv tu \mod. p^r$ with $u \in \mathcal{U}_B(\mathbf{Z}_p)$, $t \in T(\mathbf{Q}_p)$. We suppose that $v(\alpha(t)) > 0$ for all $\alpha \in R_Q^-$. Let $\mathbf{P}$ be a parabolic subgroup of $H^1$ and let us set $s = s_{\mathbf{P}} \in Cusp_\Sigma(\Gamma')$ and $s' = s_{\beta \mathbf{P} \beta^{-1}}$. Then we have two possibilities:*
a) $w_{s'}(v) = w_s(v)$ and $r_{s'}(v) > r_s(v)$.
b) *or else* $length(w_{s'}(v)) \leq length(w_s(v)) - 1$.

**Proof** In what follows, if $r_s = \infty$, $I_{v,r_s}$ means $I_{v,r}$ and $\mathcal{U}_Q^-(p^{r_s(v)}\mathbf{r}_v)$ means $\mathcal{U}_Q^-(p^r \mathbf{r}_v)$. Let us write $\beta = \delta tu$ with $\delta \in \Gamma(p^r)$ and choose by 4.7 $g_0 \in I_{v,r_s}$ such that $\pi(s) = [g_0 w_s]$. Then $\pi(s') = [\beta g_0 w_s] = [\delta tu g_0 w_s]$. By the Iwahori decomposition, we can write

$$ug_0 = u^- u^+ t_0 \text{ with } u^- \in \mathcal{U}_Q^-(p^{r_s(v)}\mathbf{r}_v),\ u^+ \in Q(\mathbf{Q}_p) \text{ and } t_0 \in T(\mathbf{Z}_p)$$

Thus $\pi(s') = [\delta(u^-)^t (u^+)^t w_s]$. We deal with two cases:
a) $(u^+)^t(v) \in Q(\mathbf{r}_v)$. Since $(u^-)^t(v) \in \mathcal{U}_Q^-(p^{r_s(v)+1})$ by the hypothesis on $t$, we have $(\delta(u^-)^t(u^+)^t)(v) \in I_{v,r_s(v)+1}$. Thus $w_{s'}(v) = w_s(v)$ and $r_{s'}(v) > r_s(v)$.
b) $(u^+)^t(v) \notin Q(\mathbf{r}_v)$. Then we can write $(u^+)^t = u_1^+ u_1^-$ with $u_1^+(v) \in \mathcal{U}_B(F_v) \setminus \mathcal{U}_B(\mathbf{r}_v)$ and $u_1^- \in \mathcal{U}_B^-(\varpi_v \mathbf{r}_v)$. Then by the lemma 4.2, we have $\pi(s') = [g_0 w']$ with $g_0 \in I_{v,1}$, $length(w'(v)) \leq length(w_s(v)) - 1$. ∎

In the sequel, we denote $r_s = inf\{r_s(v); v|p\}$ and

$$l(w) = \sum_{v|p} [F_v : \mathbf{Q}_p] \times length(w(v)).$$

We note $Cusp_{\Sigma,w}(\Gamma')$ the subset of $Cusp_\Sigma(\Gamma')$ consisting in the cusps of Weyl type $w$ and depth $\infty$.

We recall that $C_Q = M_Q/M_Q^1 = T/T \cap M_Q^1$ is the cocenter of $M_Q$ and we note $i_Q$ the canonical surjective arrow from $T$ to $C_Q$.

For any $s_{\mathbf{P}_i} \in Cusp_\Sigma(\Gamma)$, we set $Cusp_{\Sigma,w}(s_{\mathbf{P}_i}, \Gamma')$ to be the set of cusps $s \in Cusp_{\Sigma,w}(\Gamma')$ lying over $s_{\mathbf{P}_i}$. Let us observe that as a consequence of the Bruhat decomposition relative to the parabolic subgroups $Q$ and $\mathcal{P}_\Sigma$, the set $Cusp_{\Sigma,w}(s_{\mathbf{P}_i}, \Gamma_0(p^r))$ consists in only one element. We choose a representative $\mathbf{P}_{i,w}$ of this class; the class itself is therefore denoted by $s_{\mathbf{P}_{i,w}}$.

**Definition 4.2** *We put:* $M_\Sigma^w = wM_\Sigma w^{-1}$, $T_{\Sigma,w}^1 = T \cap M_{\Sigma,w}^1$ *while $E_{i,w}(\Gamma)$ is defined as the image of $\Gamma \cap T(\mathbf{Z}_p)/\Gamma \cap \mathbf{P}_{i,w}^1 \cap T(\mathbf{Z}_p)$ in $C_Q(\mathbf{Z}/p^r\mathbf{Z})$*

**Lemma 4.4** *Let $s_{\mathbf{P}_i} \in Cusp_\Sigma(\Gamma)$. Then*

- $Cusp_{\Sigma,w}(s_{\mathbf{P}_i}, \Gamma_0(p^r)) = \{s_{\mathbf{P}_{i,w}}\}$
- $Cusp_{\Sigma,w}(s_{\mathbf{P}_i}, \Gamma_1(p^r)) \cong C_Q(\mathbf{Z}/p^r\mathbf{Z})/i_Q(T_{\Sigma,w}^1(\mathbf{Z}/p^r\mathbf{Z})\overline{E}_{i,w}(\Gamma))$

**Proof:** The first point has already been noticed. The second follows from the calculation of the stabilizer of $s_{\mathbf{P}_{i,w}}$ in $C_Q(\mathbf{Z}/p^r\mathbf{Z}) = \Gamma_0(p^r)/\Gamma_1(p^r)$. ∎

## 4.3 On the p-ordinary cohomology associated to a type $\Sigma$

For any $s = s_\mathbf{P} \in Cusp_\Sigma(\Gamma')$, we set $\Gamma'_s = \mathbf{P} \cap \Gamma'$. We begin now the study of the ordinary part of the following cohomology groups:

$$G^q_\Sigma(\Gamma'; M) = \bigoplus_{s_\mathbf{P} \in Cusp_\Sigma(\Gamma')} H^q(\Gamma'_{s_\mathbf{P}}; M)$$

If $\xi \in H^1(\mathbf{Q})$ is such that $\xi\Gamma'\xi^{-1} \cap \Gamma'$ is of finite index in $\Gamma$ then the double class $[\Gamma'\xi\Gamma']$ acts on this cohomology group as follows. Let us first decompose $\Gamma'\xi\Gamma'_s = \coprod_j \Gamma'\xi(j)$. Then if $c = \oplus c_{s_\mathbf{P}} \in G^q_\Sigma(\Gamma'; M)$ then

$$(c|[\Gamma'\xi\Gamma'])_{s_\mathbf{P}} = \oplus_j c_{s_{\xi(j)\mathbf{P}\xi(j)^{-1}}} |[\Gamma'_{s_{\xi(j)P\xi(j)^{-1}}}\xi(j)\Gamma'_{s_\mathbf{P}}]$$

Using definition 4.1, we define the following cohomology subgroups:

$$G^q_{\Sigma,l,t}(\Gamma'; M) = \bigoplus_{\substack{s_\mathbf{P} \in Cusp_\Sigma(\Gamma') \\ r_{s_\mathbf{P}} \leq t,\ length(w_{s_\mathbf{P}}) = l \\ or\ length(w_{s_\mathbf{P}}) > l}} H^q(\Gamma'_{s_\mathbf{P}}; M)$$

$$G^q_{\Sigma,w}(\Gamma'; M) = \bigoplus_{Cusp_{\Sigma,w}(\Gamma')} H^q(\Gamma'_{s_\mathbf{P}}; M)$$

By lemma 4.3, the Hecke operators $[\Gamma'\xi\Gamma']$ act on these cohomology subgroups and we have:

**Corollary 4.1** *The Hida idempotent $e_Q$ annihilates the quotients:*

$$G^q_{\Sigma,l-1,r-1}(\Gamma'; M)/G^q_{\Sigma,l,\infty}(\Gamma'; M)$$

*and we have*

$$e_Q G^q_\Sigma(\Gamma'; M) = \bigoplus_{w \in \underline{W}_{Q\Sigma}} e_Q G^q_{\Sigma,w}(\Gamma', M)$$

**Proof** Let $n$ be an integer; for the $\xi_Q$ defined at the end of section 2.4, consider the right coset decomposition

$$\Gamma'\xi_Q^{n!}\Gamma'_{s_\mathbf{P}} = \coprod_i \Gamma'\xi(j)$$

where $s_\mathbf{P} \in Cusp_\Sigma(\Gamma')$ has Weyl type of length $l - 1$. Then by lemma 4.3, $s_{\xi(j)\mathbf{P}\xi(j)^{-1}}$ has Weyl type of length less than $l - 2$ or has a depth $\infty$ if $n$ is chosen sufficiently large. Thus, for $c \in G^q_{\Sigma,l-1,r-1}(\Gamma'; M)$, $(c|[\Gamma'\xi(j)\Gamma'])_{s_\mathbf{P}} = 0$ thus $c|[\Gamma'\xi(j)\Gamma'] \in G^q_{\Sigma,l,\infty}(\Gamma'; M)$ as desired. For the second point, let us note that we have a filtration:

$$\ldots \subset G^q_{\Sigma,l,\infty}(\Gamma'; M) \subset G^q_{\Sigma,l-1,r-1}(\Gamma'; M) \subset$$

$$\subset G^q_{\Sigma,l-1,\infty}(\Gamma';M) \subset G^q_{\Sigma,l-2,r-1}(\Gamma';M) \subset \ldots$$

the ordinary part of which

$$\ldots \subset e_Q G^q_{\Sigma,l,\infty}(\Gamma';M) = e_Q G^q_{\Sigma,l-1,r-1}(\Gamma';M) \subset$$

$$e_Q G^q_{\Sigma,l-1,\infty}(\Gamma';M) = e_Q G^q_{\Sigma,l-2,r-1}(\Gamma';M) \subset \ldots$$

is still a filtration of $e_Q G^q_\Sigma(\Gamma';M)$ with quotient isomorphic to

$$e_Q(G^q_{\Sigma,l-1,\infty}(\Gamma';M)/G^q_{\Sigma,l-1,r-1}(\Gamma';M))$$

since $e_Q$ is an idempotent. The result follows from this observation. ∎

On the other hand, by the last part of lemma 4.3, we can decompose the action of $e_Q$ on $G^q_{\Sigma,w}(\Gamma')$ into an action on each $H^q(\Gamma'_s;M)$ with $w_s = w$ and $r_s = \infty$ and thus:

$$e_Q G^q_{\Sigma,w}(\Gamma';M) = \bigoplus_{\substack{s \in Cusp_\Sigma(\Gamma) \\ r_s = \infty,\ w_s = w}} e_Q(\Gamma'_s) H^q(\Gamma'_s;M)$$

where $e_Q(\Gamma'_s) = \lim_{n \to \infty} [\Gamma'_s \xi_s \Gamma'_s]^{n!}$ with $\xi_s$ normalizing $\mathbf{P}_s$.

What we are aiming at is the study of the action of $Z_Q(\mathbf{Z}/p^r\mathbf{Z})$ on $e_Q G^p_\Sigma(\Gamma_1(p^r);M)$ via the canonical isogeny :

$$Z_Q(\mathbf{Z}/p^r\mathbf{Z}) \to C_Q(\mathbf{Z}/p^r\mathbf{Z}) \cong \Gamma_0(p^r)/\Gamma_1(p^r)$$

$$z \longmapsto \bar{\gamma}_z$$

given by the natural action of $\Gamma_0(p^r)/\Gamma_1(p^r)$ on our cohomology groups. If $\gamma \in \Gamma_0(p^r)$ normalizes $\mathbf{P}_{i,w}$ (i.e $\bar{\gamma} \in i_Q(T^1_{\Sigma,w}(\mathbf{Z}/p^r\mathbf{Z})\overline{E}_{i,w}(\Gamma)))$, it operates on the cohomology $H^p(\Gamma_1(p^r) \cap \mathbf{P}_{i,w};M))$.

**Proposition 4.1** *With the previous notations, we have:*

$$e_Q G^\bullet_{\Sigma,w}(\Gamma_1(p^r);M) = \bigoplus_{s_{\mathbf{P}_i} \in Cusp_\Sigma(\Gamma)} Ind^{C_Q(\mathbf{Z}/p^r\mathbf{Z})}_{i_Q(T^1_{\Sigma,w}(\mathbf{Z}/p^r\mathbf{Z})\overline{E}_{i,w}(\Gamma))} e_Q H^\bullet(\Gamma_1(p^r) \cap \mathbf{P}_{i,w};M))$$

**Proof** For any $s_{\mathbf{P}} \in Cusp_{\Sigma,w}(s_{\mathbf{P}_i}, \Gamma_1(p^r))$, there exists $\gamma \in \Gamma_0(p^r)$ such that $\mathbf{P} = \mathbf{P}^\gamma_{i,w}$. Since $\Gamma_0(p^r)$ normalizes $\Gamma_1(p^r)$, we have $\Gamma_1(p^r) \cap \mathbf{P} = (\Gamma_1(p^r) \cap \mathbf{P}_{i,w})^\gamma$ and thus $H^q(\Gamma_1(p^r) \cap \mathbf{P};M) \cong H^q(\Gamma_1(p^r) \cap \mathbf{P}_{i,w};M)$. Then the proposition follows easily from this observation and the lemma 4.4. ∎

In the sequel we fix an arithmetic character $\chi = \varepsilon\psi$ of $C_Q$ where $\psi$ is algebraic and $\varepsilon$ factors through $C_Q(\mathbf{Z}/p^r\mathbf{Z})$ ($p^r$ is a level of $\chi$ cf. Definition 3.2). We can now prove the independence of the weight for the $\Sigma$-stratum ordinary cohomology:

**Corollary 4.2** *With the previous notations, we have a canonical isomorphism:*
$$e_Q G^\bullet_{\Sigma,w}(\Gamma_1(p^\infty), L^a(\rho \otimes \chi; A)) = \varinjlim_r e_Q G^\bullet_{\Sigma,w}(\Gamma_1(p^r); L^a(\rho \otimes \chi; A))$$
$$\cong e_Q G^\bullet_{\Sigma,w}(\Gamma_0(p), \mathcal{C}(\rho^1; A))$$

**Proof:** By the same proof of corollary 3.2, we can prove that
$$e_Q H^\bullet(\mathbf{P}_{i,w} \cap \Gamma_1(p^s), L^a(\rho \otimes \chi); A_r) \cong e_Q H^\bullet(\mathbf{P}_{i,w} \cap \Gamma_1(p^s), V(\rho) \otimes A_r)$$
$$e_Q H^\bullet(\mathbf{P}_{i,w} \cap \Gamma_0(p^r), \mathcal{C}(\rho^1; A_r)) \cong e_Q H^\bullet(\mathbf{P}_{i,w} \cap \Gamma_0(p^r), \mathcal{C}(M'(\mathbf{Z}_p/p^r\mathbf{Z}_p), \rho^1; A_r))$$

On the other hand, by Shapiro's lemma, we have still
$$e_Q H^\bullet(\mathbf{P}_{i,w} \cap \Gamma_1(p^r), V(\rho) \otimes A_r) = e_Q H^\bullet(\mathbf{P}_{i,w} \cap \Gamma_0(p^r); Ind^{\mathbf{P}_{i,w} \cap \Gamma_0(p^r)}_{\mathbf{P}_{i,w} \cap \Gamma_1(p^r)} V(\rho) \otimes A_r)$$

Therefore by the first isomorphism, we have
$$e_Q G^\bullet_{\Sigma,w}(\Gamma_1(p^r), L^a(\rho \otimes \chi; A_r)) =$$
$$Ind^{C_Q(\mathbf{Z}/p^r\mathbf{Z})}_{i_Q(T^1_{\Sigma,w}(\mathbf{Z}/p^r\mathbf{Z})\overline{E}_{i,w}(\Gamma))} e_Q H^\bullet(\mathbf{P}_{i,w} \cap \Gamma_0(p^r); Ind^{\mathbf{P}_{i,w} \cap \Gamma_0(p^r)}_{\mathbf{P}_{i,w} \cap \Gamma_1(p^r)} V(\rho) \otimes A_r))$$
$$= H^\bullet(\mathbf{P}_{i,w} \cap \Gamma_0(p^r), Ind^{C_Q(\mathbf{Z}/p^r\mathbf{Z})}_{i_Q(T^1_{\Sigma,w}(\mathbf{Z}/p^r\mathbf{Z})\overline{E}_{i,w}(\Gamma))} Ind^{\mathbf{P}_{i,w} \cap \Gamma_0(p^r)}_{\mathbf{P}_{i,w} \cap \Gamma_1(p^r)} V(\rho) \otimes A_r))$$

By transitivity of the induction process, we have:
$$\mathcal{C}(M'(\mathbf{Z}_p/p^r\mathbf{Z}_p), \rho^1; A_r) = Ind^{C_Q(\mathbf{Z}/p^r\mathbf{Z})}_{i_Q(T^1_{\Sigma,w}(\mathbf{Z}/p^r\mathbf{Z})\overline{E}_{i,w}(\Gamma))} Ind^{\mathbf{P}_{i,w} \cap \Gamma_0(p^r)}_{\mathbf{P}_{i,w} \cap \Gamma_1(p^r)} V(\rho) \otimes A_r$$

We deduce now our result from these previous isomorphisms and taking the inductive limit. ■

By the same arguments of the proof of theorem 3.2, one can prove the following lemma. We leave it to the reader:

**Lemma 4.5** *For all $\chi$ of level $p^r$, we have a canonical isomorphism*
$$e_Q G^\bullet_{\Sigma,w}(\Gamma_0(p^r), L^a(\rho \otimes \chi; A)) \cong e_Q G^\bullet_{\Sigma,w}(\Gamma_0(p^r); \mathcal{C}(\rho; A)[\omega_\chi]) \tag{4.8}$$

In section 3, we proved a control theorem for the total cohomology in degree $3d$ using an abstract control criterion due to Hida. To apply this criterion, we started from the identification 3.9 of which 4.8 is an analogue. Then, we used vanishing results of the cohomology groups in degree less than $3d$. Unfortunately, this last point is not satisfied for the $\Sigma$-stratum cohomology or for the boundary cohomology. We thus need to go further in the computation of this cohomology. In particular, we will see in Section 5 how to compute them explicitly for $F = \mathbf{Q}$ and to get the control for the boundary cohomology. The general case is treated only modulo torsion and we will not obtain a control theorem for the boundary cohomology but only a weaker result, which nevertheless proves sufficient to imply a control theorem for the interior cohomology.

## 4.4 Ordinary cohomology of parabolic subgroups

### 4.4.1 An abstract lemma

We begin with some abstract considerations. Let $\mathfrak{g}$ be a group and $\mathfrak{Z}$ be a normal subgroup of $\mathfrak{g}$; we set $\bar{\mathfrak{g}} = \mathfrak{g}/\mathfrak{Z}$ and $g \mapsto \bar{g}$ the reduction map. Let $\mathfrak{N}$ be a subgroup of $\mathfrak{g}$ such that $\mathfrak{N}_1 = \mathfrak{N} \cap \mathfrak{Z}$ is a non trivial subgroup of $\mathfrak{N}$. Let $\eta$ be a element of $\mathfrak{g}$ such that $\mathfrak{N} \cap \eta\mathfrak{N}\eta^{-1}$ is of finite index in $\mathfrak{N}$ and in $\eta\mathfrak{N}\eta^{-1}$. We suppose that there exist elements $\eta_1$ and $\eta'$ in $\mathfrak{g}$ such that:

(1) $\bar{\eta}_1$ normalizes $\mathfrak{N}_0 = \mathfrak{N}/\mathfrak{N}_1$

(2) $\eta'$ normalizes $\mathfrak{N}_1$

(3) $\eta = \eta'\eta_1 = \eta_1\eta'$

(4) $\eta\mathfrak{Z} \cap \mathfrak{N}\eta\mathfrak{N} = \eta'\mathfrak{N}_1\eta_1\mathfrak{N}_1$

We say that the quadruple $(\mathfrak{N}, \mathfrak{N}_1, \eta, \eta_1)$ is admissible. We can now consider the following double classes:

$$[\mathfrak{N}\eta\mathfrak{N}], \ [\mathfrak{N}_1\eta_1\mathfrak{N}_1] \ and \ [\mathfrak{N}_0\eta_0\mathfrak{N}_0]$$

where $\eta_0 = \bar{\eta}'$. Let $M$ be a $\mathfrak{N}$-module (resp. $N$ a $\mathfrak{N}_0$-module), with an action by $\eta^{-1}$ and $\eta_1^{-1}$ (resp. $\eta_0^{-1}$). Then these double classes act respectively on the cohomology groups

$$H^*(\mathfrak{N}; M), \ H^*(\mathfrak{N}_1; M) \ (\text{resp. } H^*(\mathfrak{N}_0; N)).$$

and if these modules are finite or cofinite over $\mathbf{Z}_p$, we can consider the respective idempotents $e$, $e_1$ and $e_0$; this assumption will be understood in the following lemma. For a group $G$ and a $G$-module $A$, let $C^\bullet(G, A) = \mathrm{Hom}_{\mathbf{Z}}(\mathbf{Z}G^{\bullet+1}, A)$ (the standard homogeneous complex giving rise to $H^\bullet(G, A)$ is the fixed part $C^\bullet(G, A)^G$). Let us consider the double complex which gives the Hochschild-Serre spectral sequence:

$$L^{p,q} = C^p(\mathfrak{N}/\mathfrak{N}_1; C^q(\mathfrak{N}; M)^{\mathfrak{N}_1})^{\mathfrak{N}/\mathfrak{N}_1}.$$

**Lemma 4.6** *Let $(\mathfrak{N}, \mathfrak{N}_1, \eta, \eta_1)$ be an admissible quadruple.*

(i) *If $\mathfrak{N}_1\eta_1\mathfrak{N}_1 = \coprod_i \mathfrak{N}_1\eta_1 z_i$ and $\mathfrak{N}_0\eta_0\mathfrak{N}_0 = \coprod_j \mathfrak{N}_0\eta_0\bar{n}_j$ with $n_j \in \mathfrak{N}$ then*

$$\mathfrak{N}\eta\mathfrak{N} = \coprod_{i,j} \mathfrak{N}\eta_1 z_i \eta' n_j.$$

(ii) *Let $c \in L^{p,q}$. Then we define an action on the double complex $L^{p,q}$ by*

$$c|[\mathfrak{N}\theta\mathfrak{N}](\bar{x}_1, \ldots, \bar{x}_p)(y_1, \ldots, y_q) =$$
$$\sum_{i,j} (\eta_1 z_i \eta' n_j)^{-1} . c(\bar{x}_1^{(j)}, \ldots, \bar{x}_p^{(j)})(y_1^{(i,j)}, \ldots, y_q^{(i,j)})$$

*with $\eta_0\bar{n}_j\bar{x} = \bar{x}^{(j)}\eta_0\bar{n}_{j'}$ and $\eta_1 z_i \eta' n_j y = y^{(i,j)}\eta_1 z_{i'}\eta' n_{j'}$.*

*(iii) The Hochschild-Serre spectral sequence induces a converging spectral sequence:*

$$e_0 H^p(\mathfrak{N}_0; e_1 H^q(\mathfrak{N}_1; M)) \Rightarrow eH^{p+q}(\mathfrak{N}; M).$$

**Proof** (i) Let $x \in \mathfrak{N}\eta\mathfrak{N}$. Since $\bar{\eta}_1$ normalizes $\mathfrak{N}_0 = \mathfrak{N}/\mathfrak{N}_1$, we have

$$\begin{aligned}\mathfrak{N}_0 \bar{\eta} \mathfrak{N}_0 &= \mathfrak{N}_0 \bar{\eta}_1 \eta_0 \mathfrak{N}_0 \\ &= \mathfrak{N}_0 \eta_0 \mathfrak{N}_0 \bar{\eta}_1 \\ &= \coprod_j \mathfrak{N}_0 \eta_0 \bar{n}_j \bar{\eta}_1.\end{aligned}$$

Thus there exists $j$ such that $\bar{x} = \bar{n}\eta_0 \bar{n}_j \bar{\eta}_1$. Then $y = (n\eta')^{-1} x n_j^{-1} \in \eta_1 \mathfrak{Z} \cap \eta'^{-1} \mathfrak{N}\eta\mathfrak{N} = \mathfrak{N}_1 \eta_1 \mathfrak{N}_1$. Thus there exists $z, z' \in \mathfrak{N}_1$ such that $y = z\eta_1 z'$ and

$$x = n\eta' z \eta_1 z' n_j = n(\eta' z \eta'^{-1})\eta_1(\eta' z' \eta'^{-1})\eta' n_j$$

By (2) $(\eta' z \eta'^{-1})$ and $(\eta' z' \eta'^{-1})$ belong to $\mathfrak{N}_1$. Thus $x \in n\mathfrak{N}_1 \eta_1 \mathfrak{N}_1 \eta' n_j$ and there exists $i$ such that $x \in \mathfrak{N}\eta_1 z_i \eta' n_j$. We have proved that $\mathfrak{N}\eta\mathfrak{N} = \bigcup_{i,j} \mathfrak{N}\eta_1 z_i \eta' n_j$; the fact that the union is disjoint is obvious.

The second point consists in verifying that the action preserves $L^{p,q}$ and that it commutes with the differential maps; the reader can consult for this [16] (cf. p. 299). For the third point, it is enough to check that

$$(c|[\mathfrak{N}_1 \eta_1 \mathfrak{N}_1])|[\mathfrak{N}_0 \eta_0 \mathfrak{N}_0] = c|[\mathfrak{N}\eta\mathfrak{N}]$$

and remark that the action defined in (iii) on $L^{0,q}$ is the usual action of $[\mathfrak{N}\eta\mathfrak{N}]$. ∎

### 4.4.2 Application to the ordinary cohomology of unipotent subgroups

For all $v|p$ let $R_v = R((Sp_4) \otimes \mathbf{r}_v, T \otimes \mathbf{r}_v)$ be the set of roots of $(Sp_4) \otimes \mathbf{r}_v$ with respect to $T \otimes \mathbf{r}_v$ and let $R_p = \coprod_{v|p} R_v$. For any subset $R \subset R_p$ such that $(R + R) \cap R_p \subset R$, we denote by $\mathcal{U}_R$ the unipotent subgroup of $H(\mathbf{Z}_p)$ associated with $R$.

For any subset $R_0 \subset R_\mathbf{Z} = R(H_{/\mathbf{Z}}, T_{H/\mathbf{Z}})$ of the set of rational roots, for each $v|p$, let $R_0(v)$ be the image of $R_0$ in $R_v$ by restriction to $T \otimes \mathbf{r}_v \subset T_H \otimes \mathbf{Z}_p$; let $R_0(p) = \coprod_{v|p} R_0(v)$. We denote by $U_{R_0}(p^r) = \Gamma_1(p^r) \cap \mathcal{U}_{R_0(p)}$.

We denote $R^+$ ( resp. $R^-$) the subset of $R$ of positive roots (resp. negative roots) (according to the choice of the standard Borel subgroup). And we set

$$D_{Q,R} = \{t \in T_H(\mathbf{Q}_p) \cap \prod_{v|p} M_4(\mathbf{r}_v) \,;\, v(\alpha(t)) > 0 \ \forall v|p, \forall \alpha \in R \cap R_Q^-\}$$

For any $t \in D_{Q,R}$, one can choose $\xi_t \in H(\mathbf{Q})$ such that $t^{-1}\xi_t \equiv 1 \bmod p^r$ (modulo the center $Z_H$). The construction of such an element is similar to that of Section 2.4. We consider the double class $[U_R(p^r)\xi_t U_R(p^r)]$ and $e_{Q,R}$ the associated idempotent (it is easy to prove that idempotent is independent of $t$). Sometimes we note it $e_Q$ if there is no possible confusion.

Let $R_1 \subset R \subset R_{\mathbf{Q}}$ verifying $(R+R_1) \cap R = \emptyset$. Then $\mathcal{U}_{R_1}$ is central in $\mathcal{U}_R$. If $t \in D_{Q,R}$, one can find $t_1 \in \Delta_{Q,R_1}$ such that

$$v(\alpha(t_1)) = 0 \text{ for all } v|p \text{ and } \alpha \in R - R_1$$

The quadruple $(U_R(p^r), U_{R_1}(p^r), \xi_t, \xi_{t_1})$ is admissible (i.e. it verifies the four conditions stated before the Lemma 4.6). Therefore, if $M$ is a $(D_{Q,R}^{-1}, U_R(p^r))$-module, by lemma 4.6, we have the ordinary part of the Hochschild-Serre spectral sequence:

$$\bar{e}_{Q,R}H^p(\overline{U_R(p^r)}, e_{Q,R_1}H^q(U_{R_1}(p^r), M)) \Rightarrow e_{Q,R}H^{p+q}(U_R(p^r), M)$$

with $\overline{U_R(p^r)} = U_R(p^r)/U_{R_1}(p^r)$.

**Definition 4.3** *For all $R \subset R_p$, we set*

$$|R| := \sum_{v|p}[F_v : \mathbf{Q}_p] \times Card(R \cap R_v)$$

Let us denote $R_\Sigma$ (resp. $R_Q$) the subset of positive roots corresponding to the unipotent radical of $\mathcal{P}_\Sigma$ (resp. $Q$). We recall also that $\Delta_Q$ denotes the set of roots of $M_Q$ with respect to $T$.

**Lemma 4.7** *Let $M$ be any $\mathbf{Z}_p[U_{R_0}]$-module with a trivial action of $U_{R_0} \cap \mathcal{U}_Q$ and $U_{R_0} \cap \mathcal{U}_Q^-$. Then*

$$e_{Q,R_0}H^q(U_{R_0}(p^r), M) = H^{q-q'_0}(\mathcal{U}_{R_0} \cap M_Q, M)$$

*with $q'_0 = |R_0(p) \cap R_Q|$ and $q_0 = |R_0(p)| - |R_0(p) \cap R_Q^-|$; thus this group is trivial for $q$ outside $[q'_0, q_0]$.*

**Remark:** If $R_0(p) \subset \Delta_Q$ there is nothing to prove. If $R_0(p) \subset R_Q \cup R_Q^-$ then the theorem states that $e_0 H^q(U_{R_0}(p^r), M) \neq 0$ only for $q = q_0 = q'_0$ and in this case $e_0 H^q(U_{R_0}(p^r), M) = M$.

**Proof:** Let us prove this lemma by induction on the cardinality of $R_0$, using the $p$-ordinary Hochschild-Serre spectral sequence. In the general case, let $\alpha \in R_0$ such that $U_{\{\alpha\}}(p^r)$ is contained in the center of $U_{R_0}(p^r)$. Then as consequence of Theorem 3.6 of [16] (cf p. 293-297) one has

$$e_{Q,\{\alpha\}}H^q(U_{\{\alpha\}}(p^r), M) = \begin{cases} H^{q-q'_{\{\alpha\}}}(\mathcal{U}_{\{\alpha\}} \cap M_Q, M) & \text{if } q \in [q'_{\{\alpha\}}, q_{\{\alpha\}}], \\ 0 & \text{else.} \end{cases}$$

By the induction hypothesis and the spectral sequence above with $R_1 = \{\alpha\}$, we get our result.∎

Let $U_s(\Gamma') = \Gamma' \cap \mathcal{U}_{\mathbf{P}}$ for $s = s_{\mathbf{P}}$. If we apply the previous lemma to our situation we get:

**Corollary 4.3** *Let $w \in \underline{W}_{Q\Sigma} = W_Q \backslash W_p / W_\Sigma$, $s \in Cusp_{\Sigma,w}(\Gamma_1(p^r))$ and $M$ be a trivial finite or cofinite $G(\mathbf{Q})$ − module. Then*

$$e_{Q,w(R_\Sigma)} H^q(U_s(\Gamma_1(p^r)), M) = H^{q-q'_{\Sigma,w}}(\mathcal{U}_{w(R_\Sigma)} \cap M_Q, M)$$

*where $q'_{\Sigma,w} = |R_Q \cap w(R_\Sigma)|$ and $q_{\Sigma,w} = |R_\Sigma| - |R_Q^- \cap w(R_\Sigma)|$; thus these groups are trivial for $q$ outside $[q'_{\Sigma,w}, q_{\Sigma,w}]$.*

One can prove easily the following lemma:

**Lemma 4.8** *Let $\xi$ with $t^{-1}\xi \equiv 1 \mod p^r$ and $v(\alpha(t)) > (r-s)e_v \geq 0$ for all $\alpha \in R_Q^-$ and $v(\alpha(t)) = 0$ for all $\alpha \in \Delta_Q$. Then for any $R_0$, the canonical reduction map induces the following isomorphism:*

$$U_{R_0}(p^s)/(\xi^{-1}U_{R_0}(p^r)\xi \cap U_{R_0}(p^s)) \cong \overline{\mathcal{U}_{R_0 \cap R_Q^+}}/(\overline{t^{-1}\mathcal{U}_{R_0 \cap R_Q^+}t \cap \mathcal{U}_{R_0 \cap R_Q^+}})$$

*where the overline means the reduction modulo $p^s$ and where $R_0 \cap R_Q^+$ stands for an abbreviation of $R_0(p) \cap R_Q^+$.*

**Lemma 4.9** *For $r \geq s \geq 1$ and any $\Gamma$-module $M$, the restriction map induces a canonical isomorphism:*

$$e_Q H^\bullet(U_{R_0}(p^s), M) \cong e_Q H^\bullet(U_{R_0}(p^r), M).$$

**Proof** By the lemma 4.8, the double classes $[U_{R_0}(p^s)\xi U_{R_0}(p^s)], [U_{R_0}(p^r)\xi U_{R_0}(p^r)]$ and $T = [U_{R_0}(p^r)\xi U_{R_0}(p^s)]$ have the same representants. Thus we have the canonical commutative diagram:

$$\begin{array}{ccc} H^\bullet(U_{R_0}(p^s), M) & \longrightarrow & H^\bullet(U_{R_0}(p^r), M) \\ \downarrow [U_{R_0}(p^s)\xi U_{R_0}(p^s)] & \swarrow T & \downarrow [U_{R_0}(p^r)\xi U_{R_0}(p^r)] \\ H^\bullet(U_{R_0}(p^s), M) & \longrightarrow & H^\bullet(U_{R_0}(p^r), M) \end{array}$$

with $T = [U_{R_0}(p^r)\xi U_{R_0}(p^s)]$. The idempotent associated with $T$ gives the inverse of the restriction map we wanted.∎

We now make use of the lemma 4.9 in order to compute the cohomology of a unipotent subgroup for the module $M = L(\rho \otimes \chi; A)$ where $A = A_\infty = K/\mathcal{O}$ or $A = A_\alpha = p^{-\alpha}\mathcal{O}/\mathcal{O}$ (see Section 3.1).

**Proposition 4.2** *Let $s \in Cusp_{\Sigma,w}(\Gamma_i(p^r))$ for $i = 0$ or $1$. Then we have the isomorphism:*

$$e_{Q,w(R_\Sigma)} H^q(U_s(\Gamma_i(p^r)), L(\rho \otimes \chi; A_\alpha)) = H^{q-q'_{\Sigma,w}}(\mathcal{U}_{w(R_\Sigma)} \cap M_Q, V(\rho \otimes \chi; A_\alpha))$$

*and in particular is zero for $q \notin [q'_{\Sigma,w}, q_{\Sigma,w}]$. Moreover this isomorphism is equivariant for the action of $M_Q \cap M_\Sigma^w$.*

**Proof** By arguments similar to those used in [17] and section 4.2 (cf. Cor. 4.1), , we can prove for $r' \gg r$ that

$$e_{Q,w(R_\Sigma)}H^q(U_s(\Gamma_i(p^{r'})), L(\rho \otimes \chi; A_r)) \cong$$
$$e_{Q,w(R_\Sigma)}H^q(U_s(\Gamma_i(p^{r'})), V(\rho \otimes \chi; A_r))$$

Moreover this isomorphism is equivariant for the action of $M_Q$. If $r' \gg r$, we can assume that $U_s(\Gamma_i(p^{r'})) \cap U_Q^-$ acts trivially on $V(\rho \otimes \chi) \otimes A_r$. Therefore considering together the last isomorphism and lemma 4.3, we have for $r' \gg r$:

$$e_{Q,w(R_\Sigma)}H^q(U_s(\Gamma_i(p^{r'})), L(\rho \otimes \chi; A_r)) =$$
$$H^{q-q'_{\Sigma,w}}(M_Q \cap \mathcal{U}_{w(R_\Sigma)}, V(\rho \otimes \chi) \otimes A_r)$$

On the other hand, by the lemma 4.9, we see that both following cohomology groups are independent of $r' \geq r$.

$$e_{Q,w(R_\Sigma)}H^q(U_s(\Gamma_i(p^r)), L(\rho \otimes \chi; A_r)) \cong e_{Q,w(R_\Sigma)}H^q(U_s(\Gamma_i(p^{r'})), L(\rho \otimes \chi; A_r))$$

we deduce our result from these last two isomorphisms. ∎

### 4.4.3 The spectral sequence associated to the Levi decomposition

For any parabolic subgroup $\mathbf{P}$ of $G(\mathbf{Q})$, we write $U_\mathbf{P}$ for its unipotent radical and let $\mathbf{P} = M_\mathbf{P} U_\mathbf{P}$ be a Levi decomposition of $\mathbf{P}$. We set $M_\mathbf{P}^1$ the derived subgroup of $M_\mathbf{P}$ and $\mathbf{P}^1 = M_\mathbf{P}^1 U_\mathbf{P}$. Let $U_s(\Gamma') = \Gamma' \cap U_\mathbf{P}$ for $s = s_\mathbf{P}$. Then we have an exact sequence:

$$1 \to U_s(\Gamma') \to \Gamma' \cap \mathbf{P} \to M_s(\Gamma') \to 1.$$

where $M_s(\Gamma')$ is an arithmetic subgroup of $M_\mathbf{P}(\mathbf{Q})$. Note that when $\mathbf{P}$ is maximal, the derived subgroup $M'_\mathbf{P}(\mathbf{Q})$ of its Levi is isomorphic to $SL(2,F)$.

Besides, by the Hochschild-Serre spectral sequence we have:

$$H^p(M_{\mathbf{P}_{i,w}}(\Gamma'), H^q(U_{\mathbf{P}_{i,w}}(\Gamma'), M)) \Rightarrow H^{p+q}(\mathbf{P}_{i,w} \cap \Gamma', M)$$

Let

$$t_1 = \prod_{\alpha \in w(\Delta_\Sigma) \cap R_Q^-} \alpha^\vee(p) \text{ and } t' = \prod_{\alpha \in w(R_\Sigma) \cap R_Q^-} \alpha^\vee(p)$$

where $\Delta_\Sigma$ is the set of roots of $M_\Sigma$. Let $\eta_1$ resp. $\eta'$ elements of $H(\mathbf{Q})$ defined as in Section 2.4 such that $(\mathbf{P}_{i,w} \cap \Gamma', U_{\mathbf{P}_{i,w}}(\Gamma'), \eta_1\eta', \eta_1)$ is an admissible quadruple in the sense of section 5.4.1; therefore, we can take the ordinary part of this spectral sequence and get

$$e_{Q \cap M_\Sigma^w} H^p(M_{\mathbf{P}_{i,w}}(\Gamma'), e_Q H^q(U_{\mathbf{P}_{i,w}}(\Gamma'), M) \quad (4.9)$$
$$\Rightarrow e_Q H^{p+q}(\mathbf{P}_{i,w} \cap \Gamma', M)$$

**Definition 4.4** *We set:*

$$V_{\Sigma,w,q}(\rho) = H^{q-q'_{\Sigma,w}}(M_Q \cap \mathcal{U}_{w(R_\Sigma)}, V(\rho))$$
$$L_{\Sigma,w,q}(\rho) = Ind_{M_Q \cap M_\Sigma^w}^{M_\Sigma^w} V_{\Sigma,w,q}(\rho)$$

For any $w \in W_Q$ and any weight $\lambda \in X$, we set $\underline{w}.\lambda = w(\lambda + \varrho_{M_Q}) - \varrho_{M_Q}$ where $\varrho_{M_Q}$ is the half-sum of the positive roots of $M_Q$. For any algebraic subgroup $M$ of $H$ containing $T_H$ and any dominant weight $\lambda$ we denote by $E_M^\lambda$ the irreducible algebraic representation of $M$ of highest weight $\lambda$.

**Proposition 4.3** *Let $w \in \underline{W}_{Q\Sigma} = W_Q \backslash W_p / W_\Sigma$, then one has natural isogenies:*

*(i)* $V_{\Sigma,w,q}(\rho) \to \bigoplus\limits_{\substack{w_Q \in W_{Q,\Sigma,w} \\ l(w_Q) = q - q'_w}} E_{M_\Sigma^w \cap M_Q}^{\underline{w_Q}.(\lambda_\rho \chi)}(\mathcal{O}) \otimes A$

*(ii)* $L_{\Sigma,w,q}(\rho) \to \bigoplus\limits_{\substack{w_Q \in W_{Q,\Sigma,w} \\ l(w_Q) = q - q'_w}} E_{M_\Sigma^w}^{\underline{w_Q}.(\lambda_\rho \chi)}(\mathcal{O}) \otimes A$

*where $W_{Q,\Sigma,w} = \{w_Q \in W_Q \text{ such that } \Delta_Q^+ \cap w_Q(\Delta_Q^-) \subset w(R_\Sigma) \cap \Delta_Q\}$.*

*Moreover, the isogenies (i) and (ii) are isomorphisms if $p$ is chosen larger than the semisimple weight of $\rho$.*

**Proof:** If $U_{/\mathbf{Z}_p}$ is a connected unipotent subgroup of the unipotent radical of a parabolic group $P$ of $G^1$ with Levi subgroup $M$ and $V$ a finite free $\mathcal{O}$-module with an action of $G^1$ making it an irreducible module with highest weight $\mu \in X^+$, one has

$$H^r(U(\mathbf{Z}_p), V) = H^r(Lie(U), V)$$

and

$$H^r(Lie(U), V) \sim \bigoplus_{v \in W_U, l(v) = r} E_M^{\underline{v}.\mu}$$

where $W_U = \{v \in W \text{ such that } R^+ \cap v(R^-) \subset R_U\}$ and $\sim$ means that there is a natural isogeny. The first identity comes from the isomorphism between $Lie(U)$ and $U$ given by the exponential map and the second one is an integral version of a result of Kostant (see [43] ch. 3.2). It is actually a direct calculation since $Lie(U)$, when non zero, is $r_p$, acting on $V(\rho) = \otimes_{v|p} Sym^n(r_v^2)$. So we apply this result to $H^{q-q'_w}(\mathcal{U}_{w(R_\Sigma)} \cap M_Q, V(\rho))$; this yields (i). Then,(ii) follows from (i) by transitivity of induction.∎

**Theorem 4.1** *For all $i$, $w$ and $\Sigma$, we have the following spectral sequences:*

$$(i) \quad e_{Q \cap M_\Sigma^w} H^p(M_{\mathbf{P}_{i,w}}(\Gamma_0(p^r)); L_{\Sigma,w,q}(\rho \otimes \chi; A_\alpha)) \Longrightarrow$$
$$e_Q H^{p+q}(\mathbf{P}_{i,w} \cap \Gamma_0(p^r), L(\rho \otimes \chi; A_\alpha))$$

$$(ii) \quad e_{Q \cap M_\Sigma^w} H^p(M_{\mathbf{P}_{i,w}}(\Gamma_1(p^\infty)); L_{\Sigma,w,q}(\rho \otimes \chi)) \Longrightarrow$$
$$e_Q H^{p+q}(\mathbf{P}_{i,w} \cap \Gamma_1(p^\infty); L^a(\rho \otimes \chi; A))$$

**Proof:** We apply (4.9) to $M = L(\rho \otimes \chi; A_\alpha)$ and $\Gamma' = \Gamma_?(p^{r'})$ (? = 0 or 1) with $r' \gg r$. Using Proposition 4.2 we obtain a converging spectral sequence:

$$E_2^{p,q}(\Gamma_?(p^{r'})) = e_{Q \cap M_\Sigma^w} H^p(M_{\mathbf{P}_{i,w}}(\Gamma_?(p^{r'})), V_{\Sigma,w,q}(\rho \otimes \chi) \otimes A_\alpha)$$
$$\Rightarrow e_Q H^{p+q}(\mathbf{P}_{i,w} \cap \Gamma_?(p^{r'}), L(\rho \otimes \chi; A_\alpha))$$

By the same argument as in Proposition 3.1 and Corollary 3.2, we see

$$E_2^{p,q}(\Gamma_?(p^{r'})) \cong e_{Q \cap M_\Sigma^w} H^p(M_{\mathbf{P}_{i,w}}(\Gamma_?(p^{r'})), L_{\Sigma,w,q}(\rho) \otimes A_\alpha).$$

Let us prove (i). We take ? = 0 in the above formula; a slight variant of Lemma 3.1 yields:

$$E_2^{p,q}(\Gamma_0(p^{r'})) = E_2^{p,q}(\Gamma_0(p^r))$$
$$e_Q H^{p+q}(\mathbf{P}_{i,w} \cap \Gamma_0(p^{r'}), L(\rho \otimes \chi; A_\alpha)) = e_Q H^{p+q}(\mathbf{P}_{i,w} \cap \Gamma_0(p^r), L(\rho \otimes \chi; A_\alpha)))$$

This is what we wanted.

To prove (ii), we simply take the limit over $r'$ and over $r$ and notice that forming spectral sequences commutes to inductive limits.∎

**Definition 4.5** *For $Q = B$, for any $w \in W_p/W_\Sigma$, note that $q_{\Sigma,w} = q'_{\Sigma,w}$; let $n_{\Sigma,w}$ be this common value and let us abbreviate $L_{\Sigma,w,n_w}(\rho \otimes \chi)$ as $L_{\Sigma,w}(\rho \otimes \chi)$.*

**Corollary 4.4** *We have the following equality:*

$$e_B H^p(\mathbf{P}_{i,w} \cap \Gamma_0(p^r), L(\rho \otimes \chi; A)) =$$
$$e_{B \cap M_\Sigma^w} H^{p-n_{\Sigma,w}}(M_{\mathbf{P}_{i,w}}(\Gamma_0(p^r)), L_{\Sigma,w}(\rho \otimes \chi))$$

**Corollary 4.5** *The following canonical map has finite kernel and cokernel*

$$e_Q H^p(\mathbf{P}_{i,w} \cap \Gamma_0(p^r), L(\rho \otimes \chi; A)) \to \oplus_{k+l=p} e_{Q \cap M_\Sigma^w} H^k(M_{\mathbf{P}_{i,w}}(\Gamma_0(p^r)), L_{\Sigma,w,l}(\rho \otimes \chi; A))$$

**proof** It is well-known that the spectral sequence (4.9) over $\mathbf{C}$, degenerates when the weight $\lambda_\rho \otimes \chi$ is regular (by [29] Th.2.7, p.51). The corollary follows easily from this observation.∎

### 4.4.4 The center of the Levi subgroup

¿From now on, we assume that $P_\Sigma$ is maximal and, for simplicity, that $Q$ is rational (i.e. $Q_v$ is independent of $v$). In order to complete the computation of the ordinary cohomology of the strata, we shall consider the action of the center of the Levi component for each cusp. Let us first pose a definition.

**Definition 4.6** *A weight $\lambda = (x_\sigma, y_\sigma; z_\sigma)_{\sigma \in I_F}$ is called separable (resp. sufficiently separable) if for any $\sigma, \sigma' \in I_F$ with $\sigma \neq \sigma'$, we have $x_\sigma \neq y_{\sigma'}$ and $x_\sigma + y_\sigma \neq x_{\sigma'} - y_{\sigma'}$ (resp. $|x_\sigma - y_{\sigma'}| > 3$ and $|x_\sigma + y_\sigma - (x_{\sigma'} - y_{\sigma'})| > 3$).*

One can easily see that if $Inf\{x_\sigma\ ; \sigma \in I_F\} > Sup\{y_\sigma\ ; \sigma \in I_F\}$ then $\lambda$ is separable and of course regular. The utility of the above definition appears in the following lemma:

**Lemma 4.10** *Let $\lambda$ be a separable (resp. sufficiently separable) weight. Let $w \in W^\Sigma$ for a fixed type of parabolic $\Sigma$. Let $v, v' \in S_p$ with $\sigma \neq \sigma'$; assume that $(w.\lambda)|_{Z_\Sigma(F_v)} = (w.\lambda)|_{Z_\Sigma(F_{v'})}$ (resp. $|(w.\lambda)|_{Z_\Sigma(F_v)} - (w.\lambda)|_{Z_\Sigma(F_{v'})}| \leq 3$), then we have $w(v) = w(v')$.*

For any $s = s_\mathbf{P} \in Cusp_{\Sigma,w}(\Gamma')$, we consider $\delta \in H(\mathbf{Q})$ such that $\delta \mathbf{P} \delta^{-1} = P_\Sigma$ and $\delta \equiv w$ mod. $p^r$ (see Lemma 4.3 above). Then we set

$$E_s(\Gamma') = \delta(\mathbf{P} \cap \Gamma')\delta^{-1} \cap Z_\Sigma(\mathbf{Q})$$

where $Z_\Sigma$ is the center of $M_\Sigma$. Therefore $E_s(\Gamma')$ can be seen as a subgroup of $\mathbf{r}^\times$.

**Lemma 4.11** *(i) Via this identification, $E_s(\Gamma')$ acts on the maximal p-divisible subgroup of $L_{\Sigma,w,q}(\rho \otimes \chi)$ by the characters*

$$\omega_{\Sigma,w,w_Q} = w^{-1}\underline{w_Q}.(\lambda_\rho + \chi) + \sum_{\alpha \in R_\Sigma \cap w^{-1}(R_Q)} \alpha$$

*where $w_Q$ runs in $W_{Q,\Sigma,w}$ and $q - q'_{w,\Sigma} = l(w_Q)$.*

*(ii) If $\omega_{\Sigma,w,w_Q}$ is trivial on $E_s(\Gamma_1(p^r))$ and $\lambda_\rho \otimes \chi$, is regular and sufficiently separable, then $w(v) \in W(v)/W_\Sigma(v)$ is independent from the v's dividing $p$ and $q, q_w, q'_w$ are multiples of $d$.*

**proof** The first point (i) follows from Proposition 4.3. Assertion (ii) follows from the remark that an algebraic character is trivial on a subgroup of finite index of the units of $F$ if and only if it is a multiple of the norm. Then $w^{-1}(\underline{w_Q}.\lambda_\rho + \chi) + \sum_{\alpha \in R_\Sigma \cap w^{-1}(R_Q)} \alpha)(v)$ is independent of $v|p$. Now by Lemma 4.10, the character $\lambda_\rho \chi$ regular and sufficiently separable has the following property: given two places $v$ and $v'$ above $p$, if $w(v) \neq w(v')$ then the characters $(\lambda_\rho \chi)^w_\sigma(v)$ and $(\lambda_\rho \chi)^w_{v'}(v')$ are different and even far one from the other. The lemma follows.∎

# 5 p-Ordinary cohomology of the boundary of the Borel-Serre compactification

## 5.1 Review on the Borel-Serre compactification

In this section, we recall some classical facts about the Borel-Serre compactification constructed in [2]. Let G be a reductive algebraic over $\mathbf{Q}$ with positive $\mathbf{Q}$-rank. Let $X = G(\mathbf{R})/Z(\mathbf{R})^0 K_\infty$ with $K_\infty$ a maximal compact subgroup of $G(\mathbf{R})$ and $Z(\mathbf{R})^0$ the connected componnent of the center of $G(\mathbf{R})$. Borel and Serre constructed then a contractible topological space $\tilde{X}$ with an action of $G(\mathbf{Q})$ such that for any arithmetic subgroup $\Gamma$, $\Gamma\backslash\tilde{X}$ is a compactification of $\Gamma\backslash X$, whose boundary has the homotopy type of the quotient by $\Gamma$ of the Tits building of parabolic $\mathbf{Q}$-subgroups of $G$. For any parabolic $\mathbf{Q}$-subgroup $Q$, we let $A_Q$ be the identity component of the real points of the center of a Levi subgroup of $Q$. Then we have $\tilde{X} = \bigcup_Q e(Q)$ with $e(Q) = A_Q\backslash X$. Therefore, we have

$$\Gamma\backslash\tilde{X} = \bigsqcup_{\substack{\Gamma - conjugacy\ classes \\ of\ parabolic \\ subgroups \\ of\ G}} (\Gamma \cap Q)\backslash e(Q)$$

$$\partial(\Gamma\backslash\tilde{X}) = \bigsqcup_{\substack{\Gamma - conjugacy\ classes \\ of\ proper\ parabolic \\ subgroups\ of\ G}} (\Gamma \cap Q)\backslash e(Q)$$

Let $\mathfrak{P}_G(\Gamma)$ be the set of $\Gamma$-conjugacy classes of proper maximal parabolic $\mathbf{Q}$-subgroups of $G$ and let $X_Q(\Gamma)$ the closure of $(\Gamma \cap Q)\backslash e(Q)$ for any parabolic $\mathbf{Q}$-subgroup $\mathcal{Q}$. Then $X_Q(\Gamma) = \bigcup_{P \subset Q} (\Gamma \cap P)\backslash e(P)$ and the family $(X_Q(\Gamma))_{Q \in \mathfrak{P}_G(\Gamma)}$ forms a finite cover of the boundary $\partial(\Gamma\backslash\tilde{X})$ and we can calculate its cohomology by using the spectral sequence of Leray (cf. Th 5.24 p.209 of [11]): the term $E_1^{p,q}$ is

$$\bigoplus_{\substack{(P_1,\ldots,P_p) \\ with\ P_i \in \mathfrak{P}_G(\Gamma)}} H^q(\cap_{i=1}^p X_{P_i}(\Gamma), \mathcal{F})$$

for any sheaf $\mathcal{F}$. For any subsequence $(Q_1,\ldots,Q_s)$ of $(P_1,\ldots,P_r)$ we note $Res_{(Q_1,\ldots,Q_s)(P_1,\ldots,P_r)}$ the obvious map

$$H^q(\cap_{i=1}^s X_{Q_i}(\Gamma), \mathcal{F}) \to H^q(\cap_{i=1}^r X_{P_i}(\Gamma), \mathcal{F}).$$

Then differential map from $E_1^{p,q}$ to $E_1^{p+1,q}$ is given by:

$$d(\alpha) = \sum_{Q \in \mathfrak{P}_G(\Gamma)} (-1)^i Res_{(P_1,\ldots,P_p)(P_1,\ldots,P_{i-1},Q,P_i,\ldots,P_p)}(\alpha)$$

I f $G = Res_{\mathbf{Q}}^F G$ with $rk_F G = r$ with $G_{/F}$ split, we can see easily that $E_2^{p,q} = 0$ if $p > r = rank_{\mathbf{Q}}(G)$ and therefore this spectral sequence degenerates at least in $E_r$.

Now we consider the case $r = 2$; we fix a Borel subgroup $B$ and let $P_1$ and $P_2$ the two maximal parabolic subgroups containing $B$. Then the spectral sequence degenerates in $E_2$ and may be seen as the Mayer-Vietoris sequence:

$$\ldots \to H^{q-1}(\partial_B(\Gamma\backslash\tilde{X});\mathcal{F}) \to H^q((\partial(\Gamma\backslash\tilde{X});\mathcal{F}) \to$$
$$\to H^q((\partial_{P_1}(\Gamma\backslash\tilde{X});\mathcal{F}) \oplus H^q((\partial_{P_2}(\Gamma\backslash\tilde{X});\mathcal{F}) \to H^q((\partial_B(\Gamma\backslash\tilde{X});\mathcal{F}) \to \ldots$$

where $H^q((\partial_{\mathbf{P}_i}(\Gamma\backslash\tilde{X});\mathcal{F})$ (respectively $H^q((\partial_B(\Gamma\backslash\tilde{X});\mathcal{F}))$ is the sum of the groups $H^q(X_P(\Gamma);\mathcal{F})$ when $P$ varies in the set of $\Gamma$-conjugacy classes of maximal parabolic subgroups of $G$ conjugated with $P_i$ (respectively Borel subgroups of $G$).

## 5.2 Ordinary cohomology of the boundary

For any $\Gamma$-representation M, we can construct a sheaf $\tilde{M}$ over $\Gamma\backslash\tilde{X}$ as the sheaf of locally constant functions on $\Gamma\backslash\tilde{X}$ with values in M. Moreover for any $g \in G(\mathbf{Q})$ which acts on M, it is well known that we can define an action of the double class $[\Gamma g\Gamma]$ on the cohomology groups $H^q(\Gamma\backslash\tilde{X};\tilde{M}), H^q(\partial(\Gamma\backslash\tilde{X});\tilde{M})$ and $H^q(\cap_{i=1}^p X_i(\Gamma);\tilde{M})$ such that the restriction maps $Res_{(Q_1,\ldots,Q_s)(P_1,\ldots,P_r)}$ commute with this action. From this, we deduce that the action of the Hecke operators commute with the differential map of the previous spectral sequence. In particular, this point proves that $e_Q E_1^{p,q} \Rightarrow e_Q H^{p+q}(\partial(\Gamma\backslash\tilde{X});\tilde{M})$. Therefore, when the rank is 2, the above Mayer-Vietoris sequence induces, by application of the "$p$-ordinary" idempotent, a (Hecke-equivariant) long exact sequence for $\Gamma'$:

$$\ldots \to e_Q G_B^{p-1}(\Gamma';M) \to$$
$$\to e_Q H^p(\partial(\Gamma'\backslash\tilde{X});\tilde{M}) \to e_Q G_{P^*}^p(\Gamma';M) \oplus e_Q G_P^p(\Gamma';M) \xrightarrow{r^p} e_Q G_B^p(\Gamma';M) \to \ldots$$

**Theorem 5.1** *(Independence of the weight) With the previous Notation, we have a canonical isomorphism:*

$$(i) \; e_Q H^\bullet(\partial(\Gamma_1(p^\infty)\backslash\tilde{X}), L^a(\rho \otimes \chi; A)) = \varinjlim_r e_Q H^\bullet(\partial(\Gamma_1(p^r)\backslash\tilde{X}), L^a(\rho \otimes \chi; A))$$
$$\cong e_Q H^\bullet(\partial(\Gamma_0(p)\backslash\tilde{X}), \mathcal{C}(\rho^1; A))$$
$$(ii) \; e_Q H_!^\bullet \Gamma_1(p^\infty), L^a(\rho \otimes \chi; A)) = \varinjlim_r e_Q H_!^\bullet \Gamma_1(p^r), L^a(\rho \otimes \chi; A))$$
$$\cong e_Q H_!^\bullet(\Gamma_0(p), \mathcal{C}(\rho^1; A))$$

**Proof:** The first point is deduced from the ordinary Mayer-Vietoris exact sequence for $M = L^a(\rho \otimes \chi; A)$ and $M = \mathcal{C}(\rho^1; A)$ and the corollary 4.2. The second point follows from the first one and the corollary 3.2.∎

**Remark:** This theorem is true even if we do not assume that the rank is 2. Indeed, we have the isomorphism on the term $E_1^{p,q}$ of the spectral sequences abutting on the boundary cohomology $e_Q H^p(\partial(\Gamma_1(p^\infty)\backslash \tilde{X}), L^a(\rho \otimes \chi; A))$ and $e_Q H^p(\partial(\Gamma_0(p)\backslash \tilde{X}), \mathcal{C}(\rho^1; A))$, since the corollary 4.2 does not require any assumption on the rank of our group. This remark holds as well for the following lemma, deduced from lemma 4.5.

**Lemma 5.1** *For any $\chi$ of level $p^r$, we have a canonical isomorphism*

$$e_Q H^p(\partial(\Gamma_0(p^r)\backslash \tilde{X}), L^a(\rho \otimes \chi; A)) \cong e_Q H^p(\partial(\Gamma_0(p^r)\backslash \tilde{X}), \mathcal{C}(\rho^1; A)[\omega_\chi])$$

One can decompose the cohomology groups $e_Q G_\Sigma^p$ in terms of the $e_Q G_{\Sigma,w}^p$'s (where $w$ runs over the Weyl types). Let us examine whether the morphism $r^p$ preserves these decompositions. Consider $\pi_\Sigma$ the canonical surjection from $W_Q\backslash W_p$ onto $\underline{W}_{Q\Sigma} = W_Q\backslash W_p/W_\Sigma$. If we consider for all $(w, w') \in \underline{W}_{QP} \times \underline{W}_{QP^*}$ the map:

$$r_{(w,w')}^p : \left(e_Q G_{P,w}^p(\Gamma'; M) \oplus e_Q G_{P^*,w'}^p(\Gamma'; M)\right) \to \bigoplus_{v \in \pi_P^{-1}(w) \cup \pi_{P^*}^{-1}(w')} e_Q G_{B,v}^p(\Gamma'; M).$$

then we have

$$r^p = \sum_{(w,w') \in \underline{W}_{QP} \times \underline{W}_{QP^*}} r_{(w,w')}^p.$$

We want now to compute the $r^p$'s in terms of the Hochschild-Serre spectral sequences of Th.4.1. More precisely, consider the map

$$r_{i,w}^p : e_Q H^p(\mathbf{P}_{i,w} \cap \Gamma_0(p^r), L^a(\rho \otimes \chi, A)) \to$$

$$\to \bigoplus_{\substack{s_{\mathbf{B}_k} \in Cusp_B(\Gamma) \\ \mathbf{P}_i \supset \mathbf{B}_k}} \bigoplus_{v \in \pi_\Sigma^{-1}(w)} e_Q H^p(\mathbf{B}_{k,v} \cap \Gamma_0(p^r), L^a(\rho \otimes \chi, A))$$

There are similar maps at the level of the $E_2^{p,q}$-terms of spectral sequences of Th.4.1 abutting to the source, resp. target of $r_{i,w}^{p+q}$. Let us first consider the $\Gamma_0$-case:

$$R_{i,w}^{p,q} : e_{Q \cap M_\Sigma^w} H^p(M_{\mathbf{P}_{i,w}}(\Gamma_0(p^r)); L_{\Sigma,w,q}(\rho \otimes \chi)) \longrightarrow$$

$$\bigoplus_{\substack{s_{\mathbf{B}_k} \in Cusp_B(\Gamma) \\ \mathbf{P}_i \supset \mathbf{B}_k}} \bigoplus_{v \in \pi_\Sigma^{-1}(w)} \oplus_{t=0}^p H^{p-t}(M_{\mathbf{B}_{k,v}} \cap \Gamma_0(p^r), L_{\Sigma,v,q+t}(\rho \otimes \chi)) \tag{5.1}$$

It is a theorem of Hida (see Th.3.12 of [16], proven there in the $SL_2$-case, but easily generalized to the $GL_2$-case) that the right hand side module of 5.1 is canonically isomorphic to

$$e_{Q \cap M_\Sigma^w} H_\partial^p(M_{\mathbf{P}_{i,w}}(\Gamma_0(p^r)), L_{\Sigma,w}(\rho \otimes \chi))$$

Let us set

$$R_{\Sigma,w}^{p,q} = \bigoplus_{s_{\mathbf{P}_i} \in Cusp_\Sigma(\Gamma)} R_{i,w}^{p,q}$$

and let us introduce

$$r'^p_{(w,w')} = \sum_{q=q'_{P,w}}^{q_{P,w}} R_{P,w}^{p-q,q} + \sum_{q=q'_{P^*,w'}}^{q_{P^*,w'}} R_{P^*,w'}^{p-q,q}$$

and

$$r'^p = \sum_{(w,w') \in \underline{W}_{QP} \times \underline{W}_{QP^*}} r'^p_{(w,w')}$$

**Proposition 5.1** *With the previous Notation, $r'^p$ and $r^p$ induce the same map from the $Gr^*(e_Q G_{P^*}^p(\Gamma_0(p^r); M) \oplus e_Q G_P^p(\Gamma_0(p^r); M))$ into $Gr^* e_Q G_B^p(\Gamma_0(p^r); M)$*

## 5.3 The case $(Sp_4)_{/\mathbf{Q}}$:

Let $\Gamma$ be an arithmetic subgroup of level prime to $p$ and without torsion.

Let $e_Q H_\partial^q(\Gamma_1(p^\infty); L^a(\rho; A)) = \varinjlim_r e_Q H^q(\partial(\Gamma_1(p^r) \backslash \tilde{X}); L^a(\rho; A))$. Then, the control theorem for the boundary cohomology is the following:

**Theorem 5.2** *For any arithmetic character $\chi$ of $C_Q(\mathbf{Z}_p)$ we have :*

$$e_Q H^q(\partial(\Gamma_0(p^r) \backslash \tilde{X}), L^a(\rho \otimes \chi; A)) \cong e_Q H^q(\partial(\Gamma_0(p) \backslash \tilde{X}), \mathcal{C}(\rho^1; A))[\omega_\chi]$$

*Moreover this group is cofree in the case $Q = B$.*

This theorem will be a direct consequence of the calculations of subsections 5.3.1-5.3.2 below, together with Hida theory for $SL(2, \mathbf{Q})$ (Theorem 1.9 [14]).

**Corollary 5.1** *For any arithmetic character $\chi$ of $C_{M'}$ which is regular, dominant with respect to $\rho$, there is a canonical map with finite kernel or cokernel:*

$$H^3_{!,Q-no}(\Gamma_0(p^r), L^a(\rho \otimes \chi; A)) \to H^3_{!,Q-no}(\Gamma_1(p^\infty), L^a(\rho; A))[\omega_\chi]$$

Recall that $\alpha_1$ (resp. $\alpha_2$) denotes the short root (resp. the long one) and $s_i$ the element of the Weyl group associated with $\alpha_i$. We summarize the results in the following tables:

### 5.3.1 The case $Q = B$

| $P_\Sigma$ | $R_\Sigma$ | $\underline{W}_{Q\Sigma}$ |
|---|---|---|
| $P$ | $\alpha_2, \alpha_1 + \alpha_2, 2\alpha_1 + \alpha_2$ | $id, -id, s_2, -s_2$ |
| $P^*$ | $\alpha_1, \alpha_1 + \alpha_2, 2\alpha_1 + \alpha_2$ | $id, -id, s_1, -s_1,$ |
| $B$ | $\alpha_1, \alpha_2, \alpha_1 + \alpha_2, 2\alpha_1 + \alpha_2$ | $id, -id, s_1, -s_1, s_2, -s_2, s_1s_2, s_2s_1$ |

The three following tables give the values of $n_w = q_w = q'_w$ for the different $\Sigma$:

| $P_\Sigma = P$ | |
|---|---|
| $w$ | $n_w$ |
| $id$ | 3 |
| $-id$ | 0 |
| $s_2$ | 2 |
| $-s_2$ | 1 |

| $P_\Sigma = P^*$ | |
|---|---|
| $w$ | $n_w$ |
| $id$ | 3 |
| $-id$ | 0 |
| $s_1$ | 2 |
| $-s_1$ | 1 |

| $P_\Sigma = B$ | |
|---|---|
| $w$ | $n_w$ |
| $id$ | 4 |
| $-id$ | 0 |
| $s_2$ | 3 |
| $-s_2$ | 1 |
| $s_1$ | 3 |
| $-s_1$ | 1 |
| $s_1s_2$ | 2 |
| $s_2s_1$ | 2 |

$M'_Q = T'$ is the subgroup of diagonal matrix isomorphic to $\mathbb{G}_m \times \mathbb{G}_m$ and $\rho$ is nothing but an algebraic $\mathcal{O}$-valued character. We know by the corollary 4.4, we have for any dominant $\mathcal{O}$-valued arithmetic character $\chi$:

$$e_B G_\Sigma^q(\Gamma_0(p^r); L^a(\rho\chi; A)) = \bigoplus_{s_{\mathbf{P}_i} \in Cusp_\Sigma(\Gamma)} \bigoplus_{w \in W/W_\Sigma} e_{B \cap M_\Sigma^w} H^{q-n_w}(M_{\mathbf{P}_{i,w}}(\Gamma_0(p^r)); L_{\Sigma,w}(\rho \otimes \chi))$$

Now, if $P_\Sigma \neq B$ the subgroups $M_{\mathbf{P}_{i,w}}(\Gamma_1(p^r))$ are congruence subgroups of $SL_2(\mathbf{Z})$; we know by results of Hida that the ordinary part of their cohomology is non trivial only in degree 1. We set

| $q$ | $w_{P,q}$ | $w_{P^*,q}$ |
|---|---|---|
| 1 | $-id$ | $-id$ |
| 2 | $-s_2$ | $-s_1$ |
| 3 | $s_2$ | $s_1$ |
| 4 | $id$ | $id$ |

Given a pair $(\Sigma, q)$, where $q \in [1, 4]$, let us abbreviate the unique element $w_{\Sigma,q}$ defined in the tables above by $w$. We see from these remarks and the tables that:

$$e_B G_\Sigma^q(\Gamma_0(p^r), L^a(\rho \otimes \chi; A)) = e_B G_{\Sigma,w}^q(\Gamma_0(p^r), L^a(\rho \otimes \chi; A)) =$$
$$\bigoplus_{s_{\mathbf{P}_i} \in Cusp_\Sigma(\Gamma)} e_{B \cap M_\Sigma^w} H^1(M_{\mathbf{P}_{i,w}}(\Gamma_0(p^r)), L^a(\rho \otimes \chi; A))$$

and

$$e_B G_\Sigma^q(\Gamma_0(p^r), L^a(\rho \otimes \chi; A)) = \bigoplus_{s_{\mathbf{P}_i} \in Cusp_\Sigma(\Gamma)} e_{B \cap M_\Sigma^w} H^1(M_{\mathbf{P}_{i,w}}(\Gamma_0(p^r)), L_{\Sigma,w}(\rho \otimes \chi; A))$$

and

$$e_B G_\Sigma^q(\Gamma_0(p^r), L^a(\rho \otimes \chi; A)) = 0 \text{ for } q = 0, 5$$

For $P_\Sigma = B$, we have similarly:

$$e_B G_B^q(\Gamma_0(p^r); L^a(\rho \otimes \chi; A)) = \bigoplus_{\substack{s_{B_k} \in \\ Cusp_B(\Gamma)}} \bigoplus_{w,\, n_w = q} L_{B,w}(\rho \otimes \chi)$$

In order to study the properties for the boundary cohomology, we consider the maps $r_q(w, w')$. By the tables and the previous discussion, we can see that for each $q$, the only couple $(w, w')$ for which $r_{(w,w')}^q$ is non trivial is $(w_{P,q}, w_{P^*,q})$; moreover,

$$r_{(w_{P,q}, w_{P^*,q})}^q = R_{P, w_{P,q}}^{1,q} + R_{P^*, w_{P^*,q}}^{1,q}$$

Moreover for $q = 2, 3$, this is a direct sum. Noting $\phi_P^q = R_{P, w_{P,q}}^{1,q}$ and $\phi_{P^*}^q = R_{P^*, w_{P^*,q}}^{1,q}$, we have therefore by the ordinary Mayer-Vietoris exact sequence:

$$e_B H^1(\partial(\Gamma_0(p^r) \backslash \tilde{X}); M) = ker(\phi_P^1 + \phi_{P^*}^1)$$

$$0 \to coker(\phi_P^1 + \phi_{P^*}^1) \to e_B H^2(\partial(\Gamma_0(p^r) \backslash \tilde{X}); M) \to ker\phi_P^2 \oplus ker\phi_{P^*}^2 \to 0$$

$$0 \to coker\phi_P^2 \oplus coker\phi_{P^*}^2 \to e_B H^3(\partial(\Gamma_0(p^r) \backslash \tilde{X}); M) \to ker\phi_P^3 \oplus ker\phi_{P^*}^3 \to 0$$

$$0 \to coker\phi_P^3 \oplus coker\phi_{P^*}^3 \to e_B H^4(\partial(\Gamma_0(p^r) \backslash \tilde{X}); M) \to ker(\phi_P^4 + \phi_{P^*}^4) \to 0$$

Recall that we put $w$ for the unique $w_{\Sigma,q}$ defined in the tables above. By looking at how the restriction occurs in the Hochschild-Serre spectral sequence

$$ker\phi_\Sigma^q = \bigoplus_{s_{\mathbf{P}_i} \in Cusp_\Sigma(\Gamma)} eH_{cusp}^1(M_{\mathbf{P}_{i,w}}(\Gamma_0(p^r)), L_{\Sigma,w}(\rho \otimes \chi; A))$$

is cofree over $\mathcal{O}$ and $coker\phi_\Sigma^q = 0$ by classical Hida theory cf. [14].

For $q = 1, 4$, we can see that $ker(\phi_P^q + \phi_{P^*}^q) \cong ker(\phi_P^q) \oplus ker(\phi_{P^*}^q) \oplus e_B G_B^q(\Gamma_0(p^r), L^a(\rho \otimes \chi; A))$. We summarize our results in the following

**Theorem 5.3** *The p-ordinary cohomology of the boundary of the Borel-Serre compactification is described by:*

- $e_B H^q(\partial(\Gamma_0(p^r) \backslash \tilde{X}), L^a(\rho \otimes \chi; A)) = 0$ *for* $q = 0, 5$.

- $e_B H^q(\partial(\Gamma_0(p^r)\backslash \tilde{X}), L^a(\rho \otimes \chi; A)) =$
$$\bigoplus_{\substack{s_{\mathbf{B}_k} \in \\ Cusp_B(\Gamma)}} \bigoplus_{w, n_w = q} L_{B,w}(\rho \otimes \chi; A) \oplus$$
$$\bigoplus_{\Sigma = P, P^*} \bigoplus_{\substack{s_{\mathbf{P}_i} \in \\ Cusp_\Sigma(\Gamma)}} eH^1_{cusp}(M_{\mathbf{P}_{i,w_q}}(\Gamma_0(p^r)), L_{\Sigma, w_q}(\rho \otimes \chi; A))$$

  where in the last part of the sum, for $q = 1, 4$, the element $w_q$ of the Weyl group is defined by $w_1 = -id$ and $w_4 = id$.

- For $q = 2, 3$, abbreviating again $w_{\Sigma, q}$ as $w$, we have

$$e_B H^q(\partial(\Gamma_0(p^r)\backslash \tilde{X}), L(\rho \otimes \chi; A)) =$$
$$\bigoplus_{\Sigma = P, P^*} \bigoplus_{\substack{s_{\mathbf{P}_i} \in \\ Cusp_\Sigma(\Gamma)}} eH^1_{cusp}(M_{\mathbf{P}_{i,w}}(\Gamma_0(p^r)), L_{\Sigma, w}(\rho \otimes \chi; A)).$$

Similarly, calculations using the same vanishing results given by the above tables yield a theorem for $\Gamma_1(p^r)$-type groups:

**Theorem 5.4** *The p-ordinary cohomology of the boundary of the Borel-Serre compactification is described as $T(\mathbf{Z}/p^r\mathbf{Z})$-module by:*

- $e_B H^q(\partial(\Gamma_1(p^r)\backslash \tilde{X}); L^a(\rho \otimes \chi; A)) = 0$ for $q = 0, 5$.

- $e_B H^q(\partial(\Gamma_1(p^r)\backslash \tilde{X}), L^a(\rho \otimes \chi; A)) =$
$$\bigoplus_{\substack{s_{\mathbf{B}_k} \in \\ Cusp_B(\Gamma)}} \bigoplus_{w, n_w = q} Ind_{\bar{E}_{\mathbf{B}_k}(\Gamma)}^{T(\mathbf{Z}/p^r\mathbf{Z})} L_{B,w}(\rho \otimes \chi; A) \quad \oplus$$
$$\bigoplus_{\Sigma = P, P^*} \bigoplus_{\substack{s_{\mathbf{P}_i} \in \\ Cusp_\Sigma(\Gamma)}} Ind_{T^1_{\Sigma, w_q}(\mathbf{Z}/p^r\mathbf{Z})\bar{E}_{\mathbf{P}_i}(\Gamma)}^{T(\mathbf{Z}/p^r\mathbf{Z})} eH^1_{cusp}(M_{\mathbf{P}_{i,w_q}}(\Gamma_1(p^r)), L_{\Sigma, w_q}(\rho \otimes \chi; A))$$

  where in the last part of the sum, for $q = 1, 4$, the element $w_q$ of the Weyl group is defined by $w_1 = -id$ and $w_4 = id$.

- For $q = 2, 3$, abbreviating again $w_{\Sigma, q}$ as $w$, we have

$$e_B H^q(\partial(\Gamma_1(p^r)\backslash \tilde{X}), L(\rho \otimes \chi; A)) =$$
$$\bigoplus_{\Sigma = P, P^*} \bigoplus_{\substack{s_{\mathbf{P}_i} \in \\ Cusp_\Sigma(\Gamma)}} Ind_{T^{1,w}_\Sigma(\mathbf{Z}/p^r\mathbf{Z})\bar{E}_{\mathbf{P}_i}(\Gamma)}^{T(\mathbf{Z}/p^r\mathbf{Z})} eH^1_{cusp}(M_{\mathbf{P}_{i,w}}(\Gamma_1(p^r)), L_{\Sigma, w}(\rho \otimes \chi; A)).$$

### 5.3.2 The case $Q$ maximal

We begin by $Q = P$. Then we have:

| $P_\Sigma$ | $R_\Sigma$ | $W_{Q\Sigma}$ |
|---|---|---|
| $P$ | $\alpha_2, \alpha_1 + \alpha_2, 2\alpha_1 + \alpha_2$ | $id, -id, s_2$ |
| $P^*$ | $\alpha_1, \alpha_1 + \alpha_2, 2\alpha_1 + \alpha_2$ | $id, -id$ |
| $B$ | $\alpha_1, \alpha_2, \alpha_1 + \alpha_2, 2\alpha_1 + \alpha_2$ | $id, -id, s_2, -s_2$ |

The three following tables give the values of $q_w$ and $q'_w$ for the different $\Sigma$:

| $P_\Sigma = P$ | | |
|---|---|---|
| $w$ | $q'_w$ | $q_w$ |
| id | 3 | 3 |
| -id | 0 | 0 |
| $s_2$ | 1 | 2 |

| $P_\Sigma = P^*$ | | |
|---|---|---|
| $w$ | $q'_w$ | $q_w$ |
| id | 2 | 3 |
| -id | 0 | 1 |

| $P_\Sigma = B$ | | |
|---|---|---|
| $w$ | $q'_w$ | $q_w$ |
| id | 3 | 4 |
| -id | 0 | 1 |
| $s_2$ | 2 | 3 |
| $-s_2$ | 1 | 2 |

Since the spectral sequence of the theorem 4.1 still degenerates in $E_2$, we can make analogous computation and get

**Theorem 5.5** *Let $Q = P$; there exists a filtration of the degree three ordinary boundary cohomology whose associated graded module is given by:*

- $Gr\ e_P H^3(\partial(\Gamma_0(p^r)\backslash \tilde{X}), L^a(\rho \otimes \chi; A)) = \bigoplus\limits_{\substack{s\mathbf{P}_i \in \\ Cusp_P(\Gamma)}}$
  $eH^1_{cusp}(M_{\mathbf{P}_{i,s_2}}(\Gamma_0(p^r)), L_{P,s_2,2}(\rho \otimes \chi)) \oplus H^0(M_{\mathbf{P}_{i,id}}(\Gamma), L_{P,id,3}(\rho \otimes \chi))$
  $\bigoplus\limits_{\substack{s\mathbf{P}_i \in \\ Cusp_{P^*}(\Gamma)}} eH^1_{cusp}(M_{\mathbf{P}_{i,id}}(\Gamma_0(p^r)), L_{P^*,id,2}(\rho \otimes \chi))$

- $Gr\ e_P H^3(\partial(\Gamma_1(p^r)\backslash \tilde{X}), L^a(\rho \otimes \chi; A)) = \bigoplus\limits_{\substack{s\mathbf{P}_i \in \\ Cusp_P(\Gamma)}}$
  $Ind_{i_P(T_P^{1,s_2}(\mathbf{Z}/p^r\mathbf{Z})\bar{E}_{\mathbf{P}_i}(\Gamma))}^{C_P(\mathbf{Z}/p^r\mathbf{Z})} eH^1_{cusp}(M_{\mathbf{P}_{i,s_2}}(\Gamma_1(p^r)), L_{P,s_2,2}(\rho \otimes \chi))$
  $\oplus Ind_{i_P(\bar{E}_{\mathbf{P}_i}(\Gamma))}^{C_P(\mathbf{Z}/p^r\mathbf{Z})} H^0(M_{\mathbf{P}_{i,id}}(\Gamma), L_{P,id,3}(\rho \otimes \chi))$
  $\bigoplus\limits_{\substack{s\mathbf{P}_i \in \\ Cusp_{P^*}(\Gamma)}} Ind_{i_P(T_{P^*}^{1,id}(\mathbf{Z}/p^r\mathbf{Z})\bar{E}_{\mathbf{P}_i}(\Gamma))}^{C_P(\mathbf{Z}/p^r\mathbf{Z})} eH^1_{cusp}(M_{\mathbf{P}_{i,id}}(\Gamma_1(p^r)), L_{P^*,id,2}(\rho \otimes \chi))$

**Remark:** Note in both cases the contribution of full $H^0$ without taking ordinary part and for the level group $\Gamma$ comes from the coincidence $Q = w\mathcal{P}_\Sigma w^{-1}$ for $P_\Sigma = P$ and $w = id$. We observe these $H^0$'s are torsion although it doesn't matter here. We consider now the case $Q = P^*$

| $P_\Sigma$ | $R_\Sigma$ | $W$ |
|---|---|---|
| $P$ | $\alpha_2, \alpha_1 + \alpha_2, 2\alpha_1 + \alpha_2$ | $id, -id$ |
| $P^*$ | $\alpha_1, \alpha_1 + \alpha_2, 2\alpha_1 + \alpha_2$ | $id, -id, s_1$ |
| $B$ | $\alpha_1, \alpha_2, \alpha_1 + \alpha_2, 2\alpha_1 + \alpha_2$ | $id, -id, s_1, -s_1$ |

The three following tables give the values of $q_w$ and $q'_w$ for each $\Sigma$:

| $P_\Sigma = P^*$ | | |
|---|---|---|
| $w$ | $q'_w$ | $q_w$ |
| id | 3 | 3 |
| -id | 0 | 0 |
| $s_1$ | 1 | 2 |

| $P_\Sigma = P^*$ | | |
|---|---|---|
| $w$ | $q'_w$ | $q_w$ |
| id | 2 | 3 |
| -id | 0 | 1 |

| $P_\Sigma = B$ | | |
|---|---|---|
| $w$ | $q'_w$ | $q_w$ |
| id | 3 | 4 |
| -id | 0 | 1 |
| $s_1$ | 2 | 3 |
| $-s_1$ | 1 | 2 |

Since the spectral sequence of the theorem 4.1 still degenerates in $E_2$, we can make analogous computation and get

**Theorem 5.6** *Let $Q = P^*$, there exists a filtration of the degree three ordinary boundary cohomology whose associated graded module is given by:*

- $Gr\ e_{P^*} H^3(\partial(\Gamma_0(p^r)\backslash \tilde{X}), L^a(\rho \otimes \chi; A)) = \bigoplus_{\substack{s\mathbf{P}_i \in \\ Cusp_{P^*}(\Gamma)}}$
  $eH^1_{cusp}(M_{\mathbf{P}_{i,s_1}}(\Gamma_0(p^r)), L_{P^*,s_1,2}(\rho \otimes \chi)) \oplus H^0(M_{\mathbf{P}_{i,id}}(\Gamma), L_{P^*,id,3}(\rho \otimes \chi))$
  $\bigoplus_{\substack{s\mathbf{P}_i \in \\ Cusp_P(\Gamma)}} eH^1_{cusp}(M_{\mathbf{P}_{i,id}}(\Gamma_0(p^r)), L_{P,id,2}(\rho \otimes \chi))$

- $Gr\ e_{P^*} H^3(\partial(\Gamma_1(p^r)\backslash \tilde{X}); L^a(\rho \otimes \chi; A)) = \bigoplus_{\substack{s\mathbf{P}_i \in \\ Cusp_{P^*}(\Gamma)}}$
  $Ind^{C_P(\mathbf{Z}/p^r\mathbf{Z})}_{i_{P^*}(T^{1,s_1}_{P^*}(\mathbf{Z}/p^r\mathbf{Z})\bar{E}_{\mathbf{P}_i}(\Gamma))} eH^1_{cusp}(M_{\mathbf{P}_{i,s_1}}(\Gamma_1(p^r)), L_{P^*,s_1,2}(\rho \otimes \chi))$
  $\oplus Ind^{C_{P^*}(\mathbf{Z}/p^r\mathbf{Z})}_{i_{P^*}(\bar{E}_{\mathbf{P}_i}(\Gamma))} H^0(M_{\mathbf{P}_{i,id}}(\Gamma), L_{P^*,id,3}(\rho \otimes \chi))$
  $\bigoplus_{\substack{s\mathbf{P}_i \in \\ Cusp_P(\Gamma)}} Ind^{C_{P^*}(\mathbf{Z}/p^r\mathbf{Z})}_{i_{P^*}(T^{1,id}_P(\mathbf{Z}/p^r\mathbf{Z})\bar{E}_{\mathbf{P}_i}(\Gamma))} eH^1_{cusp}(M_{\mathbf{P}_{i,id}}(\Gamma_1(p^r)), L_{P,id,2}(\rho \otimes \chi))$

**Remark:** The remark following Theorem 5.5 applies here as well.

## 5.4 The general case

When $F$ is different from $\mathbf{Q}$, we cannot complete the computation as simply because it does not seem possible to control the torsion arising in the term of the spectral sequences associated to the Levi decompositions of the parabolic defining the strata of the boundary. As it will become clear below, we can only make the computation modulo finite kernel or cokernel. It is the reason why we are led to introduce in definition 6.1 below an *ad hoc* $C_Q(\mathbf{Z}_p)$-module which will control the boundary cohomology (in degree $3d$); namely, its $\chi$-part will be isogenous to $H^{3d}(\partial(\Gamma_0(p^r)\backslash \tilde{X}), L^a(\rho \otimes \chi; A))$.

**Theorem 5.7** *For any $\chi$ regular, dominant with respect to $\rho$ and sufficiently separable (see Definition 4.6), the following canonical homomorphisms have finite kernel and cokernel:*

(i) $e_Q H^{3d}(\partial(\Gamma_0(p^r)\backslash \tilde{X}), L^a(\rho \otimes \chi; A)) \longrightarrow$

$$\bigoplus_{\substack{{}^s\mathbf{P}_i \in \\ Cusp_{P^*}(\Gamma)}} H^d_{!,ord}(M^1_{\mathbf{P}_i,\mathbf{s}_1}(\Gamma_0(p^r)), L_{P^*,\mathbf{s}_1,2d}(\rho \otimes \chi))$$

$$\oplus \bigoplus_{\substack{{}^s\mathbf{P}_i \in \\ Cusp_P(\Gamma)}} H^d_{!,ord}(M^1_{\mathbf{P}_i,\mathbf{s}_2}(\Gamma_0(p^r)), L_{P,\mathbf{s}_2,2d}(\rho \otimes \chi))$$

(ii) $e_Q H^{3d}(\partial(\Gamma_1(p^r)\backslash \tilde{X}), L^a(\rho \otimes \chi; A)) \longrightarrow$

$$\bigoplus_{\substack{{}^s\mathbf{P}_i \in \\ Cusp_{P^*}(\Gamma)}} Ind^{C_Q(\mathbf{Z}/p^r\mathbf{Z})}_{i_Q(T_Q^{1,\mathbf{s}_1}(\mathbf{Z}/p^r\mathbf{Z})\bar{E}_{\mathbf{P}_i}(\Gamma))} H^d_{!,ord}(M^1_{\mathbf{P}_i,\mathbf{s}_1}(\Gamma_1(p^r)), L_{P^*,\mathbf{s}_1,2d}(\rho \otimes \chi))$$

$$\oplus \bigoplus_{\substack{{}^s\mathbf{P}_i \in \\ Cusp_P(\Gamma)}} Ind^{C_Q(\mathbf{Z}/p^r\mathbf{Z})}_{i_Q(T_Q^{1,\mathbf{s}_2}(\mathbf{Z}/p^r\mathbf{Z})\bar{E}_{\mathbf{P}_i}(\Gamma))} H^d_{!,ord}(M^1_{\mathbf{P}_i,\mathbf{s}_2}(\Gamma_1(p^r)), L_{P,\mathbf{s}_2,2d}(\rho \otimes \chi))$$

*where $\mathbf{s}_i(v) = s_i$ for all $v|p$. Moreover, if we have chosen $p$ outside a finite set of primes depending on $\Gamma$ and $\rho$, this map is an isomorphism.*

**Remark:** The assumption of sufficient separability for $\rho$ can be removed when $Q$ is the Borel subgroup at each place $v$ above $p$. Indeed, instead of using lemma 4.11 (ii) which proves that $w(v)$ is independent of $v$, one uses the same argument as in the end of the proof of Corollary A1 to show that if $w$ contributes to the boundary cohomology, it satisfies $2d \leq l(w)$ or $l(w) \leq d$.

**Proof:** Let us detail only the proof of (i); the second point is similar. Recall that

$$e_Q G^{3d}_{\Sigma,w}(\Gamma_0(p^r); L^a(\rho \otimes \chi; A)) = \bigoplus_{\substack{{}^s\mathbf{P}_i \in \\ Cusp_\Sigma(\Gamma)}} e_Q H^{3d}(\mathbf{P}_{i,w} \cap \Gamma_0(p^r), L^a(\rho \otimes \chi; A))$$

Corollary 4.5 computes the cohomology groups of the right hand side up to finite kernel and cokernel. From lemma 4.11, we see that in this corollary, all the terms of the sum in the right hand side of index $l$ not divisible by $d$ are torsion and therefore can be neglected. Thus, for $P_\Sigma$ maximal we obtain a canonical map with finite kernel and cokernel:

$$e_Q G^{3d}_{\Sigma,w}(\Gamma_0(p^r); L^a(\rho \otimes \chi; A)) \to$$
$$\bigoplus_{\substack{{}^s\mathbf{P}_i \in \\ Cusp_\Sigma(\Gamma)}} \left[ e_{M^w_\Sigma \cap Q} H^d(M_{\mathbf{P}_{i,w}}(\Gamma_0(p^r)), L_{\Sigma,w,2d}(\rho \otimes \chi)) \right.$$
$$\oplus e_{M^w_\Sigma \cap Q} H^{2d}(M_{\mathbf{P}_{i,w}}(\Gamma_0(p^r)), L_{\Sigma,w,d}(\rho \otimes \chi))$$
$$\left. \oplus e_{M^w_\Sigma \cap Q} H^{3d}(M_{\mathbf{P}_{i,w}}(\Gamma_0(p^r)), L_{\Sigma,w,0}(\rho \otimes \chi)) \right]$$

Note that for the terms occurring in the right hand side, the action of the center of $M_{\mathbf{P}_{i,w}}$ on the coefficients is trivial. Recall that for any $\mathbf{Z}$-module $M$ with a trivial action of $\mathbf{Z}^n$ we have

$$H^i(\mathbf{Z}^n, M) = Hom(\bigwedge^i \mathbf{Z}^n, M)$$

Therefore for $k \in \{1,2,3\}$, the following map has finite kernel and cokernel

$$e_{M^w_\Sigma \cap Q} H^{kd}(M_{\mathbf{P}_{i,w}}(\Gamma_0(p^r)), L_{\Sigma,w,(3-k)d}(\rho \otimes \chi)) \to$$
$$\oplus_{l=0}^{d-1} Hom(\bigwedge^l E_{\mathbf{P}_i}, H^{ad-l}(M^1_{\mathbf{P}_{i,w}}(\Gamma_0(p^r)), L_{\Sigma,w,(3-k)d}(\rho \otimes \chi)))$$

Therefore it is torsion for $k = 3$. Moreover for $k = 2$ and $1$, it is non torsion only if $w = \mathbf{w}_{\mathbf{\Sigma},4-k}$ with $\mathbf{w}_{\mathbf{\Sigma},4-k}(v) = w_{\Sigma,4-k}$ (see tables of paragraph 5.3). This last point implies that $\pi_P^{-1}(\mathbf{w}_{\mathbf{P},4-k}) \cap \pi_{P^*}^{-1}(\mathbf{w}_{\mathbf{P}^*,4-k}) = \emptyset$ and therefore ker $r^{3d}$ is isogenous to the cohomology below where $E_{\mathbf{P}_i}$ is the group of units defined in section 5.2 :

$$\oplus_{\Sigma=P,P^*} \bigoplus_{\substack{{}^s\mathbf{P}_i \in \\ Cusp_\Sigma(\Gamma)}} \oplus_{l=0}^{d-1} Hom(\bigwedge^l E_{\mathbf{P}_i}, H^{d-l}_{!,ord}(M^1_{\mathbf{P}_{i,\mathbf{w}_{\mathbf{\Sigma},3}}}(\Gamma_0(p^r)), L_{\Sigma,\mathbf{w}_{\mathbf{\Sigma},3},2d}(\rho \otimes \chi)))$$
$$\oplus \oplus_{l=0}^{d-1} Hom(\bigwedge^l E_{\mathbf{P}_i}, H^{2d-l}_{!,ord}(M^1_{\mathbf{P}_{i,\mathbf{w}_{\mathbf{\Sigma},2}}}(\Gamma_0(p^r)), L_{\Sigma,\mathbf{w}_{\mathbf{\Sigma},2},d}(\rho \otimes \chi)))$$

We remark now that for any arithmetic subgroup $X \subset SL_2(F)$ and any regular weight $\lambda$, we have

$$H^i_!(X, V_\lambda(\mathbf{C})) = 0$$

if $i \neq d$. Therefore in the above sum, all the groups $H^{d-l}_{!,ord}$ vanish except for $l = 0$ while all the $H^{2d-l}_{!,ord}$ vanish. This takes care of kernel of $r^{3d}$. For the cokernel of $r^{3d-1}$, the same calculations show its finiteness. This concludes the proof of (i). The last assertion follows from two facts. First, the isogeny of corollary 4.5 is in fact an isomorphism outside a finite number of prime depending only of the weights modulo $p$. This verification is left to the reader. Let us note that in the case $\mathcal{Q} = B$, the isomorphism comes from Corollary 5.4. Second, all modules arising in the proof of Th.5.7 are cofree if $p$ is outside a finite set of primes depending on $\Gamma$ and $\rho$, for this result see [17]. ∎

**Definition 5.1** *We set $\mathcal{W}'_\partial(\Gamma_1(p^r); \rho \otimes \chi)$ the right-hand side of (ii) in Theorem 5.7 and*

$$\mathcal{W}'_\partial(\Gamma_1(p^\infty); \rho) = \varinjlim_r \mathcal{W}'_\partial(\Gamma_1(p^r); \rho)$$

*and we denote by $r'$ the canonical map from $e_Q H^{3d}(\partial(\Gamma_1(p^r)\backslash \tilde{X}), L^a(\rho \otimes \chi; A))$ to $\mathcal{W}'_\partial(\Gamma_1(p^\infty); \rho)$ induced by theorem 5.7.(ii). which is equivariant for the action of $C_Q(\mathbf{Z}/p^r\mathbf{Z})$.*

¿From Hida's Control theorem 5.1 for $SL_2$ over number fields of [16] and theorem 5.7.(i), we obtain easily:

**Lemma 5.2** *For $\chi$ as in Theorem 5.7, there is a canonical map with finite kernel and cokernel*

$$e_Q H^{3d}(\partial(\Gamma_0(p^r)\backslash \tilde{X}); L^a(\rho \otimes \chi; A)) \to \mathcal{W}'_\partial(\Gamma_1(p^\infty); \rho)[\omega_\chi]$$

**Theorem 5.8** *For any arithmetic character $\chi$ of $C_{M'}$ which is regular, dominant with respect to $\rho$ and sufficiently separable (see Definition 4.6), there is a canonical map with finite kernel or cokernel:*

$$H^{3d}_{!,Q-ord}(\Gamma_0(p^r), L^a(\rho \otimes \chi; A)) \to H^{3d}_{!,Q-ord}(\Gamma_1(p^\infty), L^a(\rho; A))[\omega_\chi]$$

*Moreover this map is an isomorphism if $p$ is chosen outside a finite number of primes depending on $\Gamma$ and $\rho$.*

**Proof:** We abbreviate the notations by:

$$\begin{aligned}
\mathcal{V} &= H^{3d}_{!,Q-ord}(\Gamma_1(p^\infty), L^a(\rho; A)) = H^{3d}_{!,Q-ord}(\Gamma_0(p), \mathcal{C}(\rho^1; A)) \\
\mathcal{W} &= H^{3d}_{Q-ord}(\Gamma_1(p^\infty), L^a(\rho; A)) = H^{3d}_{Q-ord}(\Gamma_0(p), \mathcal{C}(\rho^1; A)) \\
\mathcal{W}_\partial &= H^{3d}_{Q-ord}(\partial(\Gamma_1(p^\infty)\backslash \tilde{X}), L^a(\rho; A)) = H^{3d}_{Q-ord}(\partial(\Gamma_0(p)\backslash \tilde{X}), \mathcal{C}(\rho^1; A)) \\
\mathcal{W}'_\partial &= \mathcal{W}'_\partial(\Gamma_1(p^\infty; \rho))
\end{aligned}$$

By theorem 5.1 and lemma 5.1, we have the following commutative diagram:

$$\begin{array}{ccccccc}
0 & \to & H^{3d}_{!,Q-ord}(\Gamma_0(p^r), L^a(\rho \otimes \chi; A)) & \to & H^{3d}_{Q-ord}(\Gamma_0(p^r), L^a(\rho \otimes \chi; A)) & \to & H^{3d}_{Q-ord}(\partial(\Gamma_0(p^r)\backslash \tilde{X}), L^a(\rho \otimes \chi; A)) \\
 & & \| & & \| & & \| \\
0 & \to & H^{3d}_{!,Q-ord}(\Gamma_0(p^r), \mathcal{C}(\rho^1; A)[\omega_\chi]) & \to & H^{3d}_{Q-ord}(\Gamma_0(p^r), \mathcal{C}(\rho^1; A)[\omega_\chi]) & \to & H^{3d}_{Q-ord}(\partial(\Gamma_0(p^r)\backslash \tilde{X}), \mathcal{C}(\rho^1; A)[\omega_\chi]) \\
 & & \downarrow \iota_\mathcal{V} & & \downarrow \iota_\mathcal{W} & & \downarrow \\
0 & \to & \mathcal{V}[\omega_\chi] & \to & \mathcal{W}[\omega_\chi] & \to & \mathcal{W}_\partial[\omega_\chi]
\end{array}$$

Since $\iota_\mathcal{W}$ has finite kernel, it is also true for $\iota_\mathcal{V}$. In order to obtain our theorem, we simply need to prove that

$$\operatorname{corank}((\mathcal{V}[\omega_\chi]) \leq \operatorname{corank}(H^{3d}_{!,Q-ord}(\Gamma_0(p^r), L^a(\rho \otimes \chi; A))).$$

Let us consider the following commutative diagram with exact lines and $C_Q(\mathbf{Z}_p)$-equivariant maps:

$$\begin{array}{ccccccc} 0 & \to & \mathcal{V} & \to & \mathcal{W} & \to & \mathcal{W}_\partial \\ & & \downarrow & & \| & & \downarrow r' \\ 0 & \to & \mathcal{V}' & \to & \mathcal{W} & \to & \mathcal{W}'_\partial \end{array}$$

where $\mathcal{V}'$ is defined as the kernel of $\mathcal{W} \to \mathcal{W}'_\partial$. It implies that $\mathcal{V}[\omega_\chi] \hookrightarrow \mathcal{V}'[\omega_\chi]$; therefore by the lemma 5.2, we obtain the inequality we wanted. For the last assertion, it is enough to remark that $\mathcal{W}_\partial = \mathcal{W}'_\partial$ by Theorem 5.7) and that in this case, there is an exact control theorem (Theorem 7.1 of [17].■

# 6  Nearly ordinary cuspidal cohomology and the universal Hecke algebra

In this section, we deduce from the results of the previous sections the main theorems of the paper, namely, control and freeness of the nearly ordinary part of the cohomology of the Siegel varieties $S_r$, finiteness and torsion freeness of the big nearly ordinary cuspidal Hecke algebra over the Hida-Iwasawa algebra, existence of several variable families of cuspidal Hecke eigensystems interpolating a given cuspidal Hecke eigensystem. This is a developed version of our note [38] where we announced these results for $F = \mathbf{Q}$.

## 6.1  induction from $Sp_4$ to $GSp_4$

Let $U$ be a level subgroup in $H_f$ of level prime to $p$, say sufficiently deep so that the discrete subgroup of $H(\mathbf{Q})$ associated to $U$ has no torsion and $Q$ be a standard parabolic of $H \otimes \mathbf{Z}_p$. Let $\rho : M \to GL_\mathcal{O}(V)$ a group-scheme morphism as before; for each $\chi \in X^*(C_M)$ dominant with respect to $\rho$, we consider the local system $\mathcal{L}^a(\rho \otimes \chi; A)$ over $S_r(U)$, defined by the action of $U_1(p)$ on

$$\mathbf{L}^a(\rho \otimes \chi; A) = \mathbf{L}^a(\rho \otimes \chi; \mathcal{O}) \otimes A$$

where

$$\mathbf{L}^a(\rho \otimes \chi; \mathcal{O}) = \{f : I/Q^+(\mathbf{Z}_p) \to V; f \text{ is polynomial}$$
$$\text{and } f(xm) = \rho \otimes \chi(m^{-1})f(x) \text{ for } m \in M(\mathbf{Z}_p)\}$$

Note that for any $r \geq 1$ we have:

$$H(\mathbf{A}) = \coprod_{t \in R} H(\mathbf{Q}) t U_0(p^r) H(\mathbf{R})^+$$

where $H(\mathbf{R})^+$ is the neutral component of $H(\mathbf{R})$; and where the set $R$ is a finite set of elements of $H_f$ whose components at places in the level of $U_0(p)$ are equal to 1 and such that its image $\nu(R)$ by the multiplicator $\nu: GSp_4(F_\mathbf{A}) \to F_\mathbf{A}^\times$ is a complete set of representatives of the $U$-ray-class group of field $F$: $Cl_U^+ = F^\times \backslash F_\mathbf{A}^\times / \nu(U) \times F_\infty^{\times+}$. Note that $R$ is independent of $r$ since $\nu(U_0(p))$ is unramified at $p$.

Moreover, recall that

$$U_0(p^r)/U_1(p^r) \cong C_M(\mathbf{Z}/p^r\mathbf{Z})$$

hence if $S_r$ is a complete system of representatives $s \in H(\mathbf{Z}_p)$ of $C_M(\mathbf{Z}/p^r\mathbf{Z}))$, one has

$$H(\mathbf{A}) = \coprod_{t \in R} \coprod_{s \in S_r} H(\mathbf{Q}) t s U_1(p^r) H(\mathbf{R})^+$$

Note that $s$ normalizes $U_1(p^r)$; let us put for any $t \in R$ and for $i = 0$ or 1, we set:

$$\Gamma_{i,t}(p^r) = H(\mathbf{Q}) \cap t U_i(p^r) H(\mathbf{R}) t^{-1} \text{ and } \Gamma'_{i,t}(p^r) = \Gamma_{i,t}(p^r) \cap H'(\mathbf{Q})$$

Let for $i = 0$ or 1:

$$S(\Gamma'_{i,t}(p^r)) = \Gamma'_{i,t}(p^r) \backslash H'_\infty / (U_\infty \cap H'_\infty)$$

note that for any $t \in R$,

$$\Gamma_{0,t}(p^r)/\Gamma_{1,t}(p^r) \cong C_{M'}(\mathbf{Z}/p^r\mathbf{Z}).\bar{E}_t$$

where $\bar{E}_t \subset Z_H(\mathbf{r}/p^r\mathbf{r})$ is the reduction modulo $p^r$ of the image in $Z_H(\mathbf{r})$ of $\nu(\Gamma_{0,t})$. So, since $S_r(U) = H_\mathbf{Q} \backslash H_\mathbf{A} / U_1(p^r) U_\infty$, we see that

$$S_r(U) \cong \coprod_{t \in R} \mathrm{ind}_{C_{M'}(\mathbf{Z}/p^r\mathbf{Z})\bar{E}_t}^{C_M(\mathbf{Z}/p^r\mathbf{Z})} S(\Gamma'_{1,t}(p^r))$$

hence

$$eH^q(S_r(U), \mathcal{L}^a(\rho \otimes \chi; A)) = \bigoplus_{t \in R} \mathrm{ind}_{C_{M'}(\mathbf{Z}/p^r\mathbf{Z})\bar{E}_t}^{C_M(\mathbf{Z}/p^r\mathbf{Z})} eH^q(\Gamma'_{1,r}(p^r), L^a(\rho \otimes \chi; A))$$

Note that $\nu$ induces an isomorphism

$$C_M(\mathbf{Z}/p^r\mathbf{Z})/C_{M'}(\mathbf{Z}/p^r\mathbf{Z})\bar{E}_t \cong (\mathbf{r}/p^r\mathbf{r})^\times / \bar{E}_t$$

## 6.2 Hecke operators and Hida-Iwasawa algebra

We let act on $(S_r(U), \mathcal{L}^a(\rho \otimes \chi; A))$ the so-called Hecke correspondences defined as follows. Let $S$ be the set of places of $F$ occurring in $level(U)$. Let $U' = U_1(p^r)$ and $S(U') = S_r(U)$.

- for any $h \in H_f^S \cap M_4(\hat{\mathbf{Z}})$, that is, such that $h_v = 1$ for $v \in S$ or dividing $p$, we define $[U'hU']$: let $U'' = U' \cap hU'h^{-1}$ and consider the map $[h] : (S(U'') \to (S(h^{-1}U''h)$ induced by right multiplication by $h$ on $H(\mathbf{A})$ and by pull-back by $h$ on the sheaf (without any action on the group $L^a(\rho \otimes \chi; A)$ itself). The diagram

$$\begin{array}{ccc} S(U'') & \stackrel{[h]}{\to} & S(U'') \\ \downarrow & & \downarrow \\ S(U') & & S(U') \end{array} \quad (6.1)$$

induces on cohomology the desired action of $[U'hU']$. Recall the classical Notation: for $v \notin S$ (and prime to $p$):$T_v = [U'\mu(1_2; \varpi_v)U']$, $S_v = [U'diag(\varpi_v, \varpi_v, \varpi_v, \varpi_v)U']$ and $R_v = [U'\mu^*(1; \varpi_v.1_2)U'] - (Nv^2 - 1)S_v$

- for $h_p \in D_p$, one considers similarly $[U'h_pU']$ defined by the diagram 6.1 but for $[h_p]$ acting on the $L^a(\rho \otimes \chi; A)$-bundle of $S(U'')$ by $(x, \ell) \mapsto (xh, h_p^{-1}\ell)$. Here, $D_p = \prod_{v|p} D_v$ denotes the subsemigroup of $Z_M(\mathbf{Q}_p)$ defined in Section 2.2.

- for $h \in Z_M(\mathbf{Z}_p)$, we let $[U'hU']$ act by normal action since $Z_M(\mathbf{Z}_p)$ viewed in $U_0(p^r)$ normalizes $U_1(p^r)$.

- for $z \in Z_H(\mathbf{A})$, we let $U'zU'$ act by central action.

**Remarks:**
1) The normal action factors through $Z_M(\mathbf{Z}/p^r\mathbf{Z})$.
2) Let
$$\Omega_r = \{z \in Z_H(\hat{\mathbf{r}}); z \equiv 1 \mod p^r\}$$

and put
$$\mathcal{R}_r = Z_H(\mathbf{Q})\backslash Z_H(\mathbf{A})/\Omega_r Z_H(\mathbf{R})^+$$

then the central action factors through $\mathcal{R}_r$.

Let
$$\mathcal{W}^q(\rho) = \varinjlim_r eH^q(S_r(U), \mathcal{L}^a(\rho; A)) \quad \mathcal{W}_\partial^q(\rho) = \varinjlim_r eH_\partial^q(S_r(U), \mathcal{L}^a(\rho; A))$$

where $H_\partial^q(S_r(U), \mathcal{L}^a(\rho; A))$ is the cohomology of the boundary of the Borel-Serre compactification of $S_r(U)$. We define $\mathcal{V}^q(\rho)$ by the exact sequence

$$0 \to \mathcal{V}^q(\rho) \to \mathcal{W}^q(\rho) \to \mathcal{W}_\partial^q(\rho)$$

We also introduce the interior cohomology groups $H_!^q$:

$$eH_!^q(S_r(U), \mathcal{L}^a(\rho; A)) = \operatorname{Ker}(eH^q(S_r(U), \mathcal{L}^a(\rho; A)) \to eH_\partial^q(S_r(U), \mathcal{L}^a(\rho; A)))$$

In particular,

$$\mathcal{V}^q(\rho) = \varinjlim_r eH_!^q(S_r(U), \mathcal{L}^a(\rho; A))$$

**Definition 6.1** *We define the nearly ordinary Hecke algebra $\mathbf{h}_{r,\rho}^{q,S}(\psi \otimes \varepsilon)$ of level $U_1(p^r)$ (outside $S$) by as the $\mathcal{O}$-subalgebra generated by the operators*

- $[U_1(p^r) h U_1(p^r)]$ $(h \in H_f^{S,p})$,

- $[U_1(p^r) h_p U_1(p^r)]$ $(h_p \in D_p)$,

- $[U_1(p^r) h U_1(p^r)]$ $(h \in Z_M(\mathbf{Z}/p^r\mathbf{Z}))$

- *and* $[U_1(p^r) z U_1(p^r)]$ $(z \in \mathcal{R}_r)$

*acting on $eH_!^q(S_r(U), \mathcal{L}^a(\rho; A))$.*

**Remark:**

One can also define similar Hecke operators of level $U$ (replacing $U_1(p^r)$ by $U$). Then, an important property of the isomorphism in Proposition 3.2 is that it commutes with the Hecke operators for $h \in H_f^{S,p}$; moreover for $h_p \in D_p$, the corresponding operators are congruent modulo $p$ (on the prime-to-$p$ level side, one has to divide $[U h_p U]$ by the constant $\omega(h_p)$).

When $r \geq 1$ grows, the endomorphisms $[U_1(p^r) h U_1(p^r)]$ $(h \in H_f^S)$, $[U_1(p^r) h_p U_1(p^r)]$ $(h_p \in D_p)$, $[U_1(p^r) h U_1(p^r)]$ $(h \in Z_M(\mathbf{Z}/p^r\mathbf{Z}))$ and $[U_1(p^r) z U_1(p^r)]$ $(z \in \mathcal{R}_r)$ form a compatible projective system; one can therefore consider the algebra

$$\mathbf{h}_{\rho \otimes \chi}^{q,S} = \varprojlim_r \mathbf{h}_{r,\rho}^{q,S}(\chi)$$

By definition, it acts faithfully on $\mathcal{V}^q(\rho)$.

**Definition 6.2** *For any $r \geq 1$, let $\mathbf{H}_r$ be the amalgamated sum:*

$$\mathbf{H}_r = Z_H(\mathbf{Q}) \backslash Z_H(\mathbf{A}) / \Omega_r Z_H(\mathbf{R}) \oplus_{Z_H(\mathbf{Z}/p^r\mathbf{Z})} Z_M(\mathbf{Z}/p^r\mathbf{Z})$$

*where the amalgamation is taken for $z_p \mapsto (\alpha(z_p^{-1}), \beta(z_p))$ where $\alpha$ is induced by the inclusion $\mathbf{Z}_p \subset \mathbf{A}$ and $\beta$ is given by the homomorphism $Z_H \subset Z_M$.*

We put

$$\mathbf{H} = \varprojlim_r \mathbf{H}_r$$

*We call it the Hida group of $H$.*

*The Hida-Iwasawa algebra $\Lambda$ is defined as the projective limit of the group algebras*

$$\Lambda_r = \mathcal{O}[\mathbf{H}_r]$$

*for the natural transition homomorphisms:*

$$\Lambda = \mathcal{O}[[\mathbf{H}]]$$

**Remarks:**

1) Let $\mathbf{H}(p)$ be the largest torsion-free pro-$p$ subgroup of $\mathbf{H}$ and $\Lambda(p) = \mathcal{O}[[\mathbf{H}(p)]]$; then $\mathbf{H} = \Phi \times \mathbf{H}(p)$ where $\Phi$ is a finite group and $\Lambda = \Lambda(p)[\Phi]$ is the group algebra of $\Phi$ over $\Lambda(p)$. Let

$$\mathcal{R}_r = Z_H(\mathbf{Q}) \backslash Z_H(\mathbf{A}) / \Omega_r Z_H(\mathbf{R}) \times Z_{M'}(\mathbf{Z}/p^r\mathbf{Z}) \text{ and } \mathcal{R} = \varprojlim_r \mathcal{R}_r$$

and $\mathcal{R}(p)$ its largest torsion-free pro-$p$ subgroup. Note that since $p$ is odd, there is a canonical isomorphism:

$$\mathbf{H}(p) \cong \mathcal{R}(p) \times Z_{M'}(\mathbf{Z}_p)(p)$$

given by $(z,t) \mapsto (z,t')$ where $t = z_p t'$, $z_p \in Z_H(\mathbf{Z}_p)(p)$ and $t' \in Z_{M'}(\mathbf{Z}_p)(p)$.

2) The relative Krull dimension of $\Lambda$ over $\mathcal{O}$ is equal to the $\mathbf{Z}_p$-rank $r_p(\mathbf{H})$ of $\mathbf{H}$; one has $r_p(\mathbf{H}) = 1 + \delta + \sum_{v|p} r_v d_v$ where $r_v$ is the rank of $C_{M'_v}$ (that is, 1 if $\Pi_v$ is maximal and 2 if $\Pi_v = B$) and $\delta$ is the defect to the Leopoldt conjecture for $(F,p)$. For instance, if $Q$ is the Borel subgroup, one has $r_p(\mathbf{H}) = 2d + 1 + \delta$. For the general definition of the Hida group, see [36]. The number $r_p(\mathbf{H})$ is important for us since it will be the number of $p$-adic parameters for the space of nearly ordinary deformations of a given nearly ordinary Hecke eigensystem.

For any $r \geq 1$, the group $\mathbf{H}_r$ acts naturally on $eH^q(S_r(U), \mathcal{L}^a(\rho; A))$ and $eH^q_!(S_r(U), \mathcal{L}^a(\rho; A))$ by

$$<>_{r,\rho\otimes\chi}: (z,t) \mapsto [U_1(p^r)ztU_1(p^r)]$$

These actions are compatible when $r$ grows; hence $\mathbf{H}$ acts on $\mathcal{W}^q_{\rho\otimes\chi}$ and $\mathcal{V}^q_{\rho\otimes\chi}$ These actions can be extended uniquely by linearity into $\mathcal{O}$-algebra homomorphisms

$$\Lambda_r \to \text{End}_\mathcal{O}(\mathcal{W}^q_{r,\rho\otimes\chi})$$

$$\Lambda_r \to \text{End}_\mathcal{O}(\mathcal{V}^q_{r,\rho\otimes\chi})$$

resp.

$$\Lambda \to \text{End}_\mathcal{O}(\mathcal{W}^q_{\rho\otimes\chi})$$

$$\Lambda_r \to \text{End}_\mathcal{O}(\mathcal{V}^q_{\rho\otimes\chi})$$

**Definition 6.3** *1) An arithmetic character $\theta$ of $\mathbf{H}$ is a continuous homomorphism $\theta : \mathbf{H} \to \mathcal{O}^\times$ such that its restriction to some $p$-adic open subgroup coincides with an algebraic character of the algebraic group $Z_H \times Z_M$.*

*2) It is called arithmetic dominant with respect to $\rho$ if one can write for $h = (z,t) \in \mathbf{H}$, $\theta(h) = \phi(z)\omega_{\varepsilon\chi}(t)$ where $\chi \in X^*(C_M)$ is dominant with respect to $\rho$, $\varepsilon : C_M(\mathbf{Z}_p) \to \mathcal{O}^\times$ has finite order, $\phi$ is locally algebraic on $\mathcal{R}$ (and $\phi = \omega_{\varepsilon\chi}$ on $(\alpha(z_p^{-1}), \beta(z_p))$ for any $z_p \in Z_H(\mathbf{Z}_p)$).*

**Definition 6.4** *For any arithmetic character $\theta$ of $\mathbf{H}$ dominant with respect to $\rho$, we define the so-called arithmetic prime $\wp_\theta$ of $\mathrm{Spec}\,\Lambda$ associated to $\theta$ as $\wp_\theta : \mathrm{Ker}(\Lambda \to \mathcal{O})$ where the $\mathcal{O}$-algebra homomorphism is defined by linearity and continuity from the group homomorphism: $h \mapsto \theta(h)$). We say that $\wp_\theta$ is algebraic if $\theta$ itself is algebraic.*

## 6.3 Control Theorems

As in Section 3.5, we fix a Dedekind ring $\mathcal{O}_0$ in $\overline{\mathbf{Q}}$ finite over $\mathbf{Z}$, and for a prime $p$, we fix a $p$-adic embedding $\iota_p$ of $\overline{\mathbf{Q}}$, hence of $\mathcal{O}_0$, and we denote by $\mathcal{O}$ the completion of $\iota_p(\mathcal{O}_0)$.

We obtain easily the following theorems of independence of the weight and of exact control for the cohomology:

**Theorem 6.1** *For any $\rho$ with $H$-admissible highest weight,*
*there is a canonical Hecke-equivariant isomorphism:*

$$j_{\chi\varepsilon} : \mathcal{V}_\rho^q \cong \mathcal{V}_{\rho\otimes\chi\varepsilon}^q$$

*such that for $\zeta \in C_M(\mathbf{Z}_p)$, we have*

$$j_{\chi\varepsilon}(<\zeta>c) = \chi\varepsilon(\zeta^{-1}) <\zeta>_{\chi\varepsilon} j_{\chi\varepsilon}(c)$$

Let $\theta = \phi \otimes \omega_{\varepsilon\chi}$ be an arithmetic character of $\mathbf{H}$ where $\chi \in X^*(C_M)$ is dominant with respect to $\rho$ and $\varepsilon : C_M(\mathbf{Z}/p^r\mathbf{Z}) \to \mathcal{O}^\times$ Let $\mathcal{V}_{r,\rho\otimes\varepsilon\psi}(\phi)$ resp. $\mathcal{V}_\rho^{3d}(\phi)$ be the largest submodule of $\mathcal{V}_{r,\rho\otimes\varepsilon\psi}$ resp. $\mathcal{V}_\rho^{3d}$ on which $\Phi$ acts via $\phi$. Then, note that the natural map

$$\mathcal{V}_{r,\rho\otimes\varepsilon\psi}^{3d} \to \mathcal{V}_\rho^{3d}[\omega_{\varepsilon\psi}]$$

induces

$$\mathcal{V}_{r,\rho\otimes\varepsilon\psi}^{3d}(\phi) \to \mathcal{V}_\rho^{3d}[\wp_\theta]$$

**Theorem 6.2** *Let $\rho$ be any absolutely irreducible representation of $M_Q$;*

1. **Weak Control:** Let $p \neq 2$ be an arbitrary rational prime. For any arithmetic character $\theta$ of $\mathbf{H}$ regular dominant with respect to $\rho$ and sufficiently separable (see Definition 4.6), say $\theta = \phi \otimes \omega_{\varepsilon\psi}$ for $\psi \in X^*(C_M)$ dominant with respect to $\rho$ and $\varepsilon : C_M(\mathbf{Z}/p^r\mathbf{Z}) \to \mathcal{O}^\times$, the maps
$$\begin{array}{rcl} \mathcal{V}^{3d}_{r,\rho\otimes\varepsilon\psi}(\phi) & \to & \mathcal{V}^{3d}_\rho[\wp_\theta] \\ \mathcal{W}^{3d}_{r,\rho\otimes\varepsilon\psi}(\phi) & \to & \mathcal{W}^{3d}_\rho[\wp_\theta] \end{array}$$
have finite kernel and cokernel.

2. **Strong Control:** Assume moreover that $\rho = \rho_0 \otimes_{\mathcal{O}_0} \mathcal{O}$ and that its highest weight is regular, there exist a finite set of primes $S_{U,\rho}$ such that for $p \notin S_{U,\rho}$,

(i) for any $q \in [0, 3d[$,
$$\mathcal{V}^q_\rho = \mathcal{W}^q_\rho = 0$$

(ii) for any arithmetic character $\theta = \phi \otimes \omega_{\varepsilon\psi}$ of $\mathbf{H}$ as above, with the supplementary condition that $\varepsilon\psi$ is congruent to 1 modulo $\pi$, there are canonical isomorphisms

$$\begin{array}{rcl} \mathcal{V}^{3d}_{r,\rho\otimes\varepsilon\psi} & \cong & \mathcal{V}^{3d}_\rho[\omega_{\varepsilon\psi}] \\ \mathcal{W}^{3d}_{r,\rho\otimes\varepsilon\psi} & \cong & \mathcal{W}^{3d}_\rho[\omega_{\varepsilon\psi}] \end{array}$$

**Proof:** This results from the corresponding theorem (cf. thm 5.8 for $Sp_4$) and from the induction formula

$$\mathcal{V}^{3d}_{r,\rho\varepsilon\psi} = \bigoplus_{t \in R} \mathrm{ind}^{C_M(\mathbf{Z}/p^r\mathbf{Z})}_{C_{M'}(\mathbf{Z}/p^r\mathbf{Z})E_t(p^r)} eH^{3d}_!(\Gamma'_{1,t}(p^r), L^a(\rho \otimes \chi; A))$$

∎

**Theorem 6.3** Let $F = \mathbf{Q}$. Let $\rho$ be defined over $\mathcal{O}_0$; assume either

(i) for all $t \in R$, $H^i_{Q-no}(\Gamma_{0,t}(p), L^a(\rho; K/\mathcal{O})) = 0$ for $i = 1, 2$

or

(ii) $\rho$ is regular and $p$ does not divide $\prod_{i=1}^{i=3} |H^i(\Gamma_t, L_0(\rho; \mathcal{O}_0)_{tors}|$,

then, for any arithmetic character $\theta$ of $\mathbf{H}$ dominant with respect to $\rho$, say $\theta = \phi \otimes \omega_{\varepsilon\psi}$ for $\psi \in X^*(C_M)$ dominant with respect to $\rho$, congruent to 1 mod. $\pi$ and $\varepsilon : C_M(\mathbf{Z}/p^r\mathbf{Z}) \to \mathcal{O}^\times$ a character congruent to 1 modulo $\pi$, there is a canonical isomorphism

$$\mathcal{V}^3_{r,\rho\otimes\varepsilon\psi}(\phi) \cong \mathcal{V}^3_\rho[\wp_\theta]$$

**Remark:** If $Q$ is the Klingen parabolic or the Borel subgroup, it follows from K. Buecker's thesis [6] that the part of assumption *(i)* relative to $i = 1$ is satisfied. Note also that for *(i)* to be fulfilled, it is sufficient that

$$H^i_{Q-no}(\Gamma'_{0,t}(p), L^a(\rho; \mathcal{O}/\pi\mathcal{O})) = 0$$

for $i = 1, 2$.

**Proof :** The proof is the same as for the previous Theorem 6.2, except that one uses here Th 3.3 which does not require any assumption on $p$.

Let

$$\mathbf{V}_{r,\rho\varepsilon\psi} = \mathrm{Hom}_{\mathcal{O}}(\mathcal{V}^{3d}_{r,\rho\varepsilon\psi}, K/\mathcal{O})$$

and

$$\mathbf{V}_\rho = \mathrm{Hom}_{\mathcal{O}}(\mathcal{V}^{3d}_\rho, K/\mathcal{O})$$

**Corollary 6.1** *Let $\rho$ and $\theta$ be as in Th.6.2 (Weak or Strong). The natural homomorphism*

$$\mathbf{V}_\rho/\wp_\theta \mathbf{V}_\rho \to \mathbf{V}_{r,\rho\varepsilon\psi}(\phi)$$

*has finite kernel and cokernel in weak control case and is an isomorphism in the strong control case.*

**Corollary 6.2** *Under the same assumptions as in Th.6.2, $\mathbf{V}_\rho$ is of finite type over $\Lambda$.*

**Proof:** It is a simple application of the topological Nakayama's lemma: let $\mathfrak{m}$ be the maximal ideal of $\Lambda$. We see from corollary 6.1 that $\mathbf{V}_\rho/\mathfrak{m}\mathbf{V}_\rho$ is a finite group. ∎

**Corollary 6.3** *Let $\rho$ be as in Th.6.2 (Weak or Strong). For any arithmetic character $\theta$ of $\mathbf{H}$ dominant with respect to $\rho$ and sufficiently separable, the natural homomorphism*

$$\mathbf{h}^{3d,S}_\rho \otimes \Lambda_{\wp_\theta}/\wp_\theta \to \mathbf{h}^{3d,S}_{r,\rho\varepsilon\psi} \otimes K$$

*is surjective with kernel contained in the radical.*

**Proof :** We have to prove that a non trivial idempotent of the left-hand side does not map to zero. Let $\bar{e}$ be such an idempotent. By Hensel lemma, it lifts to an idempotent $e$ of $(\mathbf{h}^{3d,S}_\rho)_{\wp_\theta}$. If $\bar{e}$ maps to zero, by Corollary 6.1, we have

$$e.(\mathbf{V}_\rho)_{\wp_\theta} \subset \wp_\theta.(\mathbf{V}_\rho)_{\wp_\theta}$$

Since $e$ is an idempotent, this implies that $e = 0$ and $\bar{e} = 0$. ∎

## 6.4 Structure over the Hida-Iwasawa algebra

We are going to prove that $\mathcal{V}_\rho^{3d}$ is co-free over the Hida-Iwasawa algebra provided that $p$ is outside an explicit finite set of primes. In order to prove that result we need the following general duality theorem. For any dominant character $\chi \in X^*(T)$ corresponding to a triple $(a, b; c)$, we put $\chi^\vee = (a, b; -c)$.

**Theorem 6.4** *The following pairing yields a perfect Pontrjagin duality*

$$(-,-)_\chi \ : \ H^q_{c,ord}(\Gamma'_{0,t}\backslash X; L^a(\rho \otimes \chi, \mathcal{O})) \otimes H^{6d-q}_{ord}(\Gamma'_{0,t}(p^r)\backslash X; L^a(\rho^\vee \otimes \chi^\vee; A)) \to A$$

*with $A = K/\mathcal{O}$ and $(x,y)_\chi = \varphi_\chi(x \cup W_{\Gamma'_{0,t}(p^r)}(y))$ where $\varphi_\chi$ is induced by the natural pairing between $L^a(\rho \otimes \chi; \mathcal{O})$ and $L^a(\rho^\vee \otimes \chi^\vee; \mathcal{O})$ and with $W_{\Gamma'_{0,t}(p^r)} \in GSp_4(F)$ normalizing $\Gamma$ such that $W_{\Gamma'_{0,t}(p^r)} \equiv \begin{pmatrix} 0_2 & 1_2 \\ p^r 1_2 & 0_2 \end{pmatrix} \mod U_1(p^{2r})$*

**Proof:** Recall $A_r = p^{-r}\mathcal{O}/\mathcal{O}$. Consider the following commutative diagram where $r \leq s$:

$$\begin{array}{ccccc}
H^q_{c,ord}(\Gamma'_{0,t}(p^r)\backslash X, L^a(\rho \otimes \chi; A_s)) & \otimes & H^{6d-q}_{ord}(\Gamma'_{0,t}(p^r)\backslash X, L^a(\rho^\vee \otimes \chi^\vee; A_s)) & \to & H^{6d}_c(\Gamma'_{0,t}(p^r)\backslash X; A_s) \cong A_s \\
\| & & \| & & \| \\
H^q_{c,ord}(\Gamma'_{0,t}(p^s)\backslash X, L^a(\rho \otimes \chi; A_s)) & \otimes & H^{6d-q}_{ord}(\Gamma'_{0,t}(p^s)\backslash X, L^a(\rho^\vee \otimes \chi^\vee; A_s)) & \to & H^{6d}_c(\Gamma'_{0,t}(p^s)\backslash X, A_s) \cong A_s \\
\| & & \| & & \| \\
H^q_{c,ord}(\Gamma'_{0,t}(p^s)\backslash X, V(\rho \otimes \chi) \otimes A_s) & \otimes & H^{6d-q}_{ord}(\Gamma'_{0,t}(p^s)\backslash X, V(\rho^\vee \otimes \chi^\vee) \otimes A_s) & \to & H^{6d}_c(\Gamma'_{0,t}(p^s)\backslash X, A_s) \cong A_s
\end{array}$$

The isomorphisms of the left hand side for compact support cohomology follow from versions for compact support cohomology of

- Hida's lemma for lowering the $p$-level (Lemma 3.1 above with the same proof).
- The contraction lemma (Proposition 3.1 above with the same proof).

The last line is exact by the Poincaré duality theorem and the fact that $(e_Q)^* = W_{\Gamma'_{0,t}(p^r)} e_Q (W_{\Gamma'_{0,t}(p^r)})^{-1}$. Now we get our result, taking inductive limite in the equality :

$$Hom(H^q_{c,ord}(\Gamma'_{0,t}(p^r)\backslash X, L^a(\rho \otimes \chi; A_s)), A) \cong H^{6d-q}_{ord}(\Gamma'_{0,t}(p^r)\backslash X, L^a(\rho^\vee \otimes \chi^\vee; A_s)).$$

∎

**Corollary 6.4** *If $p \notin S_{U,\rho} \cup S_{U,\rho^\vee}$ then $\mathcal{V}_\rho^{3d}$ is cofree of cofinite type over the Hida-Iwasawa algebra.*

**Proof:** In order to get our result, we prove that $H^{3d}_{c,ord}(\Gamma'_{0,t}(p^r)\backslash X, L^a(\rho\otimes\chi; A))$ is p-divisible for a densely populated set of sufficiently separable, regular dominant algebraic characters $\chi$. By the following exact sequence:

$$H^{3d}_{c,ord}(\Gamma'_{0,t}(p^r)\backslash X, L^a(\rho \otimes \chi; K)) \to H^{3d}_{c,ord}(\Gamma'_{0,t}(p^r)\backslash X, L^a(\rho \otimes \chi; A)) \to$$
$$\to H^{3d+1}_{c,ord}(\Gamma'_{0,t}(p^r)\backslash X, L^a(\rho \otimes \chi; \mathcal{O})) \to \ldots$$

this will be true if we know that $H^{3d+1}_{c,ord}(\Gamma'_{0,t}(p^r)\backslash X, L^a(\rho\otimes\chi;\mathcal{O}))$ is without $p$-torsion. But by the previous theorem, the Pontrjagin dual of this group is $H^{3d-1}_{ord}(\Gamma'_{0,t}(p^r)\backslash X, L^a(\rho^\vee\otimes\chi^\vee;A))$ which is zero by theorem 6.2.(i) for $p \notin S_{U,\rho^\vee}$.∎

**Corollary 6.5** *If $p \notin S_{U,\rho} \cup S_{U,\rho^\vee}$ then $\mathbf{h}^{3d,S}_\rho$ is torsion-free over the Hida-Iwasawa algebra.*

### 6.5 Families of eigensystems

**Definition 6.5** *A family of Q-nearly ordinary Siegel-Hilbert cusp eigensystems of genus 2 is the datum of an homomorphism of $\Lambda$-algebra from the universal Hecke algebra $\mathbf{h}^{3d,S}_\rho$ into a finite and flat extension $\mathbf{J}$ of $\Lambda(p)$.*
$$\lambda : \mathbf{h}^{3d,S}_\rho \to \mathbf{J}$$

The above definition is justified by the following:

**Corollary 6.6** *Let $\lambda$ be such a family. For any arithmetic character $\theta$ of $\mathbf{H}$ separable and dominant with respect to $\rho$, say $\theta = \phi \otimes \omega_{\varepsilon\psi}$ for $\psi \in X^*(C_M)$ dominant with respect to $\rho$ and $\varepsilon : C_M(\mathbf{Z}/p^r\mathbf{Z}) \to \mathcal{O}^\times$ such that $\varepsilon\psi$ is congruent to 1 modulo $\pi$. Then for any $\wp'$ prime ideal above $\wp_\theta$, there is the following commutative diagram:*

$$\begin{array}{ccc}
\mathbf{h}^{3d,S}_\rho \otimes \Lambda(p)/\wp_\theta & \stackrel{\lambda \bmod \wp'}{\to} & \mathbf{J}/\wp' \\
\downarrow & & \| \\
\mathbf{h}^{3d,S}_{r,\rho}(\phi \otimes \varepsilon\psi) & \stackrel{\lambda(\epsilon\psi)}{\to} & \mathbf{J}/\wp'
\end{array}$$

**Proof :** This results immediately from 6.3.∎

**Corollary 6.7** *Let $\pi$ be a cuspidal representation of $GSp_{4/F}$ whose archimedian component belongs to the discrete serie with cohomological regular weight $\theta_0$. Then there exists a finite number of prime $S(\pi) \supset Ram(\pi)$ such that if $p \notin S(\pi)$ and $\pi$ is Q-ordinary at $p$ there exist a family of Q-nearly ordinary Siegel-Hilbert cusp eigensystems of genus 2 whose specialisation mod. $\wp_{\theta_0}$ "is" $\lambda_\pi$ the character of the Hecke algebra corresponding to $\pi$.*

**Proof :** This results from the previous corollary, the torsion-freeness of the universal Hecke algebra and the Going-Down theorem for lifting ideal in normal extensions.∎

## 7 Application to Galois Representations

We recall below some classical conjectures on Galois representations associated to cohomological automorphic representations. Let us start by some preliminaries. Recall we have fixed in the previous sections embeddings $\iota_p$ and $\iota_\infty$ of $\overline{\mathbf{Q}}$ in $\mathbf{C}_p$ and $\mathbf{C}$. Let $\pi$ be a cohomological cuspidal representation of $GSp_{4/F}$. Then it must occur in

$$H^{3d}(S(U); \bigotimes_{\sigma \in I_F} L_{(a_\sigma, b_\sigma; c_\sigma)}(\mathbf{C}))$$

where $L_{(a_\sigma, b_\sigma; c_\sigma)}(\mathbf{C})$ is the irreducible representation of $GSp_{4/\mathbf{C}}$ of highest weight $(a_\sigma, b_\sigma; c_\sigma)$ over which $GSp_4(F)$ acts via $\iota_\infty \circ \sigma$. Note that $c_\sigma$ is independent of $\sigma$ (because it is the infinite type of the central character of $\pi$); we denote sometimes this common value by $c$. We denote by $\lambda_\pi$ the character of the Hecke algebra corresponding to $\pi$ and by $E_\pi$ the subfield of $\overline{\mathbf{Q}}$ generated by the values of $\lambda_\pi$; we embed it canonically in $\overline{\mathbf{Q}}_p$ via $\iota_p$.

**Conjecture 2** *Let $\pi$ be a cohomological cuspidal representation of $GSp_{4/F}$. There exists a continuous Galois representation $\varrho_\pi$ unramified outside $Ram(\pi) \cup S_p$:*

$$\varrho_\pi : Gal(\bar{F}/F) \to GSp_4(\bar{\mathbf{Q}}_p)$$

*such that for all prime $w \notin Ram(\pi) \cup S_p$ the characteristic polynomial of $\varrho_\pi(Frob_w)$ is given by $\lambda_\pi(Q_w(X))$ where*

$$Q_w(X) = X^4 - T_w X^3 + q_w(R_v + (1 + q_w^2)S_w)X^2 - q_w^3 T_w S_v X + q_w^6 S_w^2.$$

For $F = \mathbf{Q}$, if $\pi$ has multiplicity one this conjecture results from works of Laumon [23](for trivial coefficients) and Weissauer [46]. In general for $F = \mathbf{Q}$, Weissauer has proven the existence of a 4-dimensional Galois representation associated to $\pi$ when $\pi_\infty$ is holomorphic. It seems the general case where $\pi_\infty$ is in the cohomological discrete series can be treated in a similar manner, although details have not been written (oral communication of R. Weissauer).

**Definition 7.1** *For each $v \in S_p$, let $I_{F_v}$ the set of embeddings of $F_v$ in $\overline{\mathbf{Q}}_p$. Then by the choice of $\iota_p$, we can identify $I_F$ and $\bigsqcup_{v \in S_p} I_{F_v}$ by $\sigma \mapsto \iota_p \circ \sigma$. For any $\sigma \in I_F$, we set $v = v_\sigma \in S_p$ the place of $F$ dividing $p$ associated to $\iota_p \circ \sigma$. For each $v \in S_p$, we identify $G_{F_v} = Gal(\bar{F}_v/F_v)$ with a decomposition subgroup at $v$ of $Gal(\bar{F}/F)$. We denote by $\varrho_{\pi,v}$ the restriction of $\varrho_\pi$ to $G_{F_v}$.*

The conjectural local properties at places dividing $p$ are given by:

**Conjecture 3** *We keep the hypothesis of Conjecture 2. Then for all $v \in S_p$, we have*

*(i) $\varrho_{\pi,v}$ is Hodge-Tate, and for all $\sigma \in I_{F_v}$, the Hodge-Tate weights associated to $\sigma$ are:*

$$(a_\sigma + b_\sigma + c_\sigma)/2 + 3, \ (a_\sigma - b_\sigma + c_\sigma)/2 + 2, \ (-a_\sigma + b_\sigma + c_\sigma)/2 + 1,$$

$$\text{and} \quad (-a_\sigma - b_\sigma + c_\sigma)/2$$

*moreover the four corresponding Hodge numbers are equal.*

(ii) Assume that $\pi$ is unramified at $v$. Then $\varrho_{\pi,v}$ is crystalline in the sense of [10].

(iii) Assume that $\pi$ is unramified at $v$. Then the characteristic polynomial of the crystalline Frobenius acting on the filtered $\varphi$-module associated to $\varrho_{\pi,v}$ is $\lambda_\pi(Q_v(X))$.

**Comments:** (i) Assume that $F = \mathbf{Q}$, $\pi$ is unramified at $p$. If $\pi$ is endoscopic, it comes from cuspidal representations $(\sigma_1, \sigma_2)$ for $GL_2$ and its Galois representation is Hodge-Tate because those associated to $\sigma_1$ and $\sigma_2$ are. If not, then the existence of $\varrho_\pi$ is also known, and it is constructed as a submodule of the sum of four copies of the etale cohomology of the Siegel variety of level prime to $p$ (see [46] and [23]). Therefore, by the etale-crystalline comparison theorem of Faltings [9], it is crystalline at $p$ hence Hodge-Tate. The fact that the four Hodge-Tate weights occur should come from the stability of the $L$-packet at infinity (i.e. $\pi_f \otimes \pi_\infty^H$ is automorphic if and only if $\pi_f \otimes \pi_\infty^W$ is, assuming that $\pi$ is not endoscopic). This also implies that the four archimedean Hodge numbers are equal. This motivates statement (i). Statement (iii) is investigated in [42]. A proof thereof seems accessible.

Let $v \in S_p$; as a $p$-adic valuation of $F$, we normalize it by $v(p) = 1$. Let us denote by $\alpha_0, \alpha_1, \alpha_2, \alpha_3$ the $p$-adic valuations of the roots of $\lambda_\pi(Q_v(X))$ written in increasing order. Recall that $e_v$, resp. $f_v$, denotes the ramification index, resp. the residual degree of $F_v$. Then we have

**Lemma 7.1** *If $\pi$ is $Q$-ordinary, we have*

$$\alpha_0 + \alpha_3 = \alpha_1 + \alpha_2 = f_v(3+c)$$

*and*

- *If $Q_v = P$, $\alpha_0 = \dfrac{1}{2e_v} \sum_{\sigma \in I_{F_v}} (c - a_\sigma - b_\sigma)$.*

- *If $Q_v = P^*$, $\alpha_0 + \alpha_1 = f_v + \dfrac{1}{e_v} \sum_{\sigma \in I_{F_v}} (c - a_\sigma)$.*

- *If $Q_v = B$, $\alpha_0 = \dfrac{1}{2e_v} \sum_{\sigma \in I_{F_v}} (c - a_\sigma - b_\sigma)$ and $\alpha_1 = f_v + \dfrac{1}{2e_v} \sum_{\sigma \in I_{F_v}} (c + b_\sigma - a_s)$.*

**Proof:** It will follow from the calculation of the $v$-order of the coefficient of $\lambda_\pi(Q_v(X))$ adapted to $Q_v$. As explained in the beginning of Section 3.5, we consider the Hecke operator $T_{Q_v}^0$ defined in terms of the action by $\omega_v(d^{-1}) \cdot \rho_v(d))$ (for $d \in D_v$) on

$$H^{3d}(S(U); \bigotimes_{\sigma \in I_F} L_{(a_\sigma, b_\sigma; c_\sigma)}(\mathbf{C}))$$

The $\mathcal{Q}_v$-ordinarity condition says that the image of $T^0_{\mathcal{Q}_v}$ by $\lambda_\pi$ is a $v$-adic unit. One has $\omega_v(d_v) = e_v^{-1} \cdot \sum_{\sigma \in I_{F_v}} (c - a_\sigma - b_s)/2$ for $\mathcal{Q}_v = P$, resp. $\omega_v(d_v) = e_v^{-1} \cdot \sum_{\sigma \in I_{F_v}} (c - a_\sigma)$; for $\mathcal{Q}_v = B$, there two elements $d_v$, and the two previous formulas occur. For $\mathcal{Q}_v = P^*$, we thus have

$$v(\lambda_\pi(T_v)) = e_v^{-1} \cdot \sum_{\sigma \in I_{F_v}} (c - a_\sigma - b_s)/2$$

similarly, for $\mathcal{Q}_v = P^*$,

$$v(\lambda_\pi(q_v R_v) = f_v + e_v^{-1} \cdot \sum_{\sigma \in I_{F_v}} (c - a_\sigma)$$

and both for the Borel case ∎

**Definition 7.2** *Let $v \in S_p$. A weight $((a_\sigma, b_\sigma; c_\sigma))_{\sigma \in I_{F_v}}$ is called $v$-admissible if and only if the corresponding four Hodge-Tate weights given by the conjecture above are independent of $\sigma \in I_{F_v}$.*

**Proposition 7.1** *We assume the hypothesis and conclusions of conjectures 2 and 3 and that $\pi$ is Q-ordinary with separable weight (cf. definition 4.6). Then the local representations $\varrho_{\pi,v}$ takes values in a conjugate of $Q^*_v(L)$ where $Q^*_v$ is the Langlands'dual of $Q_v$ (i.e. $Q^*_v = P^*$ (resp. $P$ and $B$) if $Q_v = P$ (resp. $P^*$ and $B$).*

**Proof :** This is a consequence of Lemma 7.1 and Corollary B.1 (see the appendix B). Thereafter we give another proof in the case where $v$ is admissible because it is the special case for which the Newton and Hodge polygons meet. In the other cases, they never meet as it can be easily checked (see the last remark of the appendix B). We take the opportunity to thank H. Hida for pointing out to us that the assumption of $p$-admissibility is not satisfied in general for ordinary cuspidal automorphic forms. This led us to write the appendix B.

We use the terminology of [10]. Let $L$ be a field of coefficients of $\rho_\pi$, finite over $\mathbf{Q}_p$, of degree, say $\ell$. Choose $\sigma \in I_{F_v}$ and let us denote by $D_{cris}(\varrho_{\pi,v}) = D_{cris,\sigma}(\varrho_{\pi,v}) = (V(\varrho_{\pi,v}) \, ot B_{cris})^{G_{F_v}}$ the filtered $\varphi$-module constructed "à la Fontaine" where the action of $G_{F_v}$ on $B_{cris}$ is done through $\sigma$. It is an $L \otimes F_{v,0}$-module where $F_{v,0}$ is the maximal unramified extension of $\mathbf{Q}_p$ contained in $F_v$. By Conjecture 3.(ii), its $F_{v,0}$-dimension is equal to $4\ell$. An easy argument shows in fact that each eigenspace $D^{[t]}$ given by the slope $t$ of the crystalline Frobenius is free of rank 1 over $L \otimes F^\sigma_{v,0}$. By Conjecture 3.(i), the Hodge-Tate eigenspaces have same $F_v$-dimension hence are of dimension $\ell$. Thereafter, we construct the Hodge polygon, resp. Newton polygon by applying to the coordinates an homothety of factor $\ell^{-1}$. Let $\sigma$ any element in $I_{F_v}$.

After a Tate twist, we can assume that $c = a_\sigma + b_\sigma$. By Conjecture 3.(i) the Hodge-Tate weights are $(a_\sigma + b_\sigma + 3, a_\sigma + 2, b_\sigma + 1, 0)$. and the Hodge polygon of $D_{cris}(\varrho_{\pi,v})$ is the convex envelope of the set of points:

$$P_{Hodge} = \{(0,0), (1,0), (2, b_\sigma + 1), (3, a_\sigma + b_\sigma + 3), (4, 2a_\sigma + 2b_\sigma + 6)\}$$

Let us denote by $(t_i)_{0 \leq i \leq 3}$ the slopes of the absolute Frobenius $\varphi$. ¿From the point (iii) of the conjecture 3, they are given by divising the valuations $\alpha_i$ of the lemma 7.1 by $f_v$ (recall that the relative Frobenius is the $\varphi^{f_v}$). Since the weight is supposed to be $v$-admissible, it gives rise thus to the following Newton polygons:

$$\begin{aligned}
P^P_{Newton} &= \{(0,0), (1,0), (2, t_1), (3, a_\sigma + b_\sigma + 3), (4, 2a_\sigma + 2b_\sigma + 6)\} \\
P^{P*}_{Newton} &= \{(0,0), (1, t_0), (2, b_\sigma + 1), (3, a_\sigma + b_\sigma + t_0 + 3), (4, 2a_\sigma + 2b_\sigma + 6)\} \\
P^B_{Newton} &= \{(0,0), (1,0), (2, b_\sigma + 1), (3, a_\sigma + b_\sigma + 3), (4, 2a_\sigma + 2b_\sigma + 6)\}
\end{aligned}$$

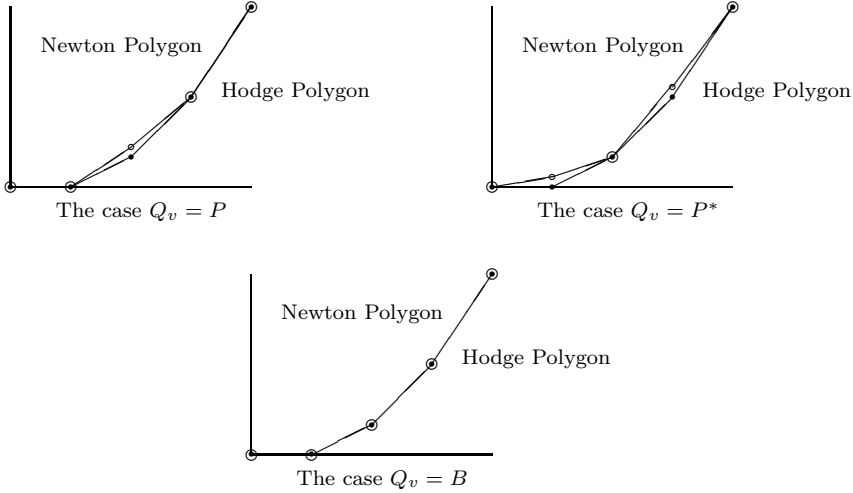

In the case, $Q_v = B$ we see immediately that $D_{cris}(\varrho_{\pi,v})$ is ordinary in the sense of [26]. In the case $Q_v = P$ (resp. $Q_v = P^*$), $Fil^{a_\sigma + b_\sigma + 3} D_{cris}(\varrho_{\pi,v})$ (resp. $Fil^{a_\sigma + 2} D_{cris}(\varrho_{\pi,v})$) is a weakly admissible submodule of $D_{cris}(\varrho_{\pi,v})$ and therefore admissible by [10] prop. 5.4.2. For $F^\sigma_{v,0}$ as above, we denote by $V_{cris} : M \mapsto (B_{cris} \otimes_{F^\sigma_{v,0}} M)_{\varphi = id}$ the quasi-inverse functor of $D_{cris}$. It is defined on the category of admissible filtered $\varphi$-module. Then,

- if $Q_v = P$, $\varrho_{\pi,v}$ leaves stable the 1-dimensional subspace

$$V_{cris}(Fil^{a_\sigma + b_\sigma + 3} D_{cris}(\varrho_{\pi,v})) \subset V(\varrho_{\pi,v}),$$

  hence takes values in $P^*(L)$,

- if $Q_v = P^*$, it leaves stable the 2-dimensional subspace

$$V_{cris}(Fil^{a_\sigma + 2} D_{cris}(\varrho_{\pi,v})) \subset V(\varrho_{\pi,v})$$

hence falls in $P(L)$,

- if $Q_v = B$, it leaves stable the three-step filtration defined by the above submodules, hence it falls in $B(L)$.

∎

Let us apply these considerations to a family of nearly $Q$-ordinary cuspidal Siegel modular forms.

**Theorem 7.1** *Let $\lambda : \mathbf{h}_\rho^{3d,S} \to \mathbf{J}$ be a family of nearly $Q$-ordinary cuspidal Siegel eigensystems such that Conjecture 2 is true for $\lambda$ mod $P$ for a densely populated set of arithmetic primes $P$. Then there exist a finite extension $\mathcal{F}$ of the fractions field $\mathcal{F}_\mathbf{J}$ of $\mathbf{J}$ and a unique semi-simple Galois representation:*

$$\varrho_\lambda : Gal(\bar{F}/F) \to GL_4(\mathcal{F})$$

*such that*

1. *It is unramified outside $S_p \cup S$ such that for all prime $v \notin S \cup S_p$ the characteristic polynomial of $\varrho_\lambda(Frob_v)$ is given by $\lambda(Q_v(X))$,*

2. *If moreover we assume Conjecture 3 for $\lambda$ mod $P$ for a densely populated set of separable algebraic primes $P$, the local representations $\varrho_{\lambda,v} = \varrho_\lambda|_{Gal(\bar{F}_v/F_v)}$ takes values in the $\mathcal{F}$-point of a parabolic subgroup of $GL_4$ whose the trace on $GSp_4$ is conjugate to $Q_v^*$ where $Q_v^*$ is the Langlands'dual of $Q_v$ (i.e. $Q_v^* = P^*$ (resp. $P$, $B$) if $Q_v = P$ (resp. $P^*$, $B$).*

**Proof :** By corollary 6.6, the existence of the representation $\varrho_\lambda$ results from the existence of the 4-dimensional representations $\varrho_{\lambda \bmod P}$ and the theory of (degree 4) pseudo-representations developed by R.Taylor (cf. [33]).

Let us prove the second statement for instance for $Q_v = P^*$, as the proof is analogous in other cases. We first establish a lemma. Let $E$ be a finite extension of $\mathbf{Q}_p$.

For any $\mathbf{J} = O_E[[T_1, \ldots, T_r]]/\mathfrak{a}$ , we put $\hat{\mathbf{J}}_\delta = \mathbf{C}_p < T_1, \ldots, T_r >_\delta /\hat{\mathfrak{a}}$

where $\mathbf{C}_p < T_1, \ldots, T_r >_\delta$ is the ring of power series converging on the closed polydisk $\mathbf{D}_p^r(\delta)$ of $\mathbf{C}_p^r$ of radius $\delta \in ]0,1[$ and $\hat{\mathfrak{a}}$ is the closure of the ideal of $\mathbf{C}_p < T_1, \ldots, T_r >_\delta$ generated by $\mathfrak{a}$. A map from $Spec(\mathbf{J})(\mathbf{C}_p) \cap \mathbf{D}_p^r(\delta)$ to $\mathbf{C}_p$ will be called analytic, if it is defined by an element of $\hat{\mathbf{J}}_\delta$.

**Lemma 7.2** *Let $E$ be a finite extension of $\mathbf{Q}_p$ and $\varrho_\mathbf{J}$ be a continuous representation of $Gal(\bar{E}/E)$ on a $\mathbf{J}$-lattice $\mathcal{L}(\varrho_\mathbf{J})$ of a $\mathcal{F}_\mathbf{J}$-vector space $V(\varrho_\mathbf{J})$. Assume there exists a densely* populated set $\mathcal{X} \subset Spec(\mathbf{J})(L)$ where $L/\mathbf{Q}_p$ is finite such that for all $P \in \mathcal{X}$, $\varrho_{\mathbf{J}}$ mod $P$ exists, is Hodge-Tate and stabilizes a subspace $V_1(P) \subset \mathcal{L}(\varrho_{\mathbf{J}}) \otimes \mathbf{J}/P$ such that the Hodge-Tate weights of $V_1(P)$ vary analytically and are different from those of $V_2(P) = (\mathcal{L}(\varrho_{\mathbf{J}}) \otimes \mathbf{J}/P)/V_1(P)$. Then there exists $V_1 \subset \mathcal{L}(\varrho_{\mathbf{J}})$ stable by $\mathrm{Gal}(\bar{E}/E)$ and such that $V_1$ mod $P \cong V_1(P)$ for all $P \in \mathcal{X}$.

**Proof :** The proof relies on results of Sen [30], [31]. Let $\mathcal{R} = End_{\mathbf{J}}(\mathcal{L}(\varrho_{\mathbf{J}}) \otimes \hat{\mathbf{J}}_\delta)$. According to [31], there exists $\varphi_{\mathbf{J}} \in \mathcal{R}$ such that for all $P \in \mathcal{X}$ $\varphi_{\mathbf{J}}$ mod $P$ is the operator defined in [30] associated to the representation $\varrho_{\mathbf{J}}$ and therefore is semisimple with integral eigenvalues corresponding to the Hodge-Tate weights. By our hypothesis, there exist elements $\mathbf{k}_1, \ldots, \mathbf{k}_t \in \hat{\mathbf{J}}_\delta$ such that $\{\mathbf{k}_1$ mod $P, \ldots, \mathbf{k}_t$ mod $P\}$ are the Hodge-Tate weights of $V_1(P)$. If we denote by $\Phi_1 = \prod_{i=1}^{t}(\varphi_{\mathbf{J}} - \mathbf{k}_i)$. Since the eigenvalues $\mathbf{k}_i(P)$ do not occur in $V_2(P)$, we have $V_1(P) \otimes \mathbf{C}_p = \mathrm{Ker}\, \Phi_1$ mod $P$. We consider $V_1(\hat{\mathbf{J}}_\delta) = \mathrm{Ker}\, \Phi_1 \subset \mathcal{L}(\varrho_{\mathbf{J}}) \otimes \hat{\mathbf{J}}_\delta$ and $V_1 = V_1(\hat{\mathbf{J}}_\delta) \cap \mathcal{L}(\varrho_{\mathbf{J}})$. Therefore for all $P \in \mathcal{X}$ such that $(V_1)_P$ and $\mathcal{L}(\varrho_{\mathbf{J}})_P$ is free over $\mathbf{J}_P$ (this is a Zariski-open condition), we have $V_1(P) = (V_1) \otimes \mathbf{J}_P/P$. Therefore $V_1 \otimes \mathbf{J}_P/P$ is stable by $\mathrm{Gal}(\bar{E}/E)$ for densely populated $P$ and thus $V_1$ is stable by $\mathrm{Gal}(\bar{E}/E)$ too.∎

Let us come back to the proof of the theorem.

Let $\varrho_{\lambda,v}$ be the restriction of $\varrho_\lambda$ to $G_{F_v}$. So, we are going to see that the hypotheses of Lemma 7.2 are satisfied by the representation $V = Ind_{\mathbf{Q}_p}^{F_v} \varrho_{\lambda,v}$. Let us denote by $\mathbf{I}$ the integral closure of $\mathbf{J}$ in $\mathcal{F}$. Let $\mathcal{L}$ be a stable $\mathbf{I}$-lattice of $V$ and $\mathcal{X} \subset Spec(\mathbf{I})(\bar{\mathbf{Q}}_p)$ be the set of sufficiently regular and separable algebraic primes (or arithmetic of level at most $p$) such that $\mathcal{L}_P$ is free over $\mathbf{I}_P$. For such $P$'s, the Hecke eigensystem $\lambda$ modulo $P$ corresponds to a cuspidal representation of level prime to $p$ by Proposition 3.2 and we can therefore apply Proposition 7.1 (in fact lemma 7.1 and corollary B.1 of the appendix): For all $P \in \mathcal{X}$, let $(a_\sigma, b_\sigma; c)_{\sigma \in I_F}$ be the algebraic and separable character associated to $\lambda$ modulo $P$. Then for $i_P = Inf_{\sigma \in I_{F_v}} a_\sigma - b_\sigma + c_\sigma)/2 + 2$

$$V_1(P) = V_{cris}(Fil^{i_P} D_{cris,\sigma}(Ind_{\mathbf{Q}_p}^{F_v} \mathcal{L}_P \otimes \mathbf{I}_P/P))$$

is a $G_{\mathbf{Q}_p}$-stable subspace of $\mathcal{L}_P \otimes \mathbf{I}_P/P$ (cf. appendix B); its Hodge-Tate weights $\{(a_\sigma + b_\sigma + c_\sigma)/2 + 3, (a_\sigma - b_\sigma + c_\sigma)/2 + 2\}_{\sigma \in I_{F_v}}$ are different from those of $V_2(P)$. One sees easily, using the series $\log_p(1+T)/log_p(1+p) \in \hat{\mathbf{J}}_\delta$, that the Hodge-Tate weights of our $V_1(P)$ are expressed analytically in terms of $P$. Thus, Lemma 7.2 applies and we note $V_1 \subset V$ the corresponding stable subspace. We can now conclude as in the end of the proof of Proposition B.1 by considering the $Res_{G_K}^{G_{\mathbf{Q}_p}} V_1 \subset Res_{G_K}^{G_{\mathbf{Q}_p}} V$.∎

**Comment:** This Theorem admits an integral version if one assumes that the reduction $\bar{\varrho}_\lambda$ of the representation $\varrho_\lambda$ is absolutely irreducible: in this situation, one can use a theorem of Nyssen [25] to construct a representation $\varrho_\lambda$ into $GSp_4(\mathbf{J})$ and even to $GSp_4(\mathbf{T})$ where $\mathbf{T}$ is the local component of $\mathbf{h}_\rho^{3d,S}$ attached to $\bar{\varrho}_\lambda$. One can then prove the analogue of Theorem 7.1. This gives rise to a $\Lambda$-algebra homomorphism from the universal ring of nearly ordinary deformations of $\bar{\varrho}_\lambda$ to $\mathbf{T}$. The possibility that it is an isomorphism for $Q = B$ was raised in [36] and more precisely in [37], and was at the origin of the present work as an analogue of [48] in the symplectic case. We hope to come back subsequently to this topic.

# A  Cohomology of the Siegel-Hilbert variety

The purpose of this appendix is to apply some results of J.Franke in order to prove some facts on the cohomology of the $3d$-dimensional Siegel-Hilbert variety. We are following very closely the formulation of [45].

Let us introduce the Notation we need. Let $\mathfrak{g}_0^1 = Lie_{\mathbf{Q}}(Res_{\mathbf{Q}}^F Sp_{4/F})$. We denote by $A_0$ be the maximal $\mathbf{Q}$-split sub-torus of the diagonal torus of $Res_{\mathbf{Q}}^F Sp_{4/F}$. Let $\mathfrak{a}_0 = Lie_{\mathbf{Q}}(A_0) \otimes_{\mathbf{Q}} \mathbf{R}$. Then $\mathfrak{a}_0^* = \mathbf{R}.\lambda_1 \oplus \mathbf{R}.\lambda_2$ where $\lambda_1(diag(x_1, x_2, -x_1, -x_2)) = x_1$ and $\lambda_2(diag(x_1, x_2, -x_1, -x_2)) = x_2$. We set $\alpha_1 = \lambda_1 - \lambda_2$ and $\alpha_2 = 2\lambda_2$.

For $P_\Sigma = P(resp. P^*, B)$ be the standard Siegel parabolic (resp. Klingen parabolic, Borel) subgroup of $G$. Let $A_\Sigma$ be the maximal $\mathbf{Q}$-split sub-torus of the center of $M_\Sigma$ the standard Levi of $P_\Sigma$. We set $\mathfrak{a}_\Sigma = Lie_{\mathbf{Q}} A_\Sigma \otimes_{\mathbf{Q}} \mathbf{R} \subset \mathfrak{a}_0$. Then

$$\mathfrak{a}_P = Ker(\alpha_1) \qquad \mathfrak{a}_{P^*} = Ker(\alpha_2) \text{ and } \mathfrak{a}_B = \mathfrak{a}_0$$

and we have

$$\mathfrak{a}_P^* = \mathbf{R}.(\lambda_1 + \lambda_2) \qquad \mathfrak{a}_{P^*}^* = \mathbf{R}.\lambda_1$$

Let $I_\infty$ the set of embedding of $F$ in $\mathbf{C}$, then

$$\mathfrak{g}^1 = Lie_{\mathbf{Q}}(Res_{\mathbf{Q}}^F Sp_{4/F}) \otimes_{\mathbf{Q}} \mathbf{C} = \oplus_{\sigma \in I_\infty} \mathfrak{sp}_4(\mathbf{C})_\sigma$$

Let $\mathfrak{h} = \oplus_{\sigma \in I_\infty} \mathfrak{h}_\sigma \subset \mathfrak{g}^1$ where $\mathfrak{h}_\sigma$ is the diagonal Cartan algebra of $sp_4(\mathbf{C})_\sigma$ we denote by $(\lambda_{1,\sigma}, \lambda_{2,\sigma})$ its canonical basis (with obvious notation). We denote also for any $\Sigma$, $\mathfrak{h}_{\Sigma,\sigma} = Lie_F(A_\Sigma \otimes F) \otimes_{F,\sigma} \mathbf{C}$.

The Weyl group $W$ is obviously isomorphic to the product $\prod_\sigma W_\sigma$ where each $W_\sigma$ is the Weyl group associated to the $\sigma$-component and is generated by $s_{\alpha_1,\sigma}$ and $s_{\alpha_2,\sigma}$. For any parabolic $P_\Sigma$, we set $W^\Sigma$ the subset of $W$ of elements $w$ such that $w^{-1}(\alpha) > 0$ for all positive roots $\alpha$ of the Levi component of $P_\Sigma$.

For any weight $\lambda \in \mathfrak{h}^*$, we denote by $\lambda_\sigma$ its canonical projection on $\mathfrak{a}_\Sigma^*$ (induced by the injection $\mathfrak{a}_\Sigma \otimes_{\mathbf{R}} \mathbf{C} \subset \mathfrak{h}$ of $Re(\lambda)$). Note that the projection $Re(\lambda)|_{\mathfrak{h}_\Sigma}$ and $\lambda_\Sigma$ are different: the first one is product over $\sigma \in I_\infty$ of the restrictions of $Re(\lambda)_\sigma$ to $\mathfrak{h}_{\Sigma,\sigma}$ while the second is a restriction of $Re(\lambda)$ to $\mathfrak{a}_\Sigma \otimes \mathbf{C} = \mathfrak{a}_0 \otimes \mathbf{C} \cap \mathfrak{h}_\Sigma$. Then if $\lambda = \sum_{\sigma \in I_\infty} x_\sigma.\lambda_{1,\sigma} + y_\sigma.\lambda_{2,\sigma}$ then

- $\lambda_B = \sum_{\sigma \in I_\infty} Re(x_\sigma)\lambda_1 + Re(y_\sigma)\lambda_2$
- $\lambda_P = \frac{1}{2} \left( \sum_{\sigma \in I_\infty} Re(x_\sigma + y_\sigma) \right).(\lambda_1 + \lambda_2)$
- $\lambda_{P^*} = \left( \sum_{\sigma \in I_\infty} Re(x_\sigma) \right).\lambda_1$

Let us denote by $\mathcal{C}$ the Weyl chamber in $\mathfrak{a}_0^*$ defined by

$$\mathcal{C} = \{(x,y) = x.\lambda_1 + y.\lambda_2 \text{ with } x \geq y \geq 0\}.$$

For $\lambda, \lambda' \in \mathfrak{a}_0^*$, we write $\lambda \succ \lambda'$ if $\lambda - \lambda' \in \mathcal{C}$. This defines a partial order; we are interested in maximal elements of finite subsets of $\mathfrak{a}^*$. For any $\lambda \in \mathfrak{a}_0^*$ we denote by $[\lambda]$ the projection of $\lambda$ on the convex $\mathcal{C}$ and for any finite subset $\Theta \subset \mathcal{C}$ we define by induction on $p \geq -1$, $\Theta^p$ the set of maximal elements of $\Theta - \cup_{i=0}^{p-1} \Theta^i$ and $\Theta^{-1} = \emptyset$. For any dominant weight $\lambda \in \mathfrak{h}^*$ and $p \in \mathbf{N}$, we set:

$$W^\Sigma(\lambda, p) = \{w \in W^\Sigma; \text{ s. t. } - w(\lambda + \rho)_\Sigma \in [W(\lambda + \rho)]^p\}.$$

$$W^\Sigma(\lambda) = \coprod_p W^\Sigma(\lambda, p)$$

Note that $W^\Sigma(\lambda)$ does not contain $id$.

Then the result of Franke we like to use in our situation is the following:

**Theorem A.1 (Franke)** *Let $K$ be an open compact subgroup of $GSp_4(\mathbb{A}_F^\infty)$ then there exist a spectral sequence whose $E_1^{p,q}$-term is given by:*

$$H_{(2)}^{p+q}(GSp_4(F)\backslash GSp_4(\mathbb{A}_F)/K.K_\mathbf{R}, E_\lambda^G) \oplus_{\Sigma \in \{B,P,P^*\}} \oplus_{\substack{w \in W^\Sigma(\lambda, p) \\ l(w) \leq p+q}}$$
$$H_{(2)}^{p+q-l(w)}(M_\Sigma(\mathbf{Q})\backslash M_\Sigma(\mathbb{A})/(K.K_\mathbf{R} \cap M_\Sigma(\mathbb{A})), E_{w(\lambda+\rho)-\rho}^{M_\Sigma}(-w(\lambda+\rho)_\Sigma))$$

*and abutting on*

$$H^{p+q}(GSp_4(F)\backslash GSp_4(\mathbb{A}_F)/K.K_\mathbf{R}, E_\lambda^G)$$

*where $H_{(2)}$ stands for $L^2$-cohomology.*

**Lemma A.1** *Let $\lambda \mathfrak{h}^*$ be a regular weight. Then for any $\Sigma$ and any $w \in W^\Sigma(\lambda)$, $w(\lambda+\rho) - \rho$ viewed as a weight of $M_\Sigma$ is regular.*

**Proof:** It is an easy calculation. Let us verify it of $P_\Sigma = P^*$. In that case, $(W^\Sigma)_\sigma = \{id, s_1, s_1 s_2, s_2 s_1 s_2\}$. We need to project for each $\sigma \in I_\infty$ $w_\sigma(\lambda_\sigma + \rho_\sigma) - \rho_\sigma$ on $\mathbf{R}.\alpha_2$ along $\mathbf{R}.\alpha_1$ and check that the projection is on the upper half-line. For $\lambda_\sigma = (x, y)$ $(x > y > 0)$ and for $w_\sigma \in (W^\Sigma)_\sigma$ of length $0, 1, 2, 3$, we find respectively $y, x+1, x+1, y$ which are strictly positive if $\lambda$ is regular.∎

**Corollary A.1** *If $\lambda$ is regular, then we have*

$$H^q(GSp_4(F)\backslash GSp_4(\mathbb{A}_F)/K.K_\mathbf{R}, E_\lambda^G) = 0$$

*for $q < 3d$.*

**Proof:** The summands of Franke's $E_1^{p,q}$ involve all Levi subgroups of $Sp_4$: first for $Sp_4$ itself, we know that the $L^2$-cohomology of $Sp_4$ with regular weight coefficients is non zero only in degree $3d$ (this results from Vogan-Zuckerman classification of unitary representations occurring in the cohomology). Then, let us consider a summand of Franke's $E_1^{p,q}$ corresponding to a type $\Sigma \in \{P, P^*, B\}$. For the Levi subgroups of the maximal parabolics, the corresponding factors vanish unless $p + q - l(w) = d$; this results from the calculation of the relative Lie algebra cohomology of $SL_2$ using the classification of cohomological unitary representations. For the Borel, the corresponding factors vanish unless $p + q - l(w) = 0$. Let us examine which $w$ does occur in the sum. For any $w \in W^\Sigma(\lambda, p)$ such that

$$H_{(2)}^{p+q-l(w)}(M_\Sigma(\mathbf{Q})\backslash M_\Sigma(\mathbb{A})/(K.K_\mathbf{R} \cap M_\Sigma(\mathbb{A})), E_{w(\lambda+\rho)-\rho}^{M_\Sigma}(-w(\lambda+\rho)_\Sigma)) \neq 0$$

the central character of $E_{w(\lambda+\rho)-\rho}^{M_\Sigma}(-w(\lambda+\rho)_\Sigma)$ needs to be trivial on the rational points of the center of $M_\Sigma$. Therefore since $F$ is totally real, the $\sigma$-components of $w(\lambda+\rho)|_{\mathfrak{h}_\Sigma}$ do not depend on $\sigma \in I_\infty$; let $\mu$ be the common value of these components. We have $w(\lambda + \rho)_\Sigma = d.Re(\mu)$. Since $w \in W^\Sigma(\lambda)$, we have $Re(\mu) < 0$, so for any $\sigma \in I_\infty$, $-w(\lambda + \rho)|_{\mathfrak{h}_{\Sigma,\sigma}} \in \mathcal{C}_\sigma$.

- For $\Sigma = P^*$, resp.$P$, we find that $w_\sigma \in \{s_1 s_2, -s_1\}$, resp. $w_\sigma \in \{s_2 s_1, -s_2\}$, so $length(w_\sigma) \geq 2$ and $l(w) \geq 2d$. On the other hand, we have $p + q = l(w) + d$; therefore we conclude $p + q \geq 3d$.

- for $\Sigma = B$, one has $-w(\lambda + \rho)_\sigma \in \mathcal{C}_\sigma$, hence $w = -id$, so $length(w_\sigma) = 4$ and $p + q = l(w) = 4d$.

This concludes the proof. ∎

# B  A remark on ordinary representations

Let $K$ be a finite extension of $\mathbf{Q}_p$ of degree $d = [K : \mathbf{Q}_p]$. Let $\rho$ be a representation of $Gal(\bar{K}/K)$ on a $E$-vector space $V$ for $E \subset \overline{\mathbf{Q}}_p$ a finite extension of $\mathbf{Q}_p$. Let $B_{cris}$ and $B_{HT}$ be the usual Fontaine's rings. We consider $\overline{\mathbf{Q}}_p$ naturally embedded in $B_{HT}$. We denote by $I_K$ the set of embeddings of $K$ in $\bar{\mathbf{Q}}_p$. Then for all $\sigma \in I_K$, we set

$$D_{HT,E,\sigma}(V) = (V \otimes_E B_{HT})^{G_K}$$

This is a $K$-vector space. We first assume that $V$ is Hodge-Tate. That means that for all $\sigma \in I_K$

$$dim_K D_{HT,E,\sigma}(V) = dim_E V = n$$

and we denote by $(a_\sigma^1, \ldots, a_\sigma^{k_\sigma})$ the corresponding Hodge-Tate weights with $a_\sigma^1 < \ldots < a_\sigma^{k_\sigma}$. We denote by with $h_\sigma^1, \ldots, h_\sigma^{k_\sigma}$ the corresponding Hodge numbers (i.e $h_\sigma^i$ is the dimension over $\mathbf{Q}_p$ of the $a_\sigma^i$-component of $D_{HT,E,\sigma}(V)$) .

Next, we assume that $V$ is cristalline. Following Fontaine, we set

$$D_{cris,\sigma}(V) = (V \otimes_{\mathbf{Q}_p} B_{cris})^{G_K}$$

where $B_{cris}$ is endowed with the action of $G_K$ via $\sigma$. Then $D_{cris,\sigma}$ is a $K_0$ vector space for $K_0$ the maximal unramified extension of $\mathbf{Q}_p$ contained in $K$. Saying that $(V,\rho)$ is cristalline means that for all $\sigma \in I_K$ we have

$$dim_{K_0} D_{cris,\sigma}(V) = dim_E(V)$$

The absolute Frobenius $\varphi$ acts on the $D_{cris,\sigma}$ and it is not difficult to see that its slopes are independent of $\sigma \in I_K$. Let us denote them by $\alpha_1 < \ldots < \alpha_l$ and let us call $d_i$ the multiplicity (which is independent of $\sigma$) of the slope $\alpha_i$ in $D_{cris,\sigma}(V)$.

**Remark:** If $V$ is cristalline, it is *a fortiori* Hodge-Tate. However note that, it does not make sense in general to consider $D_{cris,E,\sigma}$ because we cannot embed $E$ in $B_{cris}$ (e.g. if $E$ is ramified). Moreover even if $E$ is unramified this $E$ component would not be stable by the action of the absolute Frobenius $\varphi$.

In this appendix we want to overcome a difficulty mentioned in the introduction, namely that when $K \neq \mathbf{Q}_p$, it may well occur that the Hodge polygons (indexed by $\sigma \in I_K$) never meet the Newton polygon. Since our main tool for showing the ordinarity is precisely that when these polygons meet at an integral point, it yields a filtration of filtered $\varphi$-modules, this creates a problem. The idea, to overcome this, is to work with the induced representation, provided one assumes a "separability" condition which permits to relate the Hodge polygon of $V$ and that of its induction to $\mathbf{Q}_p$.

We make the following assumptions:

(Indep I) The number $k_\sigma$ of Hodge-Tate weights is independent of the embedding $\sigma$; we note it $k$.

(Indep II) The Hodge numbers are independent of the embedding $\sigma$ we denote them by $h_1, \ldots, h_k$.

These conditions are naturally satisfied if $V$ comes from the $p$-adic etale realisation of a motive. We will consider in the proposition below the following hypothesis:

(Sep $t$) For all $\sigma, \sigma' \in I_K$, we have $a_\sigma^t < a_{\sigma'}^{t+1}$.

Then we prove the following proposition:

**Proposition B.1** *We assume that for some $t$, the condition (Sep $t$) is satisfied and that*

$$\sum_{i=1}^{t} h_i \sum_{\sigma \in I_K} a_\sigma^i = [K : \mathbf{Q}_p] . \sum_{i=1}^{t} \alpha_i d_i \qquad \text{and} \qquad \sum_{i=1}^{t} h_i = \sum_{i=1}^{t} d_i$$

*Then there exists in $V$ a $E$-subspace $V'$ of dimension $\sum_{i>t} h_i$ over $\mathbf{Q}_p$ which is stable under the action of $G_K$. Moreover the Hodge-Tate weights of $V'$ associated to $\sigma \in I_K$ are $a_\sigma^{t+1}, \ldots, a_\sigma^k$ with Hodge numbers $h_{t+1}, \ldots, h_k$.*

*Proof.* Let us consider $W = \mathrm{Ind}_{G_K}^{G_{\mathbf{Q}_p}} V$ the induced representation of $G_{\mathbf{Q}_p} = \mathrm{Gal}(\overline{\mathbf{Q}}_p/\mathbf{Q}_p)$. Then the Hodge-Tate weights of $W$ are $\{a_\sigma^1, \ldots, a_\sigma^{k_\sigma}; \sigma \in I_K\}$ with Hodge numbers equal to $h_i$ for each $a_\sigma^i$, $\sigma \in I_K$ is. By the assumption (Sep $t$), we see that the point

$$([K : \mathbf{Q}_p] . \sum_{i=1}^{t} h_i, \sum_{i=1}^{t} h_i \sum_{\sigma \in I_K} a_\sigma^i)$$

is a vertex of the Hodge polygon of $W$.

The slopes of the absolute Frobenius acting on $D_{cris}(W) = (W \otimes_{\mathbf{Q}_p} B_{cris})^{G_{\mathbf{Q}_p}}$ are $\alpha_1, \ldots, \alpha_l$ with multiplicity $[K : \mathbf{Q}_p] . d_i$ for the slope $\alpha_i$. Therefore, the point

$$([K : \mathbf{Q}_p] . \sum_{i=1}^{t} d_i, [K : \mathbf{Q}_p] . \sum_{i=1}^{t} \alpha_i d_i)$$

is a vertex of the Newton polygon of $W$.

Under the assumptions of the Proposition, the two vertices described above are the same and the Newton and Hodge polygons meet at that vertex. Therefore for $i = \mathrm{Inf}_{\sigma \in I_K} a_\sigma^{t+1}$, the subspace

$$W' = V_{cris}(Fil^i D_{cris}(W)) = (B_{cris} \otimes Fil^i D_{cris}(W))^{\varphi = id}$$

is of dimension $[K : \mathbf{Q}_p](\sum_{i>t} h_i)$ and is stable under $G_{\mathbf{Q}_p}$. Its restriction to $G_K$ splits into $[K : \mathbf{Q}_p]$ subspaces permuted by the action of $Gal(\overline{\mathbf{Q}}_p/\mathbf{Q}_p)$ whose component $V'$ in $V$ (for a choosen embedding of $V$ in $Res_{G_K}^{G_{\mathbf{Q}_p}} W$) satisfies the conclusion of the proposition. The details are left to the reader.∎

The following corollary is straightforward:

**Corollary B.1** *We assume that all the Hodge numbers and all multiplicities of slopes are equal to $h$. Let $(t_i)_{1 \leq s \leq r}$ be some integers with $0 < t_1 < \ldots < t_r \leq k$ such that for all $s \in \{1, \ldots, r\}$, (Sep $t_s$) is satisfied and*

$$\sum_{i=t_s+1}^{t_{s+1}} \alpha_i = \sum_{i=t_s+1}^{t_{s+1}} \frac{1}{[K : \mathbf{Q}_p]} \sum_{\sigma \in I_K} a_\sigma^i$$

*then there exists a filtration of $E$-vector spaces $V_r \subset \ldots \subset V_0 = V$ stable under the action of $G_K$ such that for all $i$, $dim_{\mathbf{Q}_p} V_i = h.t_i$ and the Hodge-Tate weights of $V_i/V_{i+1}$ associated to $\sigma \in I_K$ are $a_\sigma^{t_i+1}, \ldots, a_\sigma^{t_{i+1}}$ with same Hodge numbers.*

**Remark** Note that this Corollary allows us to conclude that the representation $\rho$ is ordinary for the Parabolic subgroup $P$ of $GL_n$ whose Levi is $GL_{n_1} \times GL_{n_2} \times \ldots \times GL_{n_{r+1}}$ with $n_i = h(t_i - t_{i-1})/[E : \mathbf{Q}_p]$ (i.e the image of the representation falls in a conjugate of $P$), even if the Newton polygon of $\rho$ never meets the Hodge polygon associated to any embedding $\sigma \in I_K$ (that is, Newton is strictly above Hodge). This does happen if the Hodge polygons for various $\sigma$ do not coincide (that is, when the assumption of admissibility of Definition 7.2, with $F_v = K$ does not hold).